\documentclass[%
3p,%
preprint,
fleqn, 
times
]{elsarticle}

\makeindex

\usepackage{makeidx}         
\usepackage{graphicx}        
\usepackage{amsmath,amssymb}
\usepackage[babel,german=guillemets,english=american]{csquotes} 
\usepackage{dsfont} 
  \newcommand{\IP}{\ensuremath\mathds{P}}
  \newcommand{\IR}{\ensuremath\mathds{R}}
\usepackage{pstricks}
\usepackage{tikz,pgfplots}
\usepackage{tabularx}             
  \newcolumntype{R}{>{\raggedleft\arraybackslash}X}  \newcolumntype{L}{>{\raggedright\arraybackslash}X}  \newcolumntype{C}{>{\centering\arraybackslash}X}
\usepackage{booktabs} 
\usepackage{bm}    
\usepackage{mathtools}
\usepackage{fixmath}
\usepackage{hyperref}
  \hypersetup{breaklinks=true,colorlinks=true,pdfborder={0 0 0}}
  \hypersetup{linkcolor=black, anchorcolor=black, citecolor=black, filecolor=black, urlcolor=black}
\usepackage[xcolor=dvipdf, authormarkup=none]{changes} 
  \definechangesauthor[color=red,name=Vadym]{V}
  \definechangesauthor[color=red,name=Balthasar]{B}
  \definechangesauthor[color=red,name=Flo]{F}
\newcommand*{\setT}{\ensuremath{\mathcal{T}}}
\newcommand*{\setE}{\ensuremath{\mathcal{E}}}
\newcommand*{\setV}{\ensuremath{\mathcal{V}}}
\newcommand*{\vphi}{\varphi}                                     
\newcommand*{\mapEE}{\hat{\vec{\vartheta}}}                      
\renewcommand*{\vec}[1]{{\boldsymbol{#1}}}                       
\DeclareMathAlphabet{\mathbfsf}{\encodingdefault}{\sfdefault}{bx}{n}
\newcommand*{\vecc}[1]{\mathbfsf{#1}}                            
\newcommand*{\grad}{\vec{\nabla}}                                
\renewcommand*{\div}{\vec{\nabla}\cdot}                          
\newcommand*{\dd}{\mathrm{d}}                                    
\newcommand*{\abs}[1]{\ensuremath{|#1|}}                         
\newcommand*{\card}[1]{\ensuremath{\##1}}                        
\newcommand*{\transpose}[1]{{#1}^\mathrm{T}}                     
\newcommand*{\invtrans}[1]{{#1}^\mathrm{-T}}                     
\newcommand*{\jump}[1]{\left\llbracket{#1}\right\rrbracket} 
\DeclareMathOperator*{\diag}{diag}                               
\DeclareMathOperator*{\diam}{diam}                               
        
\newcommand*{\llbrace}{\left\lbrace\!\left\vert}\newcommand*{\rrbrace}{\right\vert\!\right\rbrace}
  \newcommand*{\avg}[1]{\llbrace{#1}\rrbrace}                    
\newcommand{\II}{I\!I} \newcommand{\III}{I\!I\!I} \newcommand{\IV}{I\!V} \newcommand{\V}{V} \newcommand{\VI}{V\!I} \newcommand{\VII}{V\!I\!I}

\newcommand{\Matlab}{\mbox{MATLAB}}
\newcommand{\Octave}{GNU~\mbox{Octave}}
\newcommand{\MatOct}{\Matlab\,/\,\allowbreak\Octave} 

\usepackage{textcomp} 
\usepackage{listings}
  \definecolor{keywordcolor}{rgb}{0, 0.25, 0.5}
  \definecolor{commentcolor}{rgb}{0.2, 0.5, 0.2}
  \definecolor{stringcolor}{rgb}{0.5, 0.5, 0.2}
\newcommand*{\code}[1]{\mbox{\lstinline[basicstyle=\ttfamily\small]{#1}}}

\usepackage{caption}
\begin{document}

\title{FESTUNG: A~\MatOct~toolbox for the discontinuous Galerkin method. Part I: Diffusion operator}
\author[Rice]{Florian Frank}
  \ead{florian.frank@rice.edu}
\author[FAU]{Balthasar Reuter}
  \ead{reuter@math.fau.de}
\author[FAU]{Vadym Aizinger\corref{cor}}
  \ead{aizinger@math.fau.de}
\author[FAU]{Peter Knabner}
  \ead{knabner@math.fau.de}
\address[Rice]{Rice University, Department of Computational and Applied Mathematics, 
6100 Main Street -- MS~134, Houston, TX 77005-1892, USA}
\address[FAU]{Friedrich-Alexander University of Erlangen--N\"urnberg, Department of Mathematics, 
Cauerstra{\ss}e~11, 91058~Erlangen, Germany}
\cortext[cor]{Corresponding author}

\begin{abstract}
This is the first in a~series of papers on implementing a~discontinuous Galerkin (DG) method 
as an~open source \MatOct~toolbox.
The intention of this ongoing project is to provide a rapid prototyping package for application
development using DG methods. 
The implementation relies on fully vectorized matrix\,/\,vector operations and is carefully documented;
in addition, a~direct mapping between discretization terms and code routines is maintained throughout.
The present work focuses on a~two-dimensional time-dependent diffusion equation with 
space\,/\,time-varying coefficients.  
The spatial discretization is based on the local discontinuous Galerkin formulation. 
Approximations of orders zero through four based on orthogonal polynomials have been 
implemented; more spaces of arbitrary type and order can be easily accommodated by the code structure.  
\end{abstract}

\begin{keyword}
\Matlab\sep\Octave\sep local discontinuous Galerkin method\sep vectorization\sep open source
\end{keyword}

\maketitle

 
\section{Introduction}
The discontinuous Galerkin~\mbox{(DG)} methods first introduced in 
\cite{ReedHill1973} for a~hyperbolic equation started gaining in popularity
with the appearance of techniques to deal with second order terms such as Laplace operators.  
Three different approaches to the discretization of second order terms are known in the literature.
The oldest originates from the interior penalty~\mbox{(IP)} methods introduced in the late~\mbox{1970s} 
and early~\mbox{1980s} for elliptic and parabolic equations~(cf.~\cite{ArnoldBCM2002} 
for an overview).  The IP~methods discretize the second order operators directly, similarly to the classical 
finite element method. To produce a~stable scheme, however, they need additional stabilization terms 
in the discrete formulation. 
\par
In the most recent developments, staggered DG~methods were proposed in which, in addition to element
degrees of freedom, some discontinuous vertex~\cite{LiuSTZ2007} or edge\,/\,face~\cite{ChungLee2012} 
basis functions are employed.
\par
In our \Matlab~\cite{MATLAB}\,/\,\Octave~\cite{Octave} implementation \textsl{FESTUNG} (\textsl{F}inite \textsl{E}lement 
\textsl{S}imulation \textsl{T}oolbox for \textsl{UN}structured \textsl{G}rids) available at~\cite{FESTUNG},
we rely on the local discontinuous Galer\-kin~\mbox{(LDG)} method first proposed in~\cite{CockburnShu1998} and further developed 
in~\cite{AizingerDCC2000, AizingerDawson2002}. 
The LDG~scheme utilizes a~mixed formulation in which each second order equation is replaced by
two first order equations introducing in the process an auxiliary flux variable. 
As opposed to methods from the IP~family the LDG~method is also consistent for piecewise constant approximation spaces.
\par
In developing this toolbox we pursue a~number of goals:
\begin{enumerate}
\item Design a~general-purpose software package using the DG~method for a~range of standard applications 
and provide this toolbox as a~research and learning tool in the open source format~(cf.~\cite{FESTUNG}).
\item Supply a~well-documented, intuitive user-interface to 
ease adoption by a~wider community of application and engineering professionals.
\item Relying on the vectorization capabilities of \MatOct~optimize the computational performance 
of the toolbox components and demonstrate these software development strategies.
\item Maintain throughout full compatibility with \Octave~to support users of open source software.
\end{enumerate}
The need for this kind of numerical tool appears to be very urgent right now. 
On the one hand, the DG~methods take a~significant amount of time to implement in a~computationally 
efficient manner---this hinders wider adoption of this method in the science and 
engineering community in spite of the many advantages of this type of discretization. On the other hand, 
a~number of performance optimizations, including multi-thread and GPU enhancements, combined
with a~user-friendly interface make \Matlab~and \Octave~ideal candidates for a~general purpose
toolbox simple enough to be used in students' projects but versatile enough to 
be employed by researchers and engineers to produce 'proof-of-concept' type applications 
and compute simple benchmarks.
The proposed development is by no means intended as a replacement for the traditional programming languages 
(FORTRAN, C/C++) and parallelization libraries (MPI/OpenMP) in the area of application development. 
The idea is rather to speed up the application development cycle by utilizing the rapid prototyping 
potential of \MatOct.

\subsection{Overview of existing \MatOct~DG codes}
The authors were unable to find published DG~codes running in~\Octave, 
and the number of DG~codes using \Matlab~is rather small:
A~\Matlab~code for different IP~discretizations of the one-dimensional 
Poisson~equation can be found in~\cite{Riviere2008}. 
In~\cite{LarsonBengzon2013}, an~IP~implementation for the Poisson~equation 
with homogeneous boundary conditions in two dimensions is presented.  
A~few other (unpublished) DG~\Matlab~programs can be found online, mostly small educational codes.
\par
A~special mention in this context must go to the book of Hesthaven and Warburton~\cite{HesthavenWarburton2008} 
on nodal DG~methods. A~large number of classic systems of partial differential equations in one, 
two, and three space dimensions are covered by the collection of \Matlab~codes accompanying the book.  
Most of the algorithms are time-explicit or matrix-free, but the assembly of a~full system is also presented.
The codes are available for download from~\cite{NUDG}. The book and the codes also utilize the \mbox{LDG}
framework for diffusion operators; however, the nodal basis functions used in that implementation
differ in many important ways from orthonormal modal bases adopted in the present work.
\par
None of the DG codes cited above use full vectorization in the assembly of global systems. 
A~recent preprint~\cite{Fu2013} discusses the vectorized assembly in some detail for an implementation 
of the hybridizable DG~method; however, no code has been provided 
either in the paper or as a~separate download.  
A~few \Matlab~toolboxes for the classical finite element method exist in the literature that 
focus on computationally efficient application of vectorization techniques, such as \mbox{iFEM}~\cite{Chen2009}
or \mbox{p1afem}~\cite{Funken2011}. The vectorized assembly of global matrices is 
also demonstrated in~\cite{RahmanJan2013} for the case of linear continuous elements.

\subsection{Structure of this article}
The rest of this paper is organized as follows:
We introduce the model problem in the remainder of this section and describe its
discretization using the LDG scheme in Sec.~\ref{sec:discretization}.
Implementation specific details such as data structures, reformulation and assembly of matrix
blocks, and performance studies follow in Sec.~\ref{sec:implementation}.
All routines mentioned in this work are listed and documented in Sec.~\ref{sec:routines}.
Some conclusions and an outlook of future work wrap up this publication.

\subsection{Model problem} Let~$J\coloneqq\,(0,t_\mathrm{end})\,$ be a~finite time interval and~$\Omega\subset\mathds{R}^2$ 
a~polygonally bounded domain with boundary~$\partial\Omega$ subdivided into Dirichlet~$\partial\Omega_\mathrm{D}$ 
and Neumann~$\partial\Omega_\mathrm{N}$ parts.  
We consider the \emph{diffusion equation}
\begin{subequations}\label{eq:diffusion:nonmixed}
\begin{equation}\label{eq:diffusion:nonmixed:c}
\partial_t c(t, \vec{x})  -  \div \big( d(t, \vec{x})\, \grad c(t, \vec{x}) \big) \;=\; f(t, \vec{x}) \qquad~~\text{in}~J\times\Omega
\end{equation}
with space\,/\,time-varying coefficients~$d:J\times\Omega\rightarrow \IR^+$ and $f:J\times\Omega\rightarrow \IR$.
A~prototype application of~\eqref{eq:diffusion:nonmixed:c} is the diffusive transport in fluids,
in which case the primary unknown~$c$ denotes the concentration of a~solute, 
$d$~is the diffusion coefficient, and~$f$ accounts for generation or degradation of~$c$, e.\,g., 
by chemical reactions.  
Equation~\eqref{eq:diffusion:nonmixed:c} is complemented by the following boundary and initial conditions, $\vec{\nu}$ denoting the outward unit normal:
\begin{align}
c  &\;=\; c_\mathrm{D}                         &&\text{on}~J\times{\partial\Omega}_{\mathrm{D}}\;, &&\label{eq:diffusion:nonmixed:D}\\
- \grad c\cdot\vec{\nu}   &\;=\; g_\mathrm{N}  &&\text{on}~J\times{\partial\Omega}_\mathrm{N}\;,\label{eq:diffusion:nonmixed:N}&&\\
c   &\;=\; c^0                                 &&\text{on}~\{0\}\times\Omega&&
\end{align}
\end{subequations}
with given initial $c^0:\Omega\rightarrow\IR^+_0$ and boundary data~$c_\mathrm{D}:J\times \partial\Omega_\mathrm{D}\rightarrow\IR_0^+, \;g_\mathrm{N}:J\times\partial\Omega_\mathrm{N}\rightarrow\IR$.

\section{Discretization}\label{sec:discretization}
\subsection{Notation}
Before describing the LDG~scheme for~\eqref{eq:diffusion:nonmixed} we introduce some notation; an overview can be found in the section~\enquote{Index of notation}.
\par
Let $\setT_h=\{T\}$ be a~regular family of non-overlapping partitions of~$\Omega$ into~$K$ closed 
triangles~$T$ of characteristic size~$h$ such that $\displaystyle \overline{\Omega}=\cup T$.
For~$T\in\setT_h$, let $\vec{\nu}_T$ denote the unit normal on~$\partial T$ exterior to~$T$.  
Let $\setE_\Omega$ denote the set of interior edges, $\setE_{\partial\Omega}$ the set of boundary edges, 
and $\setE \coloneqq \setE_\Omega\cup\setE_{\partial\Omega}=\lbrace E\rbrace$ 
the set of all edges (the subscript~$h$ is suppressed here). 
We subdivide further the boundary edges into Dirichlet $\setE_{\mathrm{D}}$ and Neumann $\setE_{\mathrm{N}}$ edges.
\par
For an interior edge $E\in\setE_\Omega$ shared by triangles $T^-$ and $T^+$, and for $\vec{x}\in E$, 
we define the one-sided values of a~scalar quantity~$w=w(\vec{x})$ by
\begin{equation*}
w^-(\vec{x})\;\coloneqq\;\lim_{\varepsilon \to 0^{+}} w(\vec{x} - \varepsilon\,\vec{\nu}_{T^-})
\qquad\text{and}\qquad
w^+(\vec{x})\;\coloneqq\;\lim_{\varepsilon \to 0^{+}} w(\vec{x} - \varepsilon\,\vec{\nu}_{T^+})\;,
\end{equation*}
respectively.  For a~boundary edge~$E\in\setE_{\partial\Omega}$, only the definition on the left is meaningful. 
The one-sided values of a~vector-valued quantity~$\vec{y}$ are defined analogously.  
The \emph{average} and the \emph{jump} of~$w$ on~$E$ are then given by
\begin{equation*}
\avg{w} \;\coloneqq\; (w^- + w^+)/2\qquad  \mbox{and} \qquad \jump{w} \;\coloneqq\; w^- \vec{\nu}_{T^-} + w^+ \vec{\nu}_{T^+} \;=\; (w^- - w^+)\,\vec{\nu}_{T^-}\,,
\end{equation*}
respectively. Note that $\jump{w}$ is a~vector-valued quantity.

\subsection{Mixed formulation}
To formulate an LDG~scheme we first introduce an~auxiliary vector-valued unknown~$\vec{z}\coloneqq -\grad c$ and 
re-write~\eqref{eq:diffusion:nonmixed} in mixed form, also introducing the necessary changes to the boundary conditions:
\begin{subequations}\label{eq:diffusion}
\begin{align}
\vec{z}                                            &\;=\; - \grad c                                   &&\text{in}~J\times\Omega\;,                           \label{eq:diffusion:z}\\
\partial_t c  + \div (d\,\vec{z})                  &\;=\; f                                           &&\text{in}~J\times\Omega\;,                           \label{eq:diffusion:c}\\
c                                                  &\;=\; c_\mathrm{D}                                &&\text{on}~J\times{\partial\Omega}_{\mathrm{D}}\;,    \label{eq:diffusion:D}\\
\vec{z}\cdot\vec{\nu}                              &\;=\; g_\mathrm{N}                                &&\text{on}~J\times{\partial\Omega}_\mathrm{N}\;,      \label{eq:diffusion:N}\\
c                                                  &\;=\; c^0                                         &&\text{on}~\{0\}\times\Omega\;.                       \label{eq:diffusion:0}
\end{align}
\end{subequations}

\subsection{Variational formulation}
Due to the discontinuous nature of DG~approximations, we can formulate the variational system of equations on a~triangle-by-triangle basis. 
To do that we multiply both sides of~Eqns.~\eqref{eq:diffusion:z},~\eqref{eq:diffusion:c} with 
smooth test functions~$\vec{y}:T\rightarrow\IR^2$, $w:T\rightarrow\IR$, correspondingly, and integrate by parts
over element~$T \in \setT_h$. This gives us
\begin{align*}
&\int_{T} \vec{y}\cdot\vec{z}(t)  - \int_{T} \div\vec{y} \, c(t)  
+ \int_{\partial T} \vec{y}\cdot\vec{\nu}_T\,c(t)
\;=\;0 \;,
\\
&\int_{T} w\,\partial_t c(t) - \int_{T}\grad w \cdot \big( d(t)\,\vec{z}(t)\big) + \int_{\partial T} w\,d(t)\,\vec{z}(t)\cdot\vec{\nu}_T 
\;=\; \int_{T} w\,f(t)\;.
\end{align*}

\subsection{Semi-discrete formulation}\label{sec:semidiscreteformulation}
We denote by~$\IP_p(T)$ the space of polynomials of degree at most~$p$ on~$T\in\setT_h$. 
Let
\begin{equation*}
\IP_p(\setT_h) \;\coloneqq\; \Big\{ w_h:\overline{\Omega}\rightarrow \IR\,;~\forall T\in\setT_h,~ {w_h}|_T\in\IP_p(T)\Big\}
\end{equation*}
denote the broken polynomial space on the triangulation~$\setT_h$.
For the semi-discrete formulation, we assume that the coefficient functions (for $t\in J$ fixed) are approximated as:
$d_h(t),f_h(t), c^0_h \in \IP_p(\setT_h)$. 
A~specific way to compute these approximations will be given in Sec.~\ref{sec:L2projection}; here we only
state that it is done using the $L^2$-projection into $\IP_p(T)$, 
therefore the accuracy improves with increasing polynomial order~$p$. 
Incorporating the boundary conditions~\eqref{eq:diffusion:D},~\eqref{eq:diffusion:N} and adding penalty terms for the jumps in the primary unknowns, the semi-discrete formulation reads:
\par
Seek $\left(\vec{z}_h(t), c_h(t)\right)\in  [\IP_p(\setT_h)]^2\times \IP_p(\setT_h)$ 
such that the following holds for $t\in J$ and $\forall T^-\in\setT_h,\, \forall \vec{y}_h\in [\IP_p(\setT_h)]^2,\,
\forall w_h\in \IP_p(\setT_h)\,$:
\begin{subequations}\label{prob:semidiscrete}
\begin{align}
&\hspace{-8mm}\int_{T^-} \vec{y}_h\cdot\vec{z}_h(t)
\; - \int_{T^-} \div\vec{y}_h \,c_h(t)  
\; + \int_{\partial T^-} \vec{y}_h^-\cdot\vec{\nu}_{T^-}\,\left\lbrace
\begin{array}{cl}
\avg{c_h(t) }   &\text{on}~ \setE_\Omega\\
c_\mathrm{D}(t) &\text{on}~ \setE_{\mathrm{D}}\\
c_h^-(t)        &\text{on}~ \setE_\mathrm{N}
\end{array}\right\rbrace  = 0\,,
\label{prob:semidiscrete:a}
\\
&\hspace{-8mm}\int_{T^-} w_h\,\partial_tc_h(t)
\;- \int_{T^-}\grad w_h \cdot\Big(d_h(t)\,\vec{z}_h(t) \Big)
\;+ \int_{\partial T^-} w_h^-\,
\begin{Bmatrix}
\avg{d_h(t)\,  \vec{z}_h(t)}\cdot\vec{\nu}_{T^-} + \frac{\eta}{h_{T^-}}  \jump{c_h(t)}\cdot\vec{\nu}_{T^-}& \text{on}~\setE_\Omega\\
d_h^-(t)\,  \vec{z}_h^-(t)\cdot\vec{\nu}_{T^-}  + \frac{\eta}{h_{T^-}} \big( c^-_h(t) - c_\mathrm{D}(t) \big)  & \text{on}~\setE_{\mathrm{D}}\\
d_h^-(t)\,g_\mathrm{N}(t)     & \text{on}~\setE_\mathrm{N}
\end{Bmatrix} 
\;=\; \int_{T^-} w_h\,f_h(t)\,,
\label{prob:semidiscrete:b}
\end{align}
\end{subequations}
where $\eta>0$ is a~penalty coefficient, and $h_{T^-}$ denotes the size of element $T^-$. 
The penalty terms in~\eqref{prob:semidiscrete:b} are required  to ensure 
a~full rank of the system in the absence of the time derivative~\cite[Lem.~2.15]{Riviere2008}.  
For analysis purposes, the above equations are usually summed over all triangles~$T\in\setT_h$.  
In the implementation that follows, however, it is sufficient to work with local equations. 
\par
Thus far, we used an~algebraic indexing style.
In the remainder we switch to a~mixture of algebraic and numerical style: for instance, 
$E_{kn}\in\partial T_k\cap\setE_\Omega$ means \emph{all possible} combinations of element indices~$k\in\{1,\ldots,K\}$ 
and local edge indices $n\in\{1,2,3\}$ such that $E_{kn}$ lies in~$\partial T_k\cap\setE_\Omega$. 
This implicitly fixes the numerical indices which accordingly can be used to index matrices or arrays.
\par
We use a~bracket notation followed by a~subscript to index matrices and multidimensional arrays.  
Thus, for an $n$-dimensional array~$\vecc{X}$, the symbol $[\vecc{X}]_{i_1,\ldots,i_n}$ 
stands for the component of~$\vecc{X}$ with index~$i_l$ in the $l$th dimension.  
As in \MatOct, a~colon is used to abbreviate all indices within a~single dimension.  
For example, $[\vecc{X}]_{:,:,i_3,\ldots,i_n}$ is a~two-dimensional array\,/\,matrix.

\subsubsection{Local basis representation}\label{sec:basisfunctions}
In contrast to globally continuous basis functions mostly used by the standard finite element method, 
the DG~basis functions have no continuity constraints across the triangle boundaries.  
Thus a~\emph{basis function}~$\vphi_{ki}: \overline{\Omega}\rightarrow \IR$ is only supported 
on triangle~$T_k\in\setT_h$ (i.\,e., $\vphi_{ki}=0$ on $\overline{\Omega}\smallsetminus T_k$) 
and can be defined arbitrarily while ensuring
\begin{equation*}
\forall k\in\{1,\ldots,K\}\,,\quad  \IP_p(T_k)\;=\;\mathrm{span}\,\big\{ \vphi_{ki} \big\}_{i\in\{1,\ldots,N\}}\;,
\qquad
\text{where}
\quad
N \;\coloneqq\; \frac{(p+1)(p+2)}{2} \;=\; \begin{pmatrix}p+2\\p\end{pmatrix}
\end{equation*}
is the number of \emph{local degrees of freedom}.  Clearly, the number of global degrees of freedom equals~$KN$. 
Note that~$N$ may in general vary from triangle to triangle, but we assume here for simplicity 
a~uniform polynomial degree~$p$ for every triangle. 
Closed-form expressions for basis functions on the reference triangle~$\hat{T}$ 
(cf.~Sec.~\ref{sec:transformationtoThat}) employed in our implementation up to order two are given by:
\begin{equation*}
\IP_2(\hat{T}) \left\{ 
\begin{array}{rcl}
  \begin{array}{c}
    \\ \IP_1(\hat{T}) \\ \\ 
  \end{array} & 
  \left\{\begin{array}{r}
    \IP_0(\hat{T})\; \Big\{\\ \\ \\ 
  \end{array} \right.
  &\begin{array}{rcl}
    \hat{\vphi}_1(\hat{\vec{x}})&=&\sqrt{2}\,,\\
    \hat{\vphi}_2(\hat{\vec{x}})&=& 2 - 6\hat{x}^1\,,\\
    \hat{\vphi}_3(\hat{\vec{x}})&=& 2 \sqrt{3} (1 - \hat{x}^1 - 2\hat{x}^2)\,,
  \end{array}\\
  &&\begin{array}{rcl}
    \hat{\vphi}_4(\hat{\vec{x}})&=& \sqrt{6}\big((10\hat{x}^1-8)\hat{x}^1+1\big)\,,\\
    \hat{\vphi}_5(\hat{\vec{x}})&=& 
    \sqrt{3}\big((5\hat{x}^1-4)\hat{x}^1+(-15\hat{x}^2+12)\hat{x}^2-1\big)\,,\\
    \hat{\vphi}_6(\hat{\vec{x}})&=& 
    3\sqrt{5}\big((3\hat{x}^1+8\hat{x}^2-4)\hat{x}^1+(3\hat{x}^2-4)\hat{x}^2+1\big)\,,
  \end{array}
\end{array}\right.
\end{equation*}
$\IP_p(\hat{T}) = \mathrm{span}\,\big\{\hat{\vphi}_1,\ldots,\hat{\vphi}_N\big\}$.
Note that these functions are orthonormal with respect to the \mbox{$L^2$-scalar} product on~$\hat{T}$.  
The advantage of this property will become clear in the next sections. 
The basis functions up to order four are provided in the routine~\code{phi} and their 
gradients in~\code{gradPhi}.  
Bases of even higher order can be constructed, e.\,g., with the Gram--Schmidt algorithm or 
by using a~three-term recursion relation---the latter is unfortunately not trivial to derive
in the case of triangles. 
Note that these so-called \emph{modal} basis functions~$\hat{\vphi}_i$ 
do \emph{not} posses interpolation properties at nodes unlike Lagrangian\,/\,nodal basis functions, which are often
used by the continuous finite element or nodal DG~methods.
\par
Local solutions for~$c_h$ and~$\vec{z}_h$ can be represented in terms of the local basis:
\begin{equation*}
c_h(t,\vec{x})\big|_{T_k}\eqqcolon \sum_{j=1}^N C_{kj}(t)\,\vphi_{kj}(\vec{x})\,,
\qquad
\vec{z}_h(t,\vec{x})\big|_{T_k}\eqqcolon \sum_{j=1}^N
\left(
Z_{kj}^1(t) \,\begin{bmatrix}\vphi_{kj}(\vec{x})\\ 0\end{bmatrix}
\,+\, Z_{kj}^2(t)\,\begin{bmatrix}0 \\ \vphi_{kj}(\vec{x})\end{bmatrix}\, \right)
=\sum_{j=1}^N
\begin{bmatrix}
Z_{kj}^1(t) \,\vphi_{kj}(\vec{x})\\
Z_{kj}^2(t)\,\vphi_{kj}(\vec{x})
\end{bmatrix}\,.
\end{equation*}
We condense the coefficients associated with unknowns into two-dimensional 
arrays~$\vecc{C}(t)$, $\vecc{Z}^1(t)$, $\vecc{Z}^2(t)$, such that~$C_{kj}(t) \allowbreak\coloneqq  [\vecc{C}(t)]_{k,j}$ etc. 
The vectors~$[\vecc{C}]_{k,:}$ and $[\vecc{Z}^m]_{k,:}$, $m\in\{1,2\}$,
are called \emph{local representation vectors} with respect to basis functions 
$\big\{ \vphi_{ki} \big\}_{i\in\{1,\ldots,N\}}$ 
for~$c_h$ and for components of~$\vec{z}_h$, correspondingly.
In~a~similar way, we express the coefficient functions as linear combinations of the basis functions: 
On~$T_k$, we use the local representation vectors~$[\vecc{C}^0]_{k,:}$~for 
$c_h^0$, $[\vecc{D}]_{k,:}$ for~$d_h$, and~$[\vecc{F}]_{k,:}$ for~$f_h$. 

\subsubsection{System of equations}
Testing~\eqref{prob:semidiscrete:a} with 
$\vec{y}_h = \transpose{[\vphi_{ki}, 0]}, \transpose{[0, \vphi_{ki}]}$ and~\eqref{prob:semidiscrete:b}
$w_h = \vphi_{ki}$ for $i\in\{1,\ldots,N\}$ yields a~\emph{time-dependent system of equations} 
whose contribution from~$T_k$ (identified with $T_{k^-}$ in boundary integrals) reads
\begin{subequations}\label{eq:spacediscretesystem}
\begin{align}
&\underbrace{\sum_{j=1}^N Z_{kj}^m(t) \int_{T_k}\vphi_{ki}\,\vphi_{kj}}_{I}
-\underbrace{\sum_{j=1}^N  C_{kj}(t) \int_{T_k} \partial_{x^m} \vphi_{ki}\,\vphi_{kj} }_{\II}
+ \underbrace{\int_{\partial T_k}  \vphi_{k^-i}\,\nu_{k^-}^m\,
\begin{Bmatrix}
\displaystyle \frac{1}{2}\left( \sum_{j=1}^N C_{k^- j}(t)\,\vphi_{k^-j} 
+ \sum_{j=1}^N C_{k^+j}(t)\,\vphi_{k^+j} \right) & \text{on}~\setE_\Omega \\
\displaystyle c_\mathrm{D}(t)                        & \text{on}~\setE_\mathrm{D} \\
\displaystyle \sum_{j=1}^N C_{k^- j}(t)\,\vphi_{k^-j} & \text{on}~\setE_\mathrm{N}
\end{Bmatrix}}_{\III}\nonumber\\
&\;=\;0 \quad \mbox{for}~ m\in\{1,2\}\,,
\label{eq:spacediscretesystem:a}%
\\\label{eq:spacediscretesystem:b}
&
\begin{multlined}
\underbrace{\sum_{j=1}^N\partial_t C_{kj}(t)\int_{T_k}\vphi_{ki}\,\vphi_{kj}}_{\IV}
-\underbrace{\sum_{l=1}^N D_{kl}(t) \sum_{j=1}^N\sum_{m=1}^2Z_{kj}^m(t)\int_{T_k} \partial_{x^m} \vphi_{ki}\,\vphi_{kl}\,\vphi_{kj}}_{\V}
\\
+ \underbrace{\int_{\partial T_k} \vphi_{k^-i}\,
\begin{Bmatrix}
\displaystyle \frac{1}{2}\sum_{m=1}^2\nu_{k^-}^m \left(\sum_{l=1}^N  D_{k^-l}(t)\,\vphi_{k^-l} \sum_{j=1}^N Z_{k^-j}^m(t)\,\vphi_{k^-j} + \sum_{l=1}^N  D_{k^+l}(t)\,\vphi_{k^+l} \sum_{j=1}^N Z_{k^+j}^m(t)\,\vphi_{k^+j}  \right)  & {}\\
~\hfill \displaystyle +\frac{\eta}{h_{T_{k^-}}} \left(\sum_{j=1}^N C_{k^-j}(t)\,\vphi_{k^-j} - \sum_{j=1}^N C_{k^+j}(t)\,\vphi_{k^+j}  \right) & \text{on}~\setE_\Omega \\
\displaystyle \sum_{m=1}^2\nu_{k^-}^m \sum_{l=1}^N  D_{k^-l}(t)\,\vphi_{k^-l} \sum_{j=1}^N Z_{k^-j}^m(t)\,\vphi_{k^-j}    +\frac{\eta}{h_{T_{k^-}}} \left(\sum_{j=1}^N C_{k^-j}(t)\,\vphi_{k^-j} - c_\mathrm{D}(t) \right) & \text{on}~\setE_\mathrm{D} \\
\displaystyle g_\mathrm{N}(t)\, \sum_{l=1}^N  D_{k^-l}(t)\,\vphi_{k^-l} & \text{on}~\setE_\mathrm{N}
\end{Bmatrix}
}_{\VI}
\\
\;=\;
\underbrace{\sum_{l=1}^N F_{kl}(t) \int_{T_k}\vphi_{ki}\,\vphi_{kl}}_{\VII}\;,
\end{multlined}%
\end{align}%
\end{subequations}%
where we abbreviated $\vec{\nu}_{T_k}$ by $\vec{\nu}_k = \transpose{[\nu^1_k, \nu^2_k]}$. 
Written in matrix form, system~\eqref{eq:spacediscretesystem} is then given by
\begin{equation}\label{eq:timeDepSystem}
\begin{bmatrix}
\vec{0}\\\vec{0}\\\vecc{M}\,\partial_t\vec{C}
\end{bmatrix}
+
\underbrace{\begin{bmatrix}
\vecc{M}& \cdot  & -\vecc{H}^1{+}\vecc{Q}^1{+}\vecc{Q}_{\mathrm{N}}^1\\
\cdot & \vecc{M} & -\vecc{H}^2{+}\vecc{Q}^2{+}\vecc{Q}_{\mathrm{N}}^2\\
~-\vecc{G}^1{+}\vecc{R}^1 {+}\vecc{R}^1_\mathrm{D}
& -\vecc{G}^2{+}\vecc{R}^2 {+}\vecc{R}^2_\mathrm{D}
& \eta\,\big( \vecc{S} {+} \vecc{S}_\mathrm{D}\big)
\end{bmatrix}}_{\eqqcolon \;\vecc{A}(t)}
\,
\begin{bmatrix}
\vec{Z}^1\\\vec{Z}^2\\\vec{C}
\end{bmatrix}
~~=~~ 
\begin{bmatrix}
-\vec{J}_{\mathrm{D}}^1\\
-\vec{J}_{\mathrm{D}}^2\\
\eta\, \vec{K}_{\mathrm{D}} - \vec{K}_{\mathrm{N}} {+}  \vec{L}
\end{bmatrix}
\end{equation}
with the representation vectors
\begin{align*}
\vec{Z}^m(t) &\;\coloneqq\;
\transpose{\begin{bmatrix}
Z_{11}^m(t) & \cdots & Z_{1N}^m(t) & \cdots &\cdots & Z_{K1}^m(t) & \cdots & Z_{KN}^m(t)
\end{bmatrix}}\quad \text{for}~m\in\{1,2\}\;,
\\
\vec{C}(t) &\;\coloneqq\;
\transpose{\begin{bmatrix}
C_{11}(t) & \cdots & C_{1N}(t) & \cdots &\cdots & C_{K1}(t) & \cdots & C_{KN}(t)
\end{bmatrix}}\;.
\end{align*}
The block matrices and the right-hand side of~\eqref{eq:timeDepSystem} are described in~Sections~\ref{sec:defBlocks:T} and~\ref{sec:defBlocks:E}. 
Note that some of the blocks are time-dependent (we have suppressed the time~arguments here).

\subsubsection{Contributions from area terms~$I$, $\II$, $\IV$, $\V$, $\VII$}\label{sec:defBlocks:T}
The matrices in the remainder of this section have sparse block structure; by giving definitions 
for non-zero blocks we tacitly assume a~zero fill-in.
The \emph{mass matrix}~$\vecc{M}\in\IR^{KN\times KN}$ in terms 
$I$ and $\IV$ is defined component-wise as
\begin{equation*}
[\vecc{M}]_{(k-1)N+i, (k-1)N+j} \;\coloneqq\; \int_{T_k}\vphi_{ki}\,\vphi_{kj}\;.
\end{equation*}
Since the basis functions~$\vphi_{ki}$, $i\in\{1,\ldots,N\}$ are supported only on~$T_k$, 
$\vecc{M}$ has a~block-diagonal structure 
\begin{equation}\label{eq:globMlocM}
\vecc{M} \;=\;
\begin{bmatrix}
\vecc{M}_{T_1} &          & \\
               & ~\ddots~ & \\
               &          & \vecc{M}_{T_K}
\end{bmatrix}
\qquad\text{with}\qquad
\vecc{M}_{T_k} \;\coloneqq\;
\int_{T_k}\begin{bmatrix}
\vphi_{k1}\,\vphi_{k1} & \cdots & \vphi_{k1}\,\vphi_{kN} ~\\
 \vdots                & \ddots & \vdots \\
\vphi_{kN}\,\vphi_{k1} & \cdots & \vphi_{kN}\,\vphi_{kN} 
\end{bmatrix}\;,
\end{equation}
i.\,e., it consists of $K$~\emph{local mass matrices}~$\vecc{M}_{T_k}\in\IR^{N\times N}$.  
Henceforth we write~$\vecc{M} = \diag \big(\vecc{M}_{T_1},\ldots,\vecc{M}_{T_K}\big)$.
\par
The  block matrices~$\vecc{H}^m\in\IR^{KN\times KN}, \;m\in\{1,2\}$ from term $\II$ are given by
\begin{align*}
[\vecc{H}^m]_{(k-1)N+i, (k-1)N+j} \;\coloneqq\; \int_{T_k} \partial_{x^m} \vphi_{ki}\,\vphi_{kj}\,.
\end{align*}
Hence follows the reason for placing the test function left to the solution: 
otherwise, we would be assembling the transpose of~$\vecc{H}^m$ instead.
Similarly to~$\vecc{M}$, the matrices~$\vecc{H}^m = \diag \big(\vecc{H}^m_{T_1},\ldots,\vecc{H}^m_{T_K}\big)$ are block-diagonal with local matrices
\begin{equation*}
\vecc{H}^m_{T_k} \;\coloneqq\;
\int_{T_k}\begin{bmatrix}
\partial_{x^m}\vphi_{k1}\,\vphi_{k1} & \cdots & \partial_{x^m}\vphi_{k1}\,\vphi_{kN}\\
 \vdots                & \ddots & \vdots \\
\partial_{x^m}\vphi_{kN}\,\vphi_{k1} & \cdots & \partial_{x^m}\vphi_{kN}\,\vphi_{kN} 
\end{bmatrix}\;.
\end{equation*}
In fact, all block matrices for volume integrals have block-diagonal structure 
due to the local support of the integrands.
\par
The block matrices~$\vecc{G}^m\in\IR^{KN\times KN}, \;m\in\{1,2\}$ from term~$\V$ with
\begin{align*}
[\vecc{G}^m]_{(k-1)N+i, (k-1)N+j} \;\coloneqq\; 
 \sum_{l=1}^N D_{kl}(t) \int_{T_k} \partial_{x^m} \vphi_{ki}\,\vphi_{kl} \,\vphi_{kj}
\end{align*}
are similar except that we have a~non-stationary and a~stationary factor.  
The block-diagonal matrices read~$\vecc{G}^m = \diag \big(\vecc{G}^m_{T_1},\ldots,\vecc{G}^m_{T_K}\big)$ 
with local matrices
\begin{equation}\label{eq:locGm}
\vecc{G}^m_{T_k} \;\coloneqq\;
\sum_{l=1}^N D_{kl}(t)
\int_{T_k}\begin{bmatrix}
\partial_{x^m}\vphi_{k1}\,\vphi_{kl}\,\vphi_{k1} & \cdots & \partial_{x^m}\vphi_{k1}\,\vphi_{kl}\,\vphi_{kN}\\
 \vdots                & \ddots & \vdots \\
\partial_{x^m}\vphi_{kN}\,\vphi_{kl}\,\vphi_{k1} & \cdots & \partial_{x^m}\vphi_{kN}\,\vphi_{kl}\,\vphi_{kN} 
\end{bmatrix}\;.
\end{equation}
\par
Vector~$\vec{L}(t)$ resulting from $\VII$ is obtained 
by multiplication of the representation vector of~$f_h(t)$ to the global mass matrix:
\begin{equation*}
\vec{L}(t) \;=\; \vecc{M}\, \transpose{\begin{bmatrix}
F_{11}(t) & \cdots & F_{1N}(t) & \cdots &\cdots & F_{K1}(t) & \cdots & F_{KN}(t)
\end{bmatrix}}\;.
\end{equation*}

\subsubsection{Contributions from edge terms~$\III$, $\VI$}\label{sec:defBlocks:E}


\begin{figure}[t!]\centering
\begin{pspicture}(0,-1.12)(5.02,1.12)
\rput(1.4414062,0.07){$T_{k^-}$}\rput(3.6014063,-0.21){$T_{k^+}$}\rput(3,0.39){$E_{k^-n^-}$}
\psline[linewidth=0.03cm,arrowsize=0.05291667cm 4.0,arrowlength=1.4,arrowinset=0.4]{->}(2.2,-0.18)(2.84,-0.58)
\rput(2.3914063,-0.63){$\vec{\nu}_{k^-n^-} $}
\psline[linewidth=0.03cm](3.0,1.1)(1.6,-1.1)
\psline[linewidth=0.03cm](3.0,1.1)(5.0,-0.7)
\psline[linewidth=0.03cm](5.0,-0.7)(1.6,-1.1)
\psline[linewidth=0.03cm](1.6,-1.1)(0.0,0.3)
\psline[linewidth=0.03cm](0.0,0.3)(3.0,1.1)
\end{pspicture} 
\captionsetup{justification=raggedright,singlelinecheck=false}
\caption{Two triangles adjacent to edge~$E_{k^-n^-}$. It holds~$E_{k^-n^-} = E_{k^+n^+}$ and $\vec{\nu}_{k^-n^-} = -\vec{\nu}_{k^+n^+}$.}
\label{fig:T1T2}
\end{figure}
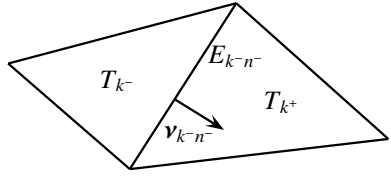
\paragraph{Interior Edges $\setE_\Omega$}
In this section, we consider a~fixed triangle $T_k=T_{k^-}$ with an~interior 
edge~$E_{k^-n^-} \in\partial T_{k^-}\cap\setE_\Omega = \partial T_{k^-} \cap \partial T_{k^+}$ 
shared by an adjacent triangle $T_{k^+}$ and associated with fixed local edge indices $n^-,n^+ \in \{1,2,3\}$~(cf.~Fig.~\ref{fig:T1T2}).
\par
First, we consider term~$\III$ in~\eqref{eq:spacediscretesystem:a}. 
For a~fixed $i\in\{1,\ldots,N\}$, 
we have a~contribution for $\vphi_{k^-i}$ in the block matrices $\vecc{Q}^m, \, m\in \{1,2\}$
\begin{equation*}
\frac{1}{2} \nu_{k^-n^-}^m \sum_{j=1}^N C_{{k^-}j}(t) \int_{E_{k^-n^-}} \vphi_{k^-i}\,\vphi_{k^-j}
+ \frac{1}{2} \nu_{k^-n^-}^m \sum_{j=1}^N C_{{k^+}j}(t) \int_{E_{k^-n^-}} \vphi_{k^-i}\,\vphi_{k^+j}\;.
\end{equation*}
Entries in diagonal blocks of $\vecc{Q}^m\in\IR^{KN\times KN}$ are then component-wise given by 
\begin{subequations}\label{eq:globQ}
\begin{equation}\label{eq:globQ:diag}
[\vecc{Q}^m]_{(k-1)N+i,(k-1)N+j} \;\coloneqq\;
\frac{1}{2} \sum_{E_{kn}\in\partial T_k\cap\setE_\Omega} \nu_{kn}^m \int_{E_{kn}} \vphi_{ki}\,\vphi_{kj}
\;=\; \frac{1}{2} \sum_{E_{kn}\in\partial T_k\cap\setE_\Omega} \nu_{kn}^m \, [\vecc{S}_{E_{kn}}]_{i,j} \;,
\end{equation}
where local matrix~$\vecc{S}_{E_{kn}}\in\IR^{N\times N}$ corresponds to interior edge~$E_{kn}$ of~$T_k, n\in\{1,2,3\}$:
\begin{equation}\label{eq:locSEkn}
\vecc{S}_{E_{kn}} = \int_{E_{kn}} \begin{bmatrix}
\vphi_{k1}\vphi_{k1} & \cdots & \vphi_{k1}\vphi_{kN} \\
\vdots & \ddots & \vdots \\
\vphi_{kN}\vphi_{k1} & \cdots & \vphi_{kN}\vphi_{kN}
\end{bmatrix}\;.
\end{equation}
Entries in off-diagonal blocks in $\vecc{Q}^m$ are only non-zero for pairs of triangles $T_{k^-}$, $T_{k^+}$ with $\partial T_{k^-}\cap\partial T_{k^+}\neq\emptyset$.
They consist of the mixed terms containing basis functions from both adjacent triangles and are given as
\begin{equation}\label{eq:globQ:offdiag}
[\vecc{Q}^m]_{({k^-}-1)N+i,({k^+}-1)N+j} \;\coloneqq\;
\frac{1}{2} \nu_{k^-n^-}^m \int_{E_{k^-n^-}} \vphi_{k^-i}\,\vphi_{k^+j}\;.
\end{equation}
\end{subequations}
Note that the local edge index $n^-$ is given implicitly since $\partial T_{k^-}\cap\partial T_{k^+}\neq\emptyset$ consist of exactly one edge~$E_{k^-n^-}=E_{k^+n^+}$.
\par
Next, consider term~$\VI$ in~\eqref{eq:spacediscretesystem:b} containing average and jump terms that produce contributions to multiple block matrices for $\vphi_{k^-i}$,
\begin{equation*}
\frac{1}{2} \nu_{k^-n^-}^m \sum_{l=1}^N D_{k^-l}(t) \sum_{j=1}^N Z_{k^-j}^m(t) 
\int_{E_{k^-n^-}} \vphi_{k^-i}\,  \vphi_{k^-l}\,  \vphi_{k^-j} + \frac{1}{2} \nu_{k^-n^-}^m \sum_{l=1}^N D_{k^+l}(t) \sum_{j=1}^N Z_{k^+j}^m(t) 
\int_{E_{k^-n^-}} \vphi_{k^-i}\, \vphi_{k^+l}\, \vphi_{k^+j}
\end{equation*}
and
\begin{equation*}
\frac{\eta}{h_{T_{k^-}}} \sum_{j=1}^N C_{k^-j}(t) \int_{E_{k^-n^-}} \vphi_{k^-i} \vphi_{k^-j} - \frac{\eta}{h_{T_{k^-}}}  \sum_{j=1}^N C_{k^+j}(t) \int_{E_{k^-n^-}} \vphi_{k^-i} \vphi_{k^+j}\;.
\end{equation*}
The first integrals are responsible for entries in the diagonal and off-diagonal blocks of a~block matrix $\vecc{R}^m\in\IR^{KN\times KN}$ that end up in the last row of system~\eqref{eq:fullSystem}. 
Entries in diagonal blocks are given component-wise by
\begin{subequations}\label{eq:globR}
\begin{equation}\label{eq:globR:diag}
[\vecc{R}^m]_{(k-1)N+i,(k-1)N+j} \;\coloneqq\;
\frac{1}{2} \sum_{E_{kn}\in\partial T_k\cap\setE_\Omega} \nu_{kn}^m \sum_{l=1}^N D_{kl}(t) 
\int_{E_{kn}} \vphi_{ki}\,\vphi_{kl}\,\vphi_{kj} \;=\; \frac{1}{2} \sum_{E_{kn}\in\partial T_k\cap\setE_\Omega} \nu_{kn}^m \, [\vecc{R}_{E_{kn}}]_{i,j}
\end{equation}
with
\begin{equation}\label{eq:locR:diag}
\vecc{R}_{E_{kn}} = \sum_{l=1}^N D_{kl}(t) \int_{E_{kn}} \begin{bmatrix}
\vphi_{k1}\vphi_{kl}\vphi_{k1} & \cdots & \vphi_{k1}\vphi_{kl}\vphi_{kN} \\
\vdots & \ddots & \vdots \\
\vphi_{kN}\vphi_{kl}\vphi_{k1} & \cdots & \vphi_{kN}\vphi_{kl}\vphi_{kN}
\end{bmatrix}\;.
\end{equation}
Once again, entries in non-zero off-diagonal blocks consist of the mixed terms:
\begin{equation}\label{eq:globR:offdiag}
[\vecc{R}^m]_{({k^-}-1)N+i,({k^+}-1)N+j} \;\coloneqq\;
\frac{1}{2} \nu_{k^-n^-}^m \sum_{l=1}^N D_{k^+l}(t) \int_{E_{k^-n^-}} \vphi_{k^-i}\,\vphi_{k^+l}\,\vphi_{k^+j}\;.
\end{equation}
\end{subequations}
All off-diagonal blocks corresponding to pairs of triangles not sharing an edge are zero. 
\par
The second integral from term~$\VI$ results in a~block matrix $\vecc{S}\in\IR^{KN\times KN}$ similar 
to $\vecc{Q}^m$. Its entries differ only in the coefficient and the lack of the normal. $h_{T_k}$~from 
the definition of the penalty term in~\eqref{prob:semidiscrete} is replaced here with the local 
edge length~$\abs{E_{kn}}$ to ensure the uniqueness of the flux over the edge 
that is necessary to ensure the local mass conservation.  
In the diagonal blocks we can reuse the previously defined $\vecc{S}_{E_{kn}}$ and have
\begin{subequations}\label{eq:globS}
\begin{equation}\label{eq:globS:diag}
[\vecc{S}]_{(k-1)N+i,(k-1)N+j} \;\coloneqq\;
 \sum_{E_{kn}\in\partial T_k\cap\setE_\Omega}\frac{1}{\abs{E_{kn}}} \int_{E_{kn}} \vphi_{ki}\vphi_{kj}\; = \;
 \sum_{E_{kn}\in\partial T_k\cap\setE_\Omega}\frac{1}{\abs{E_{kn}}}  \; [\vecc{S}_{E_{kn}}]_{i,j}\;,
\end{equation}
whereas the entries in off-diagonal blocks are given as
\begin{equation}\label{eq:globS:offdiag}
[\vecc{S}]_{({k^-}-1)N+i,({k^+}-1)N+j} \;\coloneqq\;
-\frac{1}{\abs{E_{k^-n^-}}}  \int_{E_{k^-n^-}} \vphi_{k^-i}\,\vphi_{k^+j}\;.
\end{equation}
\end{subequations}

\paragraph{Dirichlet Edges $\setE_{\mathrm{D}}$}
Consider the Dirichlet boundary~$\partial\Omega_{\mathrm{D}}$.  
The contribution of term~$\III$ of~\eqref{eq:spacediscretesystem:a} consists of a~prescribed 
data~$c_\mathrm{D}(t)$ only and consequently enters system~\eqref{eq:timeDepSystem} 
on the right-hand side as vector~$\vec{J}_\mathrm{D}^m\in\IR^{KN}$, $m\in\{1,2\}$ with
\begin{align}\label{eq:globJm:D}
[\vec{J}^m_{\mathrm{D}}]_{(k-1)N+i} &\;\coloneqq\; 
\sum_{E_{kn}\in\partial T_k\cap\setE_{\mathrm{D}}} 
\nu_{kn}^m \int_{E_{kn}} \vphi_{ki}\,c_{\mathrm{D}}(t)\;.
\end{align}
Term~$\VI$ of~\eqref{eq:spacediscretesystem:b} contains dependencies on $\vec{z}_h(t)$,~$c_h(t)$, 
and the prescribed data~$c_\mathrm{D}(t)$, thus it produces three contributions 
to system~\eqref{eq:timeDepSystem}: the left-hand side 
blocks~$\vecc{R}^m_\mathrm{D},\vecc{S}_\mathrm{D}\in\IR^{KN\times KN}$, $m\in\{1,2\}$ 
(cf.~\eqref{eq:globR:diag}, \eqref{eq:globS:diag})
\begin{align}
\label{eq:globR:D}
  [\vecc{R}^m_\mathrm{D}]_{(k-1)N+i,(k-1)N+j} & 
\;\coloneqq\;  \sum_{E_{kn}\in\partial T_k\cap\setE_{\mathrm{D}}} \nu_{kn}^m \, [\vecc{R}_{E_{kn}}]_{i,j}\;,\\
\label{eq:globS:D}
[\vecc{S}_\mathrm{D}]_{(k-1)N+i,(k-1)N+j}     & 
\;\coloneqq\;   \sum_{E_{kn}\in\partial T_k\cap\setE_{\mathrm{D}}} \frac{1}{\abs{E_{kn}}} \, [\vecc{S}_{E_{kn}}]_{i,j}\;,
\end{align}
and the right-hand side vector~$\vec{K}_\mathrm{D}\in\IR^{KN}$
\begin{equation}\label{eq:globK:D}
[\vec{K}_\mathrm{D}]_{(k-1)N+i}  
\;\coloneqq\;   \sum_{E_{kn}\in\partial T_k\cap\setE_\mathrm{D}} \frac{1}{\abs{E_{kn}}} \int_{E_{kn}} \vphi_{ki} \,c_\mathrm{D}(t)\;.
\end{equation}

\paragraph{Neumann Edges $\setE_{\mathrm{N}}$}
Consider the Neumann boundary~$\partial\Omega_{\mathrm{N}}$. 
Term~$\III$ of~\eqref{eq:spacediscretesystem:a} replaces the average of the primary variable over 
the edge by the interior value resulting in the block-diagonal 
matrix~$\vecc{Q}^m_{\mathrm{N}}\in\IR^{KN\times KN}$ (cf.~\eqref{eq:globQ:diag},~\eqref{eq:locSEkn}) with  
\begin{align}\label{eq:globQm:N}
[\vecc{Q}^m_{\mathrm{N}}]_{(k-1)N+i,(k-1)N+j} \;\coloneqq\;\sum_{E_{kn}\in\partial T_k\cap\setE_\mathrm{N}} \nu_{kn}^m \,[\vecc{S}_{E_{kn}}]_{i,j}\;.
\end{align}
Term~$\VI$ contributes to the right-hand side of system~\eqref{eq:timeDepSystem} since it contains given data only.  
The corresponding vector~$\vec{K}_\mathrm{N}\in\IR^{KN}$ reads
\begin{equation*}
[\vec{K}_\mathrm{N}]_{(k-1)N+i}  \;\coloneqq\;   
 \sum_{E_{kn}\in\partial T_k\cap\setE_\mathrm{N}} \sum_{l=1}^N D_{kl}(t)  
\int_{E_{kn}} \vphi_{ki} \,\vphi_{kl} \,g_\mathrm{N}(t)\;.
\end{equation*}

\subsection{Time discretization}
The system~\eqref{eq:timeDepSystem} is equivalent to
\begin{equation} \label{eq:fullSystem}
\vecc{W}\partial_t\vec{Y}(t) + \vecc{A}(t)\,\vec{Y}(t)  \;=\; \vec{V}(t)
\end{equation}
with $\vecc{A}(t)$ as defined in~\eqref{eq:timeDepSystem} and solution~$\vec{Y}(t) \in \IR^{3KN}$, 
right-hand-side vector~$\vec{V}(t) \in \IR^{3KN}$, and matrix~$\vecc{W}\in\IR^{3KN\times3KN}$ defined as 
\begin{equation*}
\vec{Y}(t) \;\coloneqq\;
\begin{bmatrix}\vec{Z}^1(t)\\ \vec{Z}^2(t)\\ \vec{C}(t)
\end{bmatrix}\;,
\qquad
\vec{V}(t) \;\coloneqq\;
\begin{bmatrix}
-\vec{J}_{\mathrm{D}}^1(t)\\
-\vec{J}_{\mathrm{D}}^2(t)\\
\eta \,\vec{K}_{\mathrm{D}}(t) - \vec{K}_{\mathrm{N}}(t) {+}  \vec{L}(t)
\end{bmatrix}\;,
\qquad 
\vecc{W} \;=\;
\begin{bmatrix}
\vecc{0}&\vecc{0}& \vecc{0}
\\
\vecc{0}&\vecc{0}&\vecc{0}
\\
\vecc{0}&\vecc{0}&\vecc{M}
\end{bmatrix}\;.
\end{equation*}
We discretize system~\eqref{eq:fullSystem} in time using for simplicity the implicit Euler method (generally, one 
has to note here that higher order time discretizations such as TVB (total variation bounded)
Runge--Kutta methods~\cite{CockburnShu1989} will be needed in the future
for applications to make an efficient use of high order DG~space discretizations).
Let $0= t^1<t^2<\ldots< t_\mathrm{end}$ be a~not necessarily equidistant decomposition 
of the time interval~$J$ and let $\Delta t^n \coloneqq t^{n+1} -t^n$ denote the time step size.
One step of our time discretization is formulated as
\begin{equation*} 
\left(\vecc{W} + \Delta t^n\,\vecc{A}^{n+1} \right)\,\vec{Y}^{n+1} 
\;=\; \vecc{W}\,\vec{Y}^{n} + \Delta t^n \, \vec{V}^{n+1} \,,
\end{equation*}
where we abbreviated~$\vecc{A}^n\coloneqq \vecc{A}(t^n)$, etc.

\section{Implementation}\label{sec:implementation}

We obey the following implementation conventions:
\begin{itemize}
 \item[1.] \emph{Compute every piece of information only once.} In particular, this means that stationary parts of the linear system to be solved in a~time step should be kept in the memory and not repeatedly assembled and that the evaluation of functions at quadrature points should be carried out only once.
 \item[2.] \emph{Avoid long} \code{for}~\emph{loops}. With \enquote{long} loops we mean loops that scale with the mesh size, e.\,g., loops over the triangles~$T_k\in\setT_h$ or edges~$E_k\in\setE_\Omega$.  Use \emph{vectorization} instead.
 \item[3.] \emph{Avoid changing the nonzero pattern of sparse matrices.}  Assemble global block matrices with the command \code{sparse( , , , , )}, \code{kron}, or comparable commands.
\end{itemize}
Furthermore, we try to name variables as close to the theory as possible.
Whenever we mention a~non-built-in~\MatOct~routine they are to be found in Sec.~\ref{sec:routines}.

\subsection{Grid\,/\,triangulation}\label{sec:grid}
In Sec.~\ref{sec:discretization}, we considered a~regular family of triangulations~$\setT_h$ that covers 
a~polygonally bounded~domain~$\Omega$. 
Here we fix the mesh fineness~$h$ and simply write~$\setT$ to denote the grid and also the set of triangles~$\{T\}$; the set of vertices in~$\setT$ is called~$\setV$.

\subsubsection{Data structures}
When writing \MatOct~code it is natural to use \emph{list oriented} data structures.  
Therefore, the properties of the grid~$\setT$ are stored in arrays in order to facilitate vectorization, 
in particular, by using those as index arrays.  
When we deal with a~stationary grid it is very beneficial to precompute those arrays in order to have 
access to readily usable information in the assembly routines.  
All lists describing~$\setT$ fall in two categories: \enquote{\emph{geometric data}} containing 
properties such as the coordinates of vertices~$\vec{a}_{kn}\in\setV$ or the areas of 
triangles~$T_k\in\setT$ and \enquote{\emph{topological data}} describing, e.\,g., 
the global indices~${k^+,k^-}$ of triangles sharing an edge~$E_{k^-n^-}$.  
\par
The most important lists are described in Tab.~\ref{tab:lists}.  Those and further lists are assembled by means of the routine~\code{generateGridData}, and are in some cases based on those presented in~\cite{BahriawatiCarstensen2005}.   All lists are stored in a~variable of type~\code{struct} even though it would be more efficient 
to use a~\code{class} (using \code{classdef}) instead that inherits from the class~\code{handle}. 
However, this object-oriented design strategy would go beyond the scope of this article.
\begin{table}[t!]
\begin{tabularx}{\linewidth}{@{}llX@{}}\toprule
list & dimension & description\\\midrule
\lstinline$numT$    & scalar                           & number of triangles~$\card{\setT}=K$\\
\lstinline$numE$    & scalar                           & number of edges~$\card{\setE}$\\
\lstinline$numV$    & scalar                           & number of vertices~$\card{\setV}$\\
\midrule
\lstinline$B$       & $\card{\setT}\times 2\times 2$   & transformation matrices~$\vecc{B}_k$ according to~\eqref{eq:affinemappings:Fk}\\
\lstinline$E0T$     & $\card{\setT}\times3$            & global edge indices of triangles\\
\lstinline$idE0T$   & $\card{\setT}\times3$            & edge IDs for edges~$E_{kn}$ (used to identify the interior and boundary edges as well as Dirichlet and Neumann edges)\\
\lstinline$markE0TE0T$ & $3\times 3$ (cell)            & the $(n^-,n^+)$th entry of this cell is a~sparse $K\times K$ array whose $(k^-,k^+)$th entry is one if~$E_{k^-n^-}\allowbreak=E_{k^+n^+}$\\
\lstinline$T0E$     & $\card{\setE}\times2$            & global indices of the triangles sharing edge~$E_n$ in the order dictated by the direction of the global normal on~${E_n}$ (i.\,e. $T^-, T^+$ if $\vec{\nu}_{E_n}=\vec{\nu}_{T^-}$)\\
\lstinline$V0E$     & $\card{\setE}\times2$            & global indices of the vertices sharing edge $E_n$ ordered according to the global edge orientation (the latter is given by rotating counter-clockwise by~$\pi/2$ the global edge normal to~$E_n$)\\
\lstinline$V0T$     & $\card{\setT}\times3$            & global vertex indices of triangles accounting for the counter-clockwise ordering\\
\midrule
\lstinline$areaE0T$ & $\card{\setT}\times 3$           & edge lengths~$\abs{E_{kn}}$\\
\lstinline$areaT$   & $\card{\setT}\times 1$           & triangle areas~$\abs{T_k}$\\
\lstinline$coordV0T$& $\card{\setT}\times3\times2$     & vertex coordinates~$\vec{a}_{kn}$\\
\lstinline$nuE0T$   & $\card{\setT}\times3\times2$     & local edge normals~$\vec{\nu}_{kn}$, exterior to~$T_k$\\
\bottomrule
\end{tabularx}
\caption{Arrays generated by the routine~\code{generateGridData} storing topological~\emph{(top, middle)} 
and geometric~\emph{(bottom)} grid descriptions.}
\label{tab:lists}
\end{table}

\subsubsection{Interfaces to grid generators}\label{sec:gridGeneration}
\begin{figure}[h!]
\begin{tabularx}{\linewidth}{@{}CCC@{}}
\includegraphics[width=.3\textwidth]{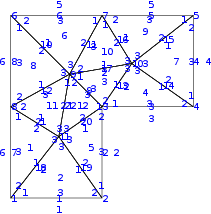} &
\includegraphics[width=.3\textwidth]{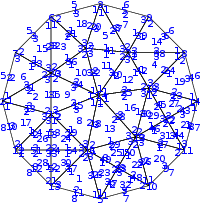} &
\includegraphics[width=.3\textwidth]{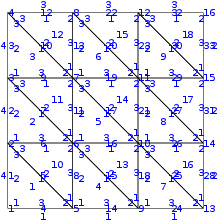} \\
\lstinline$g = domainPolygon([0 .5 .5 1 1 0], ...$ &
\lstinline$g = domainCircle(1/3)$ &
\lstinline$g = domainSquare(1/3)$ \\
\lstinline$[0 0 .5 .5 1 1], .5)$ & &
\end{tabularx}
\captionsetup{justification=raggedright,singlelinecheck=false}
\caption{Three examples for commands to build a~grid. Each is visualized
using~\code{visualizeGrid(g)}~(cf.~Sec.~\ref{sec:gridGeneration} and Sec.~\ref{sec:visualization}).}
\label{fig:gridgenerators}
\end{figure}

The routine~\code{generateGridData} requires a~list of vertex coordinates~\code{coordV} 
and an index list~\code{V0T} (cf.~Tab.~\ref{tab:lists}) to generate all further lists
for the topological and geometric description of the triangulation.  
Grid generators are a~great tool for the creation of~\code{coordV} and~\code{V0T}.
Our implementation contains at this point two interfaces to grid generators:
\begin{enumerate}
 \item The routine~\code{domainCircle} makes a~system call to the free software Gmsh~\cite{Gmsh}.
   According to the geometry description of the domain in~\code{domainCircle.geo}, Gmsh generates
   the~ASCII file~{\small\verb$domainCircle.mesh$} containing the grid, from which \code{coordV} 
   and~\code{V0T} can be extracted to call~\code{generateGridData}.
 \item \Matlab's toolbox for partial differential equations also provides a~grid generator.  
   The usage is exemplified by the routine~\code{domainPolygon} which generates 
   a~triangulation of a~polygonally bounded domain. 
\end{enumerate}
Additionally, we provide the routine~\code{domainSquare}, which produces a~Friedrichs--Keller
triangulation of a~square with given mesh size and without employing any grid generators.
An example for meshes produced by each routine is shown in Fig.~\ref{fig:gridgenerators}.

\subsection{Backtransformation to the reference triangle}\label{sec:transformationtoThat}
The computation of the volume and edge integrals in the discrete system~\eqref{eq:spacediscretesystem} is expensive when performed for each triangle~$T_k$ of the grid~$\setT_h$.  A~common practice is to transform the integrals over physical triangles~$T_k$ to a~reference triangle~$\hat{T}$ and then to compute the integrals either by numerical quadrature or analytically.  Both approaches are presented in this article.  We use the \mbox{\emph{unit}} \emph{reference triangle}~$\hat{T}$ as described in Fig.~\ref{fig:referencetriangle} and define for~$T_k\in\setT_h$ an affine \mbox{one-to-one} mapping
\begin{figure}[t!]
\centering%
\begin{pspicture}(0,-1.52)(9.32,1.52)
\psline[linecolor=black, linewidth=0.03, arrowsize=0.0529166cm 2.0,arrowlength=1.4,arrowinset=0.0]{->}(0.46,-1.02)(2.66,-1.02)
\psline[linecolor=black, linewidth=0.03, arrowsize=0.0529166cm 2.0,arrowlength=1.4,arrowinset=0.0]{->}(0.46,-1.02)(0.46,1.18)
\rput[bl](2.73,-1.32){$\hat{x}^1$}
\rput[bl](0.16,1.18){$\hat{x}^2$}

\psline[linecolor=black, linewidth=0.03](0.46,0.58)(2.06,-1.02)
\psdots[linecolor=black, dotsize=0.12](0.46,0.58)
\psdots[linecolor=black, dotsize=0.12](0.46,-1.02)
\psdots[linecolor=black, dotsize=0.12](2.06,-1.02)
\rput[bl](0.84,-0.69){$\hat{T}$}
\rput[bl](0.54,-0.92){$\hat{\vec{a}}_1$}
\rput[bl](1.09,-1.49){$\hat{E}_3$}
\rput[bl](2.0,-0.94){$\hat{\vec{a}}_2$}
\rput[bl](1.32,-0.17){$\hat{E}_1$}
\rput[bl](0.0,-0.4){$\hat{E}_2$}
\rput[bl](0.58,0.54){$\hat{\vec{a}}_3$}
\psline[linecolor=black, linewidth=0.03](0.46,-1.02)(0.46,-1.22)(0.46,-1.22)
\psline[linecolor=black, linewidth=0.03](0.46,-1.02)(0.26,-1.02)
\psline[linecolor=black, linewidth=0.03](0.46,0.58)(0.26,0.58)
\rput[bl](0.11,0.48){$1$}
\rput[bl](0.11,-1.17){$0$}
\rput[bl](2.01,-1.47){$1$}
\rput[bl](0.41,-1.47){$0$}
\psline[linecolor=black, linewidth=0.03](2.06,-1.02)(2.06,-1.22)(2.06,-1.22)

\psarc[linecolor=black, linewidth=0.022, dimen=outer, arrowsize=0.052916cm 2.0,arrowlength=1.4,arrowinset=0.4]{<-}(4.16,-2.52){2.9}{68.96249}{111.03751}
\rput[bl](4.03,0.52){$\vec{F}_k$}

\rput[bl](6.61,0.37){$E_{k2}$}
\rput[bl](7.36,1.21){{$\vec{a}_{k3}$}}
\rput[bl](7.97,0.39){$E_{k1}$}
\rput[bl](8.28,-0.51){$\vec{a}_{k2}$}
\rput[bl](7.34,-0.72){$E_{k3}$}
\rput[bl](6.49,-0.74){$\vec{a}_{k1}$}
\rput[bl](7.4,0.0){$T_k$}
\psdots[linecolor=black, dotsize=0.12](6.94,-0.36)
\psdots[linecolor=black, dotsize=0.12](8.14,-0.16)
\psdots[linecolor=black, dotsize=0.12](7.54,1.04)
\psline[linecolor=black, linewidth=0.03](7.54,1.04)(8.14,-0.16)
\psline[linecolor=black, linewidth=0.03](7.54,1.04)(6.94,-0.36)
\psline[linecolor=black, linewidth=0.03](6.94,-0.36)(8.14,-0.16)

\psline[linecolor=black, linewidth=0.03, arrowsize=0.052916cm 2.0,arrowlength=1.4,arrowinset=0.4]{->}(5.96,-1.02)(5.96,-0.32)
\psline[linecolor=black, linewidth=0.03, arrowsize=0.052916cm 2.0,arrowlength=1.4,arrowinset=0.4]{->}(5.96,-1.02)(6.66,-1.02)
\rput[bl](6.63,-1.32){$x^1$}
\rput[bl](5.66,-0.22){$x^2$}

\end{pspicture}
\caption{The affine mapping~$\vec{F}_k$ transforms the reference triangle~$\hat{T}$ with vertices~$\hat{\vec{a}}_1 = \transpose{[0,\,0]}$, $\hat{\vec{a}}_2 = \transpose{[1,\,0]}$, $\hat{\vec{a}}_3 = \transpose{[0,\,1]}$ to the physical triangle~$T_k$ with counter-clockwise-ordered vertices~$\vec{a}_{ki}$, $i\in\{1,2,3\}$.}%
\label{fig:referencetriangle}
\end{figure}
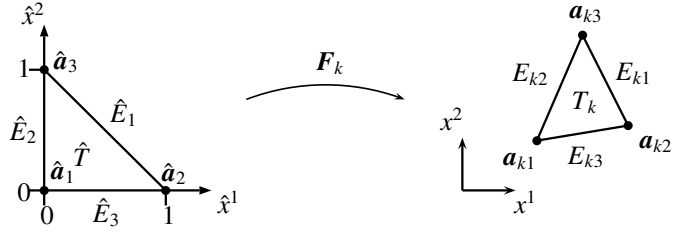
\begin{equation*}
\vec{F}_k :\quad  \hat{T} \ni \hat{\vec{x}}\mapsto \vec{x}\in T_k.
\end{equation*}
Thus any function~$\,w:T_k\rightarrow \IR\,$ implies $\hat{w}:\hat{T}\rightarrow \IR\,$ by $\,\hat{w}=w\circ \vec{F}_k\,$, i.\,e., $\,w(\vec{x}) = \hat{w}(\hat{\vec{x}})\,$.  The transformation of the gradient is obtained by the chain rule:
\begin{equation*}
\hat{\grad} \hat{w}(\hat{\vec{x}}) \;=\; \hat{\grad} w\circ\vec{F}_k(\hat{\vec{x}})  \;=\; 
\begin{bmatrix}
\partial_{x^1} w(\vec{x})\,\partial_{\hat{x}^1}F_k^1(\hat{\vec{x}}) + \partial_{x^2} w(\vec{x})\,\partial_{\hat{x}^1}F_k^2(\hat{\vec{x}}) \\
\partial_{x^1} w(\vec{x})\,\partial_{\hat{x}^2}F_k^1(\hat{\vec{x}}) + \partial_{x^2} w(\vec{x})\,\partial_{\hat{x}^2}F_k^2(\hat{\vec{x}}) 
\end{bmatrix}
\;=\;
\transpose{\big(\hat{\grad}\vec{F}_k(\hat{\vec{x}})\big)}\grad w(\vec{x})\,,
\end{equation*}
where we used the notation~$\hat{\grad} = \transpose{[\partial_{\hat{x}^1},\partial_{\hat{x}^2} ]}$, $\vec{F}_k = \transpose{[F_k^1,F_k^2]}$.  In~short,
\begin{equation}\label{eq:affinegradtrafo}
  \grad \;=\; \invtrans{\big(\hat{\grad}\vec{F}_k \big)}\,\hat{\grad}
\end{equation}
on~$T_k$.  Since~$\hat{T}$ was explicitly defined, the affine mapping can be expressed explicitly in terms of the vertices~$\vec{a}_{k1},\vec{a}_{k2},\vec{a}_{k3}$ of~$T_k$ by
\begin{subequations}\label{eq:affinemappings}

\begin{equation} \label{eq:affinemappings:Fk}
\vec{F}_k:\quad \hat{T}\ni \hat{\vec{x}} \mapsto \vecc{B}_k\,\hat{\vec{x}} +  \vec{a}_{k1}\in T_k 
\qquad \mbox{with} \qquad
\IR^{2\times2}\ni\vecc{B}_k \coloneqq \Big[  \vec{a}_{k2} - \vec{a}_{k1}\, \big|\, \vec{a}_{k3} - \vec{a}_{k1} \Big]\,.
\end{equation}
Clearly, $\hat{\grad}\vec{F}_k = \vecc{B}_k$. The inverse mapping to~$\vec{F}_k$ is easily computed:
\begin{equation} 
\vec{F}_k^{-1}:\quad T_k\ni \vec{x} \mapsto \vecc{B}_k^{-1}\,(\vec{x} - \vec{a}_{k1}) \in \hat{T}\;.
\end{equation}
\end{subequations}
Since all physical triangles have the same orientation as the reference triangle (cf.~Fig.~\ref{fig:referencetriangle}), $0<\det  \vecc{B}_k = 2\abs{T_k}$ holds.  For a~function~$w:\Omega\rightarrow\IR$, we use transformation formula
\begin{subequations}
\begin{equation}\label{eq:trafoRule:T}
\int_{T_k}w(\vec{x})\,\dd\vec{x}
\;=\;\frac{\abs{T_k}}{\abs{\hat{T}}} \int_{\hat{T}}w\circ\vec{F}_k(\hat{\vec{x}})\,\dd\hat{\vec{x}}
\;=\;2\abs{T_k} \int_{\hat{T}}w\circ\vec{F}_k(\hat{\vec{x}})\,\dd\hat{\vec{x}}
\;=\;2\abs{T_k} \int_{\hat{T}}\hat{w}(\hat{\vec{x}})\,\dd\hat{\vec{x}}
\;.
\end{equation}
The transformation rule for an integral over the edge~$E_{kn}\subset T_k$ reads
\begin{equation}\label{eq:trafoRule:E}
\int_{E_{kn}} w(\vec{x})\,\dd\vec{x} 
\;=\; \frac{\abs{E_{kn}}}{\abs{\hat{E}_n}} \int_{\hat{E}_n} w\circ\vec{F}_k(\hat{\vec{x}})\,\dd\hat{\vec{x}}
\;=\;\frac{\abs{E_{kn}}}{\abs{\hat{E}_n}} \int_{\hat{E}_n} \hat{w}(\hat{\vec{x}})\,\dd\hat{\vec{x}}
\;.
\end{equation}
The rule~\eqref{eq:trafoRule:E} is derived as follows: Denote by~$\vec{\gamma}_{kn}:[0,1]\ni s\mapsto \vec{\gamma}_{kn}(s)\in E_{kn}$ a~parametrization of the edge~$E_{kn}$ with derivative~$\vec{\gamma}_{kn}'$.  For instance $\vec{\gamma}_{k2}(s)\coloneqq (1-s)\,\vec{a}_{k3} + s\,\vec{a}_{k1}$.  Let~\mbox{$\hat{\vec{\gamma}}_{n}:[0,1]\ni s\mapsto \hat{\vec{\gamma}}_{n}(s) \in \hat{E}_{n}$} be defined analogously.  From
\begin{equation*}
\int_{E_{kn}} w(\vec{x})\,\dd\vec{x} \;=\; \int_0^1 w\circ \vec{\gamma}_{kn}(s) \,\abs{ \vec{\gamma}_{kn}'(s)}\,\dd s\;=\; \int_0^1 w\circ \vec{\gamma}_{kn}(s) \,\abs{E_{kn}}\,\dd s
\end{equation*}
and
\begin{equation*}
\int_{\hat{E}_{n}} w\circ{\vec{F}_k}(\hat{\vec{x}})\,\dd\hat{\vec{x}} \;=\; \int_0^1 w\circ{\vec{F}_k} \circ\hat{\vec{\gamma}}_{n}(s) \,\abs{ \hat{\vec{\gamma}}_{n}'(s)}\,\dd s
\;=\;\int_0^1 w\circ \vec{\gamma}_{kn}(s) \,\abs{ \hat{E}_n}\,\dd s
\end{equation*}
follows the statement in \eqref{eq:trafoRule:E}.
\end{subequations}

\subsection{Numerical integration}\label{sec:quadrature}
As an~alternative to the symbolic integration functions provided by \Matlab~we implemented a~quadrature
integration functionality for triangle and edge integrals. 
In addition, this functionality is required to produce \mbox{$L^2$-projections} 
(cf.~Sec.~\ref{sec:L2projection}) of all nonlinear functions (initial conditions, right-hand side, 
etc.) used in the system.
\par
\begin{figure}[ht!]
\centering%
\begin{tabularx}{\linewidth}{@{}CCCC@{}}
\begin{pspicture}(0,0)(2,2)
\psline[linewidth=0.03](0,0)(2,0)
\psline[linewidth=0.03](0,0)(0,2)
\psline[linewidth=0.03](0,2)(2,0)
\psdots[dotsize=0.12](0,2)
\psdots[dotsize=0.12](0,0)
\psdots[dotsize=0.12](2,0)
\psdots[dotsize=0.12,dotstyle=o,fillcolor=gray](0.3333,0.3333)
\psdots[dotsize=0.12,dotstyle=o,fillcolor=gray](1.3333,0.3333)
\psdots[dotsize=0.12,dotstyle=o,fillcolor=gray](0.3333,1.3333)
\end{pspicture}
&
\begin{pspicture}(0,0)(2,2)
\psline[linewidth=0.03](0,0)(2,0)
\psline[linewidth=0.03](0,0)(0,2)
\psline[linewidth=0.03](0,2)(2,0)
\psdots[dotsize=0.12](0,2)
\psdots[dotsize=0.12](0,0)
\psdots[dotsize=0.12](2,0)
\psdots[dotsize=0.12,dotstyle=o,fillcolor=gray](1.3327,0.3571)
\psdots[dotsize=0.12,dotstyle=o,fillcolor=gray](0.3571,1.3327)
\psdots[dotsize=0.12,dotstyle=o,fillcolor=gray](0.5600,0.1500)
\psdots[dotsize=0.12,dotstyle=o,fillcolor=gray](0.1500,0.5600)
\end{pspicture}
&
\begin{pspicture}(0,0)(2,2)
\psline[linewidth=0.03](0,0)(2,0)
\psline[linewidth=0.03](0,0)(0,2)
\psline[linewidth=0.03](0,2)(2,0)
\psdots[dotsize=0.12](0,2)
\psdots[dotsize=0.12](0,0)
\psdots[dotsize=0.12](2,0)
\psdots[dotsize=0.12,dotstyle=o,fillcolor=gray](0.8919,0.2162)
\psdots[dotsize=0.12,dotstyle=o,fillcolor=gray](0.2162,0.8919)
\psdots[dotsize=0.12,dotstyle=o,fillcolor=gray](0.8919,0.8919)
\psdots[dotsize=0.12,dotstyle=o,fillcolor=gray](0.1832,1.6337)
\psdots[dotsize=0.12,dotstyle=o,fillcolor=gray](1.6337,0.1832)
\psdots[dotsize=0.12,dotstyle=o,fillcolor=gray](0.1832,0.1832)
\end{pspicture}
&
\begin{pspicture}(0,0)(2,2)
\psline[linewidth=0.03](0,0)(2,0)
\psline[linewidth=0.03](0,0)(0,2)
\psline[linewidth=0.03](0,2)(2,0)
\psdots[dotsize=0.12](0,2)
\psdots[dotsize=0.12](0,0)
\psdots[dotsize=0.12](2,0)
\psdots[dotsize=0.12,dotstyle=o,fillcolor=gray](0.6667,0.6667)
\psdots[dotsize=0.12,dotstyle=o,fillcolor=gray](0.2026,1.5949)
\psdots[dotsize=0.12,dotstyle=o,fillcolor=gray](1.5949,0.2026)
\psdots[dotsize=0.12,dotstyle=o,fillcolor=gray](0.2026,0.2026)
\psdots[dotsize=0.12,dotstyle=o,fillcolor=gray](0.9403,0.1194)
\psdots[dotsize=0.12,dotstyle=o,fillcolor=gray](0.1194,0.9403)
\psdots[dotsize=0.12,dotstyle=o,fillcolor=gray](0.9403,0.9403)
\end{pspicture}
\\
order~$2$
&
order~$3$
&
order~$4$
&
order~$5$
\end{tabularx}
\caption{Positions of the quadrature points~$\hat{\vec{q}}_r$ on the reference triangle~$\hat{T}$ as used in the routine~\code{quadRule2D} for quadrature rules of order~$2$ to~$5$.}%
\label{fig:quadraturerules}
\end{figure}
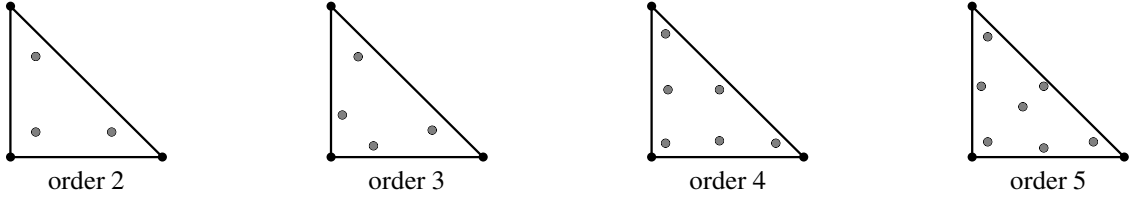

Since we transform all integrals on~$T_k\in\setT_h$ to the reference triangle~$\hat{T}$ (cf.~Sec.~\ref{sec:transformationtoThat}), it is sufficient to define the quadrature rules on~$\hat{T}$ (which, of course, can be rewritten to apply for every physical triangle~$T = \vec{F}_T(\hat{T})$):
\begin{equation}\label{eq:quadrature}
\int_{\hat{T}} \hat{g}(\hat{\vec{x}})\,\dd\hat{\vec{x}}\;\approx\; \sum_{r=1}^R \omega_r\, \hat{g}(\hat{\vec{q}}_r)
\end{equation}
with $R$~\emph{quadrature points}~$\hat{\vec{q}}_r\in\hat{T}$ and \emph{quadrature weights}~$\omega_r\in\IR$.  The \emph{order} of a~quadrature rule is the largest integer~$s$ such that~\eqref{eq:quadrature} is \emph{exact} for polynomials~$g\in\IP_s(\hat{T})$.   Note that we exclusively rely on quadrature rules with positive weights and quadrature points located strictly in the interior of~$\hat{T}$ and not on $\partial{\hat{T}}$. The rules used in the implementation are found in the routine~\code{quadRule2D}. The positions of~$\hat{\vec{q}}_r$ for some quadrature formulas are illustrated in Fig.~\ref{fig:quadraturerules}.  An~overview of quadrature rules on triangles is found in the \enquote{Encyclopaedia of Cubature Formulas}~\cite{Cools2003}. For edge integration, we rely on standard Gauss quadrature rules of required order.
\par
The integrals in~\eqref{eq:spacediscretesystem} contain integrands that are polynomials of maximum order~$3p-1$ on 
triangles and of maximum order~$3p$ on edges. Using quadrature integration one could choose rules that 
integrate all such terms exactly; however, sufficient accuracy can be achieved with quadrature rules that are exact
for polynomials of order~$2p$ on triangles and $2p+1$ on edges (cf.~\cite{CockburnShu1998b}).

\vspace{-3mm}
\subsection{Approximation of coefficient functions and initial conditions}\label{sec:L2projection}
In~Sec.~\ref{sec:semidiscreteformulation}, we assumed the coefficient functions and initial conditions 
given in piecewise polynomial spaces, for instance, $d_h(t)\in\IP_d(\setT_h)$ for $t\in J$. 
If we have an~algebraic expression for a~coefficient, say~$d$, 
we seek the representation matrix~$\vecc{D}(t)\in\IR^{K\times N}$ satisfying
\begin{equation*}
d_h(t,\vec{x})\big|_{T_k} \;=\; \sum_{j=1}^N D_{kj}(t)\,\vphi_{kj}(\vec{x})\,,
\end{equation*}
such that~$d_h(t)$ is an~adequate approximation of~$d(t)$.  
A~simple way (also used in this work) to produce $d_h$ is the \mbox{\emph{$L^2$-projection}} defined locally for $T_k\in\setT_h$ by 
\begin{equation*}
\forall w_h\in \IP_d(T)\,,\quad \int_{T_k} w_h \,d_h(t) \;=\;\int_{T_k} w_h\,d(t)\;.
\end{equation*}
Choosing $w_h = \vphi_{ki}$ for $i\in\{1\,\ldots,N\}$ and using the affine mapping~$\vec{F}_k$ we obtain
\begin{align*}
\sum_{j=1}^N D_{kj}(t) \int_{T_k} \vphi_{ki}(\vec{x})\,\vphi_{kj}(\vec{x})\,\dd\vec{x} = \int_{T_k} \vphi_{ki}(\vec{x})\,d(t,\vec{x})\,\dd\vec{x}
\Leftrightarrow    \sum_{j=1}^N D_{kj}(t) \int_{\hat{T}} \hat{\vphi}_{i}(\vec{\hat{x}})\,\hat{\vphi}_{j}(\vec{\hat{x}})\,\dd\vec{\hat{x}}
 =  \int_{\hat{T}} \hat{\vphi}_{i}(\vec{\hat{x}})\,d\big(t,\vec{F}_k(\hat{\vec{x}})\big)\,\dd\vec{\hat{x}}\;,
\end{align*}
where the factor of~$2\abs{T_k}$ canceled out.  Written in matrix form, this is equivalent to
\begin{equation*}
\hat{\vecc{M}}\,
\begin{bmatrix}
D_{k1}\\
\vdots\\
D_{kN}
\end{bmatrix}
\;=\;
\int_{\hat{T}}
\begin{bmatrix}
\hat{\vphi}_{1}(\vec{\hat{x}})\,d\big(t,\vec{F}_k(\hat{\vec{x}})\big)\\
\vdots\\
\hat{\vphi}_{N}(\vec{\hat{x}})\,d\big(t,\vec{F}_k(\hat{\vec{x}})\big)
\end{bmatrix}
\dd\vec{\hat{x}}
\end{equation*}
with local mass matrix on the reference triangle~$\hat{\vecc{M}}\in\IR^{N\times N}$ defined as in~\eqref{eq:hatM}.  
This $N\times N$ system of equations can be solved locally for every~$k\in\{1,\ldots,K\}$.  
Approximating the right-hand side by numerical quadrature~\eqref{eq:quadrature} and transposing  the equation yields
\begin{equation*}
\vecc{D}(t)~\hat{\vecc{M}}=\sum_{r=1}^R \omega_r\,
\begin{bmatrix}
d\big(t,\vec{F}_1(\hat{\vec{q}}_r)\big)\\
\vdots\\
d\big(t,\vec{F}_K(\hat{\vec{q}}_r)\big)
\end{bmatrix}\;
\arraycolsep=2.2pt
\begin{bmatrix}
\hat{\vphi}_{1}(\vec{\hat{q}}_r), &\hdots, &\hat{\vphi}_{N}(\vec{\hat{q}}_r) 
\end{bmatrix}
=
\begin{bmatrix}
d\big(t,\vec{F}_1(\hat{\vec{q}}_1)\big) & \hdots & d\big(t,\vec{F}_1(\hat{\vec{q}}_R)\big)\\
\vdots & \ddots & \vdots \\
d\big(t,\vec{F}_K(\hat{\vec{q}}_1)\big) & \hdots & d\big(t,\vec{F}_K(\hat{\vec{q}}_R)\big)
\end{bmatrix}\;
\begin{bmatrix}
\omega_1\,\hat{\vphi}_1(\vec{\hat{q}}_1) & \hdots & \omega_1\,\hat{\vphi}_N(\vec{\hat{q}}_1)\\
\vdots & \ddots & \vdots \\
\omega_R\,\hat{\vphi}_1(\vec{\hat{q}}_R) & \hdots & \omega_R\,\hat{\vphi}_N(\vec{\hat{q}}_R) 
\end{bmatrix}\;.
\end{equation*}
This is the global matrix-valued~(transposed) system of equations with unknown~$\vecc{D}(t)\in\IR^{K\times N}$ 
and a~right-hand side of dimension~$K\times N$.
The corresponding routine is~\code{projectFuncCont2DataDisc}.

\subsection{Computation of the discretization error}\label{sec:discerror}
The \emph{discretization error}~$\|c_h(t)-c(t)\|_{L^2(\Omega)}$ at time~$t\in J$ gives the $L^2$-norm of the 
difference between the discrete solution~$c_h(t)$ and the analytical solution~$c(t)$ with the latter usually specified as an~algebraic function.  
This computation is utilized in the computation of the experimental rate of convergence of the numerical scheme 
(cf.~Sec.~\ref{sec:codeverification}).  
\par
As~in the previous section we use here the numerical quadrature after transforming the integral term to the 
reference triangle~$\hat{T}$.  
The arising sums are vectorized for reasons of performance.  Suppressing the time argument, we have
\begin{align*}
\|c_h-c\|_{L^2(\Omega)}^2 
\;=\; \sum_{T_k\in\setT_h}\int_{T_k} \Big(c_h(\vec{x})-c(\vec{x})\Big)^2\,\dd\vec{x}
&=&& \hspace{-2mm} 
2 \sum_{T_k\in\setT_h} \abs{T_k} \int_{\hat{T}} \Big(\sum_{l=1}^N C_{kl}\,\hat{\vphi}_l(\hat{\vec{x}}) - c\circ\vec{F}_k(\hat{\vec{x}})\Big)^2\,\dd\hat{\vec{x}}\\
\hspace{-8mm} 
\approx 2 \sum_{T_k\in\setT_h} \abs{T_k} \sum_{r=1}^R \omega_r \Big(\sum_{l=1}^N C_{kl}\,\hat{\vphi}_l(\hat{\vec{q}}_r) - c\circ\vec{F}_k(\hat{\vec{q}}_r)\Big)^2 
&=&&\hspace{-2mm}
\arraycolsep=1.2pt
2 \begin{bmatrix}
\abs{T_1}\\\vdots\\\abs{T_K}
\end{bmatrix}
\cdot \left(
\begin{bmatrix}
C_{11}&\hdots&C_{1N}\\
\vdots&{}&\vdots\\
C_{K1}&\hdots&C_{KN}
\end{bmatrix}
\begin{bmatrix}
\hat{\vphi}_1(\hat{\vec{q}}_1) & \hdots & \hat{\vphi}_1(\hat{\vec{q}}_R)\\
 \vdots&{}&\vdots\\
\hat{\vphi}_N(\hat{\vec{q}}_1) & \hdots & \hat{\vphi}_N(\hat{\vec{q}}_R)
\end{bmatrix} - c \left(\vecc{X}^1, \vecc{X}^2 \right)\right)^2 
\begin{bmatrix}
\omega_1\\\vdots\\\omega_R
\end{bmatrix}\;,
\end{align*}
where the arguments of~$c$, $[\vecc{X}^m]_{k,r}\coloneqq F^m_k(\hat{\vec{q}}_r)$, $k\in\{1,\ldots,K\}$, $r\in\{1,\ldots,R\}$,
can be assembled using a~Kronecker product. Somewhat abusing notation we mean by~$c(\vecc{X}^1, \vecc{X}^2 )$ the $K\times R$~matrix with the entry $c([\vecc{X}^1]_{k,r}, [\vecc{X}^2]_{k,r})$ in the~$k$th row and~$r$th column.
The above procedure is implemented in the routine~\code{computeL2Error}.

\subsection{Assembly}
The aim of this section is to transform the terms required to build the block matrices in~\eqref{eq:timeDepSystem} to the reference triangle~$\hat{T}$ and then to compute those either via numerical quadrature or analytically.  
The assembly of the block matrices from the local contributions is then performed in vectorized operations.
\par
For the implementation, we need the explicit form for the components of the mappings~$\vec{F}_k:\hat{T}\rightarrow T_k$ and their inverses~$\vec{F}_k^{-1}: T_k\rightarrow \hat{T}$ as defined in~\eqref{eq:affinemappings}.
Recalling that $0<\det  \vecc{B}_k = 2\abs{T_k}$ (cf.~Sec.~\ref{sec:transformationtoThat}) we obtain
\begin{equation*}
\vec{F}_k(\hat{\vec{x}}) \;=\;
\begin{bmatrix}
B_k^{11}\,\hat{x}^1 + B_k^{12}\,\hat{x}^2 + a_{k1}^1\\
B_k^{21}\,\hat{x}^1 + B_k^{22}\,\hat{x}^2 + a_{k1}^2
\end{bmatrix}
\qquad\text{and}\qquad
\vec{F}_k^{-1}(\vec{x}) \;=\;
\frac{1}{2\,\abs{T_k}}
\begin{bmatrix}
B_k^{22}\,(x^1 - a_{k1}^1) - B_k^{12}\,(x^2 - a_{k1}^2)\\
B_k^{11}\,(x^2 - a_{k1}^2) - B_k^{21}\,(x^1 - a_{k1}^1)
\end{bmatrix}\;.
\end{equation*}
From \eqref{eq:affinegradtrafo} we obtain the component-wise rule for the gradient in~$\vec{x}\in T_k$:
\begin{equation}\label{eq:rule:gradient}
\begin{bmatrix}
\partial_{x^1}\\
\partial_{x^2}
\end{bmatrix}
\;=\;
\frac{1}{2\,\abs{T_k}}
\begin{bmatrix}
B_k^{22}\,\partial_{\hat{x}^1} - B_k^{21}\,\partial_{\hat{x}^2}\\
B_k^{11}\,\partial_{\hat{x}^2} - B_k^{12}\,\partial_{\hat{x}^1}
\end{bmatrix}.
\end{equation}
\par
In the following, we present the necessary transformation for all blocks of system~\eqref{eq:timeDepSystem} and name the corresponding \MatOct~routines that can be found in Sec.~\ref{sec:routines}.

\subsubsection{Assembly of $\vecc{M}$}\label{sec:assembly:globM}
Using the transformation rule~\eqref{eq:trafoRule:T} the following holds for the local mass matrix~$\vecc{M}_{T_k}$ as defined in~\eqref{eq:globMlocM}:
\begin{equation}\label{eq:hatM}
\vecc{M}_{T_k} \;=\; 2\abs{T_k}\,\hat{\vecc{M}}
\qquad\text{with}\qquad
\hat{\vecc{M}}\;\coloneqq\;
\int_{\hat{T}}\,\begin{bmatrix}
\hat{\vphi}_{1}\,\hat{\vphi}_{1} & \cdots & \hat{\vphi}_{1}\,\hat{\vphi}_{N} ~\\
 \vdots                & \ddots & \vdots \\
\hat{\vphi}_{N}\,\hat{\vphi}_{1} & \cdots & \hat{\vphi}_{N}\,\hat{\vphi}_{N} 
\end{bmatrix}\;,
\end{equation}
where~$\hat{\vecc{M}}\in\IR^{N\times N}$ is the representation of the local mass matrix on the reference triangle~$\hat{T}$.  
With~\eqref{eq:globMlocM} we see that the global mass matrix~$\vecc{M}$ can be expressed as a~Kronecker product of a~matrix containing the areas~$\abs{T_k}$ and the local matrix~$\hat{\vecc{M}}$:
\begin{equation*}
\vecc{M}
\;=\;
\begin{bmatrix}
\vecc{M}_{T_1} &          & \\
               & ~\ddots~ & \\
               &          & \vecc{M}_{T_K}
\end{bmatrix}
\;=\;
2 \begin{bmatrix}
\abs{T_1} &          & \\
               & ~\ddots~ & \\
               &          & \abs{T_K}
\end{bmatrix} \otimes \hat{\vecc{M}}\;.
\end{equation*}
In the corresponding assembly routine~\code{assembleMatElemPhiPhi}, the sparse block-diagonal matrix is generated using the command~\code{spdiags} with the list~\code{g.areaT} (cf.~Tab.~\ref{tab:lists}).

\subsubsection{Assembly of $\vecc{H}^m$}\label{sec:assembly:globH}
The transformation rules \eqref{eq:trafoRule:T} and \eqref{eq:rule:gradient} yield 
\begin{equation*}
\vecc{H}_{T_k}^1\;=\; B_k^{22}\,[\hat{\vecc{H}}]_{:,:,1} - B_k^{21}\,[\hat{\vecc{H}}]_{:,:,2}\qquad \text{and}\qquad
\vecc{H}_{T_k}^2\;=\; - B_k^{12}\,[\hat{\vecc{H}}]_{:,:,1} + B_k^{11}\,[\hat{\vecc{H}}]_{:,:,2} 
\end{equation*}
with
\begin{equation}\label{eq:hatH}
[\hat{\vecc{H}}]_{:,:,m}\;\coloneqq\;
\int_{\hat{T}}\,\begin{bmatrix}
\partial_{\hat{x}^m}\hat{\vphi}_{1}\,\hat{\vphi}_{1} & \cdots & \partial_{\hat{x}^m}\hat{\vphi}_{1}\,\hat{\vphi}_{N} ~\\
 \vdots                & \ddots & \vdots \\
\partial_{\hat{x}^m}\hat{\vphi}_{N}\,\hat{\vphi}_{1} & \cdots & \partial_{\hat{x}^m}\hat{\vphi}_{N}\,\hat{\vphi}_{N} 
\end{bmatrix}\in\IR^{N\times N} \quad \mbox{for }m\in\{1,2\}.
\end{equation}
Similar to~$\vecc{M}$, the global matrices~$\vecc{H}^m$ are assembled by Kronecker products in the routine~\code{assembleMatElemDphiPhi}.

\subsubsection{Assembly of $\vecc{G}^m$}\label{sec:assembly:globG}
Application of the product rule, \eqref{eq:trafoRule:T}, and \eqref{eq:rule:gradient} give us
\begin{equation*}
\int_{T_k}  \partial_{x^1}\vphi_{ki}\,\vphi_{kj} \;=\;
\phantom{-}B_k^{22}\,[\hat{\vecc{G}}]_{i,j,l,1} - B_k^{21}\,[\hat{\vecc{G}}]_{i,j,l,2}\;,
\qquad
\int_{T_k}  \partial_{x^2}\vphi_{ki}\,\vphi_{kj} \;=\;
-B_k^{12}\,[\hat{\vecc{G}}]_{i,j,l,1} + B_k^{11}\,[\hat{\vecc{G}}]_{i,j,l,2}
\end{equation*}
with a~multidimensional array~$\hat{\vecc{G}}\in\IR^{N\times N\times N\times 2}$ representing the transformed integral on the reference triangle~$\hat{T}$:
\begin{equation}\label{eq:hatG}
[\hat{\vecc{G}}]_{i,j,l,m}\;\coloneqq\;\int_{\hat{T}} \partial_{\hat{x}^m} \hat{\vphi}_i\, \hat{\vphi}_j\, \hat{\vphi}_l \quad \mbox{for }m\in\{1,2\}
\end{equation}
Now we can express the local matrix~$\vecc{G}^1_{T_k}$ from~\eqref{eq:locGm} as
\begin{align*}
\vecc{G}^1_{T_k} &=   B_k^{21}\sum_{l=1}^N D_{kl}(t)
\left(B_k^{22}\int_{\hat{T}}\begin{bmatrix}
\partial_{\hat{x}^1}\hat{\vphi}_1\hat{\vphi}_1\hat{\vphi}_l & \cdots & \partial_{\hat{x}^1}\hat{\vphi}_1\hat{\vphi}_N\hat{\vphi}_l \\
 \vdots                & \ddots & \vdots \\
 \partial_{\hat{x}^1}\hat{\vphi}_N\hat{\vphi}_1\hat{\vphi}_l & \cdots & \partial_{\hat{x}^1}\hat{\vphi}_N\hat{\vphi}_N\hat{\vphi}_l
\end{bmatrix}
-B_k^{21}\int_{\hat{T}}\begin{bmatrix}
\partial_{\hat{x}^2}\hat{\vphi}_1\hat{\vphi}_1\hat{\vphi}_l & \cdots & \partial_{\hat{x}^2}\hat{\vphi}_1\hat{\vphi}_N\hat{\vphi}_l \\
 \vdots                & \ddots & \vdots \\
 \partial_{\hat{x}^2}\hat{\vphi}_N\hat{\vphi}_1\hat{\vphi}_l & \cdots & \partial_{\hat{x}^2}\hat{\vphi}_N\hat{\vphi}_N\hat{\vphi}_l
\end{bmatrix}
\right)
\\
&= \sum_{l=1}^N D_{kl}(t)
\left(B_k^{22} [\hat{\vecc{G}}]_{:,:,l,1}-B_k^{21}  [\hat{\vecc{G}}]_{:,:,l,2}\right)
\end{align*}
and analogously~$\vecc{G}^2_{T_k}$.  
With $\vecc{G}^m=\diag(\vecc{G}_{T_1}^m,\ldots,\vecc{G}_{T_K}^m)$ we can vectorize over all triangles using the Kronecker product as done in the routine \code{assembleMatElemDphiPhiFuncDisc}.

\subsubsection{Assembly of $\vecc{S}$}\label{sec:assembly:globS}
To ease the assembly of~$\vecc{S}$ we split the global matrix as given in~\eqref{eq:globS} into a~block-diagonal part and a~remainder so that~$\vecc{S} = \vecc{S}^\mathrm{diag} + \vecc{S}^\mathrm{offdiag}$ holds.
\par
We first consider the block-diagonal entries of~$\vecc{S}$ consisting of sums of local matrices~$\vecc{S}_{E_{kn}}$, cf.~\eqref{eq:globS:diag} and~\eqref{eq:locSEkn}, respectively.   Our first goal is to transform~$\vecc{S}_{E_{kn}}$ to a~local matrix~$\hat{\vecc{S}}^\mathrm{diag}\in\IR^{N\times N\times3}$ that is independent of the physical triangle~$T_k$.  To this end, we transform the edge integral term~$\int_{E_{kn}} \vphi_{ki}\,\vphi_{kj}$ to the \mbox{$n$th} edge of the reference triangle~$\hat{E}_n$: 
\begin{equation}\label{eq:hatSdiag}
\int_{E_{kn}}\vphi_{ki}\,\vphi_{kj}
=\frac{\abs{E_{kn}}}{\abs{\hat{E}_{n}}} \int_{\hat{E}_{n}} \hat{\vphi}_i(\hat{\vec{x}})~\hat{\vphi}_j(\hat{\vec{x}})\,\dd\hat{\vec{x}} 
=\frac{\abs{E_{kn}}}{\abs{\hat{E}_{n}}} \int_0^1 \hat{\vphi}_i\circ\hat{\vec{\gamma}}_{n}(s)~~\hat{\vphi}_j\circ\hat{\vec{\gamma}}_{n}(s)\,\abs{\hat{\vec{\gamma}}_{n}'(s)}\,\dd s
=\abs{E_{kn}}  \underbrace{\int_0^1 \hat{\vphi}_i\circ\hat{\vec{\gamma}}_{n}(s)~~\hat{\vphi}_j\circ\hat{\vec{\gamma}}_{n}(s) \,\dd s}_{\eqqcolon [\hat{\vecc{S}}^\mathrm{diag}]_{i,j,n}}\;,
\end{equation}
where we used transformation rule~\eqref{eq:trafoRule:E} and $\abs{\hat{\vec{\gamma}}_n'(s)}=\abs{\hat{E}_n}$. The explicit forms of the mappings~$\hat{\vec{\gamma}}_n:[0,1]\rightarrow \hat{E}_n$ can be easily derived:
\begin{equation}\label{eq:gammaMap}
\hat{\vec{\gamma}}_1(s) \;\coloneqq\;
\begin{bmatrix}
1-s\\s
\end{bmatrix}\,,\qquad
\hat{\vec{\gamma}}_2(s) \;\coloneqq\;
\begin{bmatrix}
0\\1-s
\end{bmatrix}\,,\qquad
\hat{\vec{\gamma}}_3(s) \;\coloneqq\;
\begin{bmatrix}
s\\0
\end{bmatrix}\,.
\end{equation}
Thus, we have~$\vecc{S}_{E_{kn}}=\abs{E_{kn}} [\hat{\vecc{S}}^\mathrm{diag}]_{:,:,n}$ allowing to 
define the diagonal blocks of the global matrix~$\vecc{S}^\mathrm{diag}$ using the Kron\-ecker product:
\begin{equation*}
\vecc{S}^\mathrm{diag} \;\coloneqq\;  \sum_{n=1}^3
\begin{bmatrix}
\delta_{E_{1n}\in\setE_\Omega} &   & \\
               & ~\ddots~ & \\
               &          & \delta_{E_{Kn}\in\setE_\Omega}
\end{bmatrix} \otimes [\hat{\vecc{S}}^\mathrm{diag}]_{:,:,n}\;,
\end{equation*}
where $\delta_{E_{kn}\in\setE_\Omega}$ denotes the Kronecker delta.
\par
Next, we consider the off-diagonal blocks of~$\vecc{S}$ stored in~$\vecc{S}^\mathrm{offdiag}$.  For an~interior edge~ $E_{k^-n^-}=E_{k^+n^+}\in\partial T_{k^-}\cap \partial T_{k^+}, \; n^-,n^+\in\{1,2,3\}$ (cf.~Fig.~\ref{fig:T1T2}) we obtain analogously:
\begin{align*}
\hspace{-8mm}\int_{E_{k^-n^-}}\vphi_{k^-i}\,\vphi_{k^+j}
&\;=\;\frac{\abs{E_{k^-n^-}}}{\abs{\hat{E}_{n^-}}} \int_{\hat{E}_{n^-}} \vphi_{k^-i}\circ\vec{F}_{k^-}(\hat{\vec{x}})~~\vphi_{k^+j}\circ\overbrace{\vec{F}_{k^+}\circ\vec{F}_{k^+}^{-1}}^{=\vec{I}}\circ\vec{F}_{k^-}(\hat{\vec{x}}) \,\dd\hat{\vec{x}} 
\;=\;\frac{\abs{E_{k^-n^-}}}{\abs{\hat{E}_{n^-}}} \int_{\hat{E}_{n^-}} \hat{\vphi}_{i}(\hat{\vec{x}})~~\hat{\vphi}_{j}\circ\vec{F}_{k^+}^{-1}\circ\vec{F}_{k^-}(\hat{\vec{x}}) \,\dd\hat{\vec{x}} \\
&\;=\;\abs{E_{k^-n^-}}  \int_0^1 \hat{\vphi}_i\circ\hat{\vec{\gamma}}_{n^-}(s)~~\hat{\vphi}_j\circ\vec{F}_{k^+}^{-1}\circ\vec{F}_{k^-}\circ\hat{\vec{\gamma}}_{n^-}(s) \,\dd s\;.
\end{align*}
Note that~$\vec{F}_{k^+}^{-1}\circ\vec{F}_{k^-}$ maps from~$\hat{T}$ to~$\hat{T}$.  
Since we compute a~line integral the integration domain is further restricted to an edge~$\hat{E}_{n^-}$, $n^-\in\{1,2,3\}$ and its co-domain to an~edge~$\hat{E}_{n^+}$, $n^+\in\{1,2,3\}$.
As a~result, this integration can be boiled down to nine possible maps between the sides of the reference triangle expressed as
\begin{equation*}
\mapEE_{n^-n^+}:\quad \hat{E}_{n^-}\ni\hat{\vec{x}}\mapsto \mapEE_{n^-n^+}(\hat{\vec{x}})\,=\,\vec{F}_{k^+}^{-1}\circ\vec{F}_{k^-}(\hat{\vec{x}})\in \hat{E}_{n^+}
\end{equation*}
for an~arbitrary index pair~$\{k^-,k^+\}$ as described above. The closed-form expressions of the nine cases are:
\begin{align}\label{eq:mapEE}
\mapEE_{11}&:
\begin{bmatrix}
\hat{x}^1\\\hat{x}^2
\end{bmatrix}\mapsto
\begin{bmatrix}
1-\hat{x}^1\\1-\hat{x}^2
\end{bmatrix}\;,
&
\mapEE_{12}&:
\begin{bmatrix}
\hat{x}^1\\\hat{x}^2
\end{bmatrix}\mapsto
\begin{bmatrix}
0\\\hat{x}^2
\end{bmatrix}\;,
&
\mapEE_{13}&:
\begin{bmatrix}
\hat{x}^1\\\hat{x}^2
\end{bmatrix}\mapsto
\begin{bmatrix}
\hat{x}^1\\0
\end{bmatrix}\;,
\nonumber\\
\mapEE_{21}&:
\begin{bmatrix}
\hat{x}^1\\\hat{x}^2
\end{bmatrix}\mapsto
\begin{bmatrix}
1-\hat{x}^2\\\hat{x}^2
\end{bmatrix}\;,
&
\mapEE_{22}&:
\begin{bmatrix}
\hat{x}^1\\\hat{x}^2
\end{bmatrix}\mapsto
\begin{bmatrix}
0\\1-\hat{x}^2
\end{bmatrix}\;,
&
\mapEE_{23}&:
\begin{bmatrix}
\hat{x}^1\\\hat{x}^2
\end{bmatrix}\mapsto
\begin{bmatrix}
\hat{x}^2\\0
\end{bmatrix}\;,
\\
\mapEE_{31}&:
\begin{bmatrix}
\hat{x}^1\\\hat{x}^2
\end{bmatrix}\mapsto
\begin{bmatrix}
\hat{x}^1\\1-\hat{x}^1
\end{bmatrix}\;,
&
\mapEE_{32}&:
\begin{bmatrix}
\hat{x}^1\\\hat{x}^2
\end{bmatrix}\mapsto
\begin{bmatrix}
0\\\hat{x}^1
\end{bmatrix}\;,
&
\mapEE_{33}&:
\begin{bmatrix}
\hat{x}^1\\\hat{x}^2
\end{bmatrix}\mapsto
\begin{bmatrix}
1-\hat{x}^1\\0
\end{bmatrix}\;.
\nonumber
\end{align}
All maps~$\mapEE_{n^-n^+}$ reverse the edge orientation because an~edge shared by triangles~$T^-$ and~$T^+$ will always have different orientations when mapped by $\vec{F}_{k^-}$ and $\vec{F}_{k^+}$; this occurs due to the~counter-clockwise vertex orientation consistently maintained throughout the mesh.
We define $\hat{\vecc{S}}^\mathrm{offdiag}\in\IR^{N\times N\times 3\times 3}$ by
\begin{equation}\label{eq:hatSoffdiag}
[\hat{\vecc{S}}^\mathrm{offdiag}]_{i,j,n^-,n^+} \;\coloneqq\; \int_0^1 \hat{\vphi}_i\circ\hat{\vec{\gamma}}_{n^-}(s)~~\hat{\vphi}_j\circ\mapEE_{n^-n^+}\circ\hat{\vec{\gamma}}_{n^-}(s) \,\dd s
\end{equation}
and thus arrive at
\begin{equation*}
\vecc{S}^\mathrm{offdiag} \coloneqq - \sum_{n^-=1}^3\sum_{n^+=1}^3
\begin{bmatrix}
0                          & \delta_{E_{1n^-}=E_{2n^+}} & \hdots                & \hdots                & \delta_{E_{1n^-}=E_{Kn^+}}\\
\delta_{E_{2n^-}=E_{1n^+}} & 0                          &   \ddots              &                       & \textstyle\vdots  \\ 
\vdots                     &         \ddots             & \ddots                &          \ddots       & \vdots   \\
\vdots                     & {}                         &    \ddots             & 0                     & \delta_{E_{(K-1)n^-}=E_{Kn^+}}\\
\delta_{E_{Kn^-}=E_{1n^+}} &  \hdots                    & \hdots                & \delta_{E_{Kn^-}=E_{(K-1)n^+}} & 0
\end{bmatrix}\otimes [\hat{\vecc{S}}^\mathrm{offdiag}]_{:,:,n^-,n^+}\;.
\end{equation*}
The sparsity structure for off-diagonal blocks depends on the numbering of mesh entities and is given for each combination of $n^-$ and $n^+$ by the list~\code{markE0TE0T} (cf.~Tab.~\ref{tab:lists}).
The routine~\code{assembleMatEdgePhiPhi} assembles the matrices $\vecc{S}^\mathrm{diag}$ and $\vecc{S}^\mathrm{offdiag}$ directly into $\vecc{S}$ with a~code very similar to the formulation above.

\subsubsection{Assembly of~$\vecc{Q}^m$}\label{sec:assembly:globQ}
The assembly of $\vecc{Q}^m$ from equations~\eqref{eq:globQ:diag},~\eqref{eq:globQ:offdiag} is analogous to $\vecc{S}$ since both are constructed from the same terms only differing in constant coefficients. 
Consequently, we can choose the same approach as described in~\ref{sec:assembly:globS}. 
Again, we split the matrix into diagonal and off-diagonal blocks $\vecc{Q}^m~=~\vecc{Q}^{m,\mathrm{diag}}+\vecc{Q}^{m,\mathrm{offdiag}}$ and assemble each separately exploiting transformation rule~\eqref{eq:trafoRule:E}.
This allows to write the diagonal blocks as follows:
\begin{equation*}
\vecc{Q}^{m,\mathrm{diag}} \;\coloneqq\; \frac{1}{2} \sum_{n=1}^3
\begin{bmatrix}
\delta_{E_{1n}\in\setE_\Omega} &   & \\
               & \ddots   & \\
               &          & \delta_{E_{Kn}\in\setE_\Omega}
\end{bmatrix}
\circ
\begin{bmatrix}
\nu_{1n}^m\abs{E_{1n}} &        & \\
         & \ddots & \\
         &        & \nu_{Kn}^m\abs{E_{Kn}}
\end{bmatrix}
 \otimes [\hat{\vecc{S}}^\mathrm{diag}]_{:,:,n}\;,
\end{equation*}
where \enquote{$\circ$} is the operator for the Hadamard product.
\par
The off-diagonal blocks are assembled as before, using the mapping $\mapEE_{n^-n^+}$ from~\eqref{eq:mapEE}. This leads to to a~similar representation as for $\vecc{S}^\mathrm{offdiag}$:
\begin{multline*}
\vecc{Q}^{m,\mathrm{offdiag}} \coloneqq \frac{1}{2} \sum_{n^-=1}^3\sum_{n^+=1}^3
\begin{bmatrix}
0                          & \delta_{E_{1n^-}=E_{2n^+}} & \hdots                & \hdots                & \delta_{E_{1n^-}=E_{Kn^+}} \\
\delta_{E_{2n^-}=E_{1n^+}} & 0                          &   \ddots              &                       & \textstyle\vdots  \\ 
\vdots                     &         \ddots             & \ddots                &          \ddots       & \vdots   \\
\vdots                     & {}                         &    \ddots             & 0                     & \delta_{E_{(K-1)n^-}=E_{Kn^+}} \\
\delta_{E_{Kn^-}=E_{1n^+}}  &  \hdots                   & \hdots                & \delta_{E_{Kn^-}=E_{(K-1)n^+}}  & 0
\end{bmatrix}
\\
\circ
\begin{bmatrix}
\nu_{1n^-}^m\abs{E_{1n^-}} & \hdots & \nu_{1n^-}^m\abs{E_{1n^-}} \\
\vdots &  & \vdots \\
\nu_{Kn^-}^m\abs{E_{Kn^-}} & \hdots & \nu_{Kn^-}^m\abs{E_{Kn^-}}
\end{bmatrix}
\otimes [\hat{\vecc{S}}^\mathrm{offdiag}]_{:,:,n^-,n^+}\;.
\end{multline*}
Once again, we can use a~code close to the mathematical formulation to assemble the matrices~$\vecc{Q}^{m,\mathrm{diag}}$ and~$\vecc{Q}^{m,\mathrm{offdiag}}$.  This is realized in the routine~\code{assembleMatEdgePhiPhiNu}. 
In the implementation, the Hadamard product is replaced by a~call to the built-in function \code{bsxfun} which applies a~certain element-by-element operation (here: \code{@times}) to arrays. Since all columns in the second matrix of the product are the same this makes superfluous explicitly creating this matrix and permits the use of a~single list of all required values instead.

\subsubsection{Assembly of~$\vecc{R}^m$}\label{sec:assembly:globR}
Just as before, we split the block matrices $\vecc{R}^m$ from~\eqref{eq:globR:diag},~\eqref{eq:globR:offdiag} into diagonal and off-diagonal parts as $\vecc{R}^m~=~\vecc{R}^{m,\mathrm{diag}}+\vecc{R}^{m,\mathrm{offdiag}}$. 
Here, integrals consist of three basis functions due to the diffusion coefficient but still can be transformed in the same way.
In diagonal blocks, this takes the form
\begin{equation}\label{eq:hatRdiag}
\int_{E_{kn}} \vphi_{ki}\, \vphi_{kl} \,\vphi_{kj} \;=\; \abs{E_{kn}} \underbrace{
\int_0^1 \hat{\vphi}_i \circ \hat{\vec{\gamma}}_n(s)~~\hat{\vphi}_l\circ\hat{\vec{\gamma}}_n(s)~~\hat{\vphi}_j\circ\hat{\vec{\gamma}}_n(s)\,\dd s
}_{\eqqcolon [\hat{\vecc{R}}^\mathrm{diag}]_{i,j,l,n}}\;,
\end{equation}
which can be used to define a~common multidimensional array~$\hat{\vecc{R}}^\mathrm{diag}\in\IR^{N\times N\times N\times 3}$. This allows to re-write the local block matrix from~\eqref{eq:locR:diag} as
\begin{equation*}
\vecc{R}_{E_{kn}} = \sum_{l=1}^N D_{kl}(t) \,\abs{E_{kn}} \,[\hat{\vecc{R}}^\mathrm{diag}]_{:,:,l,n}\;.
\end{equation*}
Consequently, the assembly of $\vecc{R}^{m,\mathrm{diag}}$ can be formulated as
\begin{equation*}
\vecc{R}^{m,\mathrm{diag}} \;\coloneqq\; \frac{1}{2} \sum_{n=1}^3 \sum_{l=1}^N
\begin{bmatrix}
\delta_{E_{1n}\in\setE_\Omega} &   & \\
               & \ddots & \\
               &        & \delta_{E_{Kn}\in\setE_\Omega}  
\end{bmatrix} \circ 
\begin{bmatrix}
\nu_{1n}^m \abs{E_{1n}}  D_{1l}(t) &   & \\
               & \ddots & \\
               &        &  \nu_{Kn}^m \abs{E_{Kn}}  D_{Kl}(t)
\end{bmatrix} 
\otimes [\hat{\vecc{R}}^\mathrm{diag}]_{:,:,l,n}\;.
\end{equation*}
\par
The off-diagonal entries consist of integrals over triples of basis functions two of which belong to the adjacent triangle $T_{k^+}$, thus making it necessary to apply the mapping $\mapEE_{n^-n^+}$ from~\eqref{eq:mapEE}. Once again, this can be written as
\begin{equation}\label{eq:hatRoffdiag}
\int_{E_{k^-n^-}} \vphi_{k^-i}\, \vphi_{k^+l} \,\vphi_{k^+j} \;=\; \abs{E_{k^-n^-}} 
\underbrace{\int_0^1 \hat{\vphi}_i\circ\hat{\vec{\gamma}}_{n^-}(s)~~\hat{\vphi}_l\circ\mapEE_{n^-n^+}\circ\hat{\vec{\gamma}}_{n^-}(s)~~\hat{\vphi}_j\circ\mapEE_{n^-n^+}\circ\hat{\vec{\gamma}}_{n^-}(s) \,\dd s}_{\eqqcolon\;[\hat{\vecc{R}}^\mathrm{offdiag}]_{i,j,l,n^-,n^+}}
\end{equation}
with a~multidimensional array~$\hat{\vecc{R}}^\mathrm{offdiag}\in\IR^{N\times N\times N\times 3\times 3}$ whose help allows us to carry out the assembly of~$\vecc{R}^{m,\mathrm{offdiag}}$ (component-wise given in~\eqref{eq:globR:offdiag}) by
\begin{align*}
\vecc{R}^{m,\mathrm{offdiag}} \coloneqq \frac{1}{2} & \sum_{n^-=1}^3 \sum_{n^+=1}^3 \sum_{l=1}^N
\begin{bmatrix}
0                          & \delta_{E_{1n^-}=E_{2n^+}} & \cdots                & \cdots                & \delta_{E_{1n^-}=E_{Kn^+}}\\
\delta_{E_{2n^-}=E_{1n^+}} & 0                          &   \ddots              &                       & \textstyle\vdots  \\ 
\vdots                     &         \ddots             & \ddots                &          \ddots       & \vdots   \\
\vdots                     & {}                         &    \ddots             & 0                     & \delta_{E_{(K-1)n^-}=E_{Kn^+}}\\
\delta_{E_{Kn^-}=E_{1n^+}} &  \cdots                    & \cdots                & \delta_{E_{Kn^-}=E_{(K-1)n^+}} & 0
\end{bmatrix}\\
& \circ
\begin{bmatrix}
\nu_{1n^-}^m\abs{E_{1n^-}} & \cdots & \nu_{1n^-}^m\abs{E_{1n^-}} \\
\vdots &  & \vdots \\
\nu_{Kn^-}^m\abs{E_{Kn^-}} & \cdots & \nu_{Kn^-}^m\abs{E_{Kn^-}}
\end{bmatrix}
\circ
\begin{bmatrix}
D_{1l}(t)    & \cdots & D_{Kl}(t) \\
\vdots       &        & \vdots    \\
D_{1l}(t)    & \cdots & D_{Kl}(t)
\end{bmatrix}
\otimes [\hat{\vecc{R}}^\mathrm{offdiag}]_{:,:,l,n^-,n^+}\;.
\end{align*}
The corresponding code is found in~\code{assembleMatEdgePhiPhiFuncDiscNu}.
Again, we make use of the function \code{bsxfun} to carry out the Hadamard product; it is used twice, first to apply the row vector as before and then to apply the column vector of the diffusion coefficient.

\subsubsection{Assembly of~$\vecc{R}^m_\mathrm{D}$}\label{sec:assembly:globRD}

Since the entries of $\vecc{R}^m_\mathrm{D}$ in~\eqref{eq:globR:D} are computed in precisely the same way as for $\vecc{R}^{m,\mathrm{diag}}$ (cf.~\eqref{eq:globR:diag}), the corresponding assembly routine \code{assembleMatEdgePhiIntPhiIntFuncDiscIntNu} consists only of the part of \code{assembleMatEdgePhiPhiFuncDiscNu} which is responsible for the assembly of $\vecc{R}^{m,\mathrm{diag}}$.
It only differs in a~factor and the list of edges, namely all edges in the set $\setE_\mathrm{D}$ for which non-zero entries (here given by \code{markE0Tbdr}) are generated.

\subsubsection{Assembly of~$\vecc{S}_\mathrm{D}$}\label{sec:assembly:globSD}

The same as for $\vecc{R}^m_\mathrm{D}$ holds for $\vecc{S}_\mathrm{D}$  in~\eqref{eq:globS:D} which, in fact, has the same entries in the diagonal blocks as $\vecc{S}^\mathrm{diag}$ (cf.~\eqref{eq:globS:diag}). 
Consequently, the corresponding routine \code{assembleMatEdgePhiIntPhiInt} is again a~subset of~\code{assembleMatEdgePhiPhi}.

\subsubsection{Assembly of~$\vecc{Q}^m_\mathrm{N}$}\label{sec:assembly:globQN}

For the Neumann boundary edges in the set $\setE_\mathrm{N}$, only contributions in the block $\vecc{Q}^m_\mathrm{N}$ in~\eqref{eq:globQm:N} are generated. The responsible routine \code{assembleMatEdgePhiIntPhiIntNu} is also equivalent to the assembly routine for $\vecc{Q}^{m,\mathrm{diag}}$.

\subsubsection[Assembly of JD]{Assembly of~$\vec{J}^m_{\mathrm{D}}$}\label{sec:assembly:globJD}
The entries of $\vec{J}^m_{\mathrm{D}}$ in~\eqref{eq:globJm:D} are transformed using transformation rule~\eqref{eq:trafoRule:E}
\begin{align*}
[\vec{J}_\mathrm{D}^m]_{(k-1)N+i} &\;=\;
\sum_{E_{kn}\in\partial T_k\cap\setE_\mathrm{D}} \nu_{kn}^m \int_{E_{kn}}\vphi_{ki}\,c_{\mathrm{D}}(t) 
\;=\; \sum_{E_{kn}\in\partial T_k\cap\setE_\mathrm{D}} \nu_{kn}^m\frac{\abs{E_{kn}}}{\abs{\hat{E}_n}} \int_{\hat{E}_n} \vphi_{ki}\circ F_{k}(\hat{\vec{x}})\,c_{\mathrm{D}}\Big(t, F_{k}(\hat{\vec{x}})\Big)\,\dd\hat{\vec{x}}
\\
&\;=\; \sum_{E_{kn}\in\partial T_k\cap\setE_\mathrm{D}} \nu_{kn}^m\,\abs{E_{kn}} \int_0^1 \hat{\vphi}_i\circ \hat{\vec{\gamma}}_n(s)\,c_{\mathrm{D}}\Big(t, F_{k}\circ \hat{\vec{\gamma}}_n(s)\Big)\,\dd s\;.
\end{align*}
This integral is then approximated using a~1D~quadrature rule~\eqref{eq:quadrature} on the interval~$(0,1)$
\begin{equation*}
[\vec{J}_\mathrm{D}^m]_{(k-1)N+i} \approx \sum_{E_{kn}\in\partial T_k\cap\setE_\mathrm{D}} \nu_{kn}^m \,\abs{E_{kn}} \sum_{r=1}^R \omega_r\,\hat{\vphi}_i\circ \hat{\vec{\gamma}}_n(q_r)\,c_{\mathrm{D}}\left(t, F_{k}\circ \hat{\vec{\gamma}}_n(q_r)\right)\;,
\end{equation*}
allowing to vectorize the computation over all triangles and resulting in the routine~\code{assembleVecEdgePhiIntFuncContNu}.

\subsubsection[Assembly of KD]{Assembly of~$\vec{K}_{\mathrm{D}}$}\label{sec:assembly:globKD}
The computation of $\vec{K}_\mathrm{D}$ is done similarly to the assembly of $\vec{J}^m_\mathrm{D}$. The component-wise integrals from~\eqref{eq:globK:D} are, once again, transformed to the interval~$[0,1]$ using~\eqref{eq:trafoRule:E}:
\begin{equation*}
[\vec{K}_\mathrm{D}]_{(k-1)N+i} \;=\;
 \sum_{E_{kn}\in\partial T_k \cap \setE_\mathrm{D}} \frac{1}{\abs{E_{kn}}} \int_{E_{kn}} \vphi_{ki} \, c_\mathrm{D}(t)
\;=\;  \sum_{E_{kn}\in\partial T_k \cap \setE_\mathrm{D}} \int_0^1 \hat{\vphi}_i \circ \hat{\vec{\gamma}}_n(s) \, c_\mathrm{D} (t, \vec{F}_k \circ \hat{\vec{\gamma}}_n(s)) \, \dd s
\end{equation*}
effectively canceling out the edge length. Using a~quadrature rule and vectorization, $\vec{K}_\mathrm{D}$ is assembled in the routine~\code{assembleVecEdgePhiIntFuncCont}.

\subsubsection[Assembly of KN]{Assembly of~$\vec{K}_{\mathrm{N}}$}\label{sec:assembly:globKN}
In the integral terms of $\vec{K}_\mathrm{N}$ an~additional basis function from the diffusion
coefficient appears.
As before, the integrals are transformed using transformation 
rules~\eqref{eq:trafoRule:E},~\eqref{eq:gammaMap} and a~1D~quadrature rule~\eqref{eq:quadrature}:
\begin{align*}
\hspace{-8mm}[\vec{K}_\mathrm{N}]_{(k-1)N+i} &=
\sum_{E_{kn}\in\partial T_k \cap \setE_\mathrm{N}} \sum_{l=1}^N D_{kl}(t)
\int_{E_{kn}} \hspace{-2mm}\vphi_{ki}\,\vphi_{kl}\,g_\mathrm{N}(t) 
= \sum_{E_{kn}\in\partial T_k \cap \setE_\mathrm{N}} \hspace{-3mm}|E_{kn}|
\sum_{l=1}^N D_{kl}(t) \int_0^1 \hspace{-2mm} \hat{\vphi}_i \circ \hat{\vec{\gamma}}_n(s)\,
\hat{\vphi}_l \circ \hat{\vec{\gamma}}_n(s)\, 
g_\mathrm{N}\left(t, \vec{F}_k \circ \hat{\vec{\gamma}}_n(s)\right)\, \dd s \\
&\;\approx\; \sum_{E_{kn}\in\partial T_k \cap \setE_\mathrm{N}} |E_{kn}|
\sum_{l=1}^N D_{kl}(t) \sum_{r=1}^R \omega_r \, \hat{\vphi}_i \circ \hat{\vec{\gamma}}_n(q_r)\,
\hat{\vphi}_l \circ \hat{\vec{\gamma}}_n(q_r)\, 
g_\mathrm{N}\left(t, \vec{F}_k \circ \hat{\vec{\gamma}}_n(q_r)\right)\;.
\end{align*}
Once again, using vectorization over all triangles the assembly routine is~\code{assembleVecEdgePhiIntFuncDiscIntFuncCont}.

\subsection{Linear solver}
After assembling all blocks as described in the previous section and assembling
system~\eqref{eq:fullSystem} for time step $t^{n+1}$, a~linear system has to be solved
to yield the solution $\vec{Y}^{n+1}$.
For that we employ \MatOct's \code{mldivide}.

\subsection{Computational performance}\label{sec:performance}
As noted at the beginning of Sec.~\ref{sec:implementation}, we obey a~few implementation
conventions to improve the computational performance of our code, including the paradigm to
avoid re-computation of already existing values.
First of all, this boils down to reassembling only those linear system blocks of~\eqref{eq:timeDepSystem} that are time-dependent.
\par
Secondly, these assembly routines involve repeated evaluations of the basis functions at the
quadrature points of the reference triangle.
As stated in Sec.~\ref{sec:quadrature}, we use quadrature rules of order $2p$ on triangles and
of order $2p+1$ on edges precomputing the values of basis functions in the quadrature points.
This is done in the routine~\code{computeBasesOnQuad} for $\hat{\vphi}_i(\hat{\vec{q}}_r)$, 
$\grad\hat{\vphi}_i(\hat{\vec{q}}_r)$, $\hat{\vphi}_i\circ\hat{\vec{\gamma}}_n(s_r)$, and 
$\hat{\vphi}_i\circ\mapEE_{n^-n^+}\circ\hat{\vec{\gamma}}_n(s_r)$, with $\hat{\vec{q}}_r \in \hat{T}$, 
$s_r \in [0,1]$ given by a~2D~or 1D~quadrature rule, respectively, for all required orders.
The values are stored in global cell arrays \code{gPhi2D}, \code{gGradPhi2D}, \code{gPhi1D}, 
and \code{gThetaPhi1D}, allowing to write in the assembly routines, e.\,g., 
\code{gPhi2D\{qOrd\}(:, i)} to obtain the values of the $i$-th basis function on all 
quadrature points of a~quadrature rule of order \code{qOrd}.

\subsubsection{Estimated memory usage}
Two resources limit the problem sizes that can be solved: computational time and available memory.
The first one is a 'soft' limit---in contrast to exceeding the amount of available
memory which will cause the computation to fail.
Hence, we will give an approximate estimate of the memory requirements of the presented code to 
allow gauging the problem sizes and polynomial orders one can compute with the hardware at hand.

\paragraph{Grid data structures} 
The size of the grid data structures depends on the number of mesh entities.
To give a~rough estimate, we will assume certain simplifications: Each triangle has 
three incident edges and each edge (disregarding boundary edges) has two incident
triangles, hence it holds $3\,\card{\setT}\approx2\,\card{\setE}$. 
Additionally, the number of vertices is usually less than the number of triangles,
i.\,e., $\card{\setV}\lesssim\card{\setT}$.
Using those assumptions, the memory requirement of the grid data structures 
(cf.~Tab.~\ref{tab:lists} plus additional lists not shown there) amounts to 
$\approx89\,\card{\setT}\cdot8$~Bytes.

\paragraph{Degrees of freedom}
The memory requirements for system~\eqref{eq:fullSystem} largely depend on the sparsity structure
of the matrix which varies with the numbering of the mesh entities, the number of boundary edges, etc.
In \MatOct, the memory requirements for a~sparse matrix $A\in\IR^{n\times n}$ on a~64-bit machine 
can be approximated by $16\cdot\mathrm{nnz}(A) + 8\cdot n + 8$ Bytes, with $\mathrm{nnz}(A)$ being
the number of non-zero entries in $A$.
Coefficients, like $C_{kj}$, $D_{kl}$, $F_{kl}$, $k\in\{1,\ldots,K\}$, $j, l \in\{1,\ldots,N\}$,
and the right-hand side entries $\vec{J}_\mathrm{D}^m$, $\vec{K}_\mathrm{D}$, $\vec{K}_\mathrm{N}$, $\vec{L}$
are stored in full vectors, each of which requires $K N\cdot 8$~Bytes.
\par
The blocks of the system matrices $\vecc{A}(t)$ and $\vecc{W}$ can be divided into two groups:
(i) blocks built from element-wise integrations and (ii) blocks built from edge contributions.
For the first case we showed in Sec.~\ref{sec:defBlocks:T} that these have a~block diagonal
structure due to the local support of the basis functions.
Consequently, each contains $K$~blocks of size~$N\times N$ with nonzero entries.
Blocks from edge integrals also have block diagonal entries but additionally for each of the three 
edges of an element two nonzero blocks exist.
We neglect the blocks for boundary edges here (i.\,e., $\vecc{Q}_\mathrm{N}$, $\vecc{R}_\mathrm{D}$, $\vecc{S}_\mathrm{D}$),
since these hold only entries for edges that are not contained in the interior edge blocks.
\par
Before solving for the next time step, these blocks are assembled into the system matrices and
the right-hand side vector, effectively doubling the memory requirement.
Combining these estimates, this sums to a~memory requirement of 
$\approx (27 KN + 82 KN^2)\cdot 8$~Bytes.
Compared to $(18(KN)^2 + 9KN)\cdot 8$~Bytes alone for the assembled system, 
when using full matrices, this is still a~reasonable number and emphasizes once more
the need for sparse data structures.

\paragraph{Total memory usage}
Note that all these values are highly dependent on the connectivity graph of the mesh and additional
overhead introduced by \MatOct~(e.\,g., for GUI, interpreter, \code{cell}-data 
structures, temporary storage of built-in routines, etc.).
Other blocks, e.\,g., blocks on the reference element, like $\hat{\vecc{M}}$, counters,
helper variables, lookup tables for the basis functions, etc., don't scale with the mesh size
and are left out of these estimates.
Hence, the numbers given here should be understood as a~lower bound.
Combining the partial results for the memory usage, we obtain the total amount of
$\approx (90 + 43 N + 68 N^2)K\cdot 8$~Bytes.
This means, computing with quadratic basis functions on a grid of 10\,000~triangles
requires at least 214~MBytes of memory and should therefore be possible on any current machine. 
However, computing with a grid of half a million elements and polynomials of order~4 requires more 
than 60~GBytes of memory requiring a~high-end workstation.

\subsubsection{Computation time}
An~extensive performance model for our implementation of the DG~method 
exceeds the scope of this publication.
Instead we name the most time consuming parts of our implementation and give an insight
about computation times to be expected on current hardware for different problem sizes and
approximation orders.
\par
\Matlab's profiler is a~handy tool to investigate the runtime distribution within a~program.
We profiled a~time-dependent problem on a~grid with 872~triangles and 100~time steps on
an~Intel \mbox{Core i7-860} CPU (4~cores, 8~threads) with 8~GBytes of RAM and \mbox{\Matlab~R2014a~(8.3.0.532)}.
For low and moderate polynomial orders ($p=0,1,2$) the largest time share (50\,--\,70\,\%)
was spent in the routine~\code{assembleMatEdgePhiPhiFuncDiscNu} which assembles the contributions of the edge
integrals in term $\VI$ in~\eqref{eq:spacediscretesystem:b}.
Due to the presence of the time-dependent diffusion coefficient it has to be executed in
every timestep, and most of its time is spent in the functions~\code{bsxfun(@times,...)}
(applies the Hadamard product) and~\code{kron} (performs the assembly).
The second most expensive part is then the solver itself, for which we employ~\code{mldivide}.
When going to higher polynomial orders ($p\ge 3$) this part even becomes the most expensive one,
simply due to the larger number of degrees of freedom.
In such cases the routine~\code{assembleMatElemDphiPhiFuncDisc} also takes a~share worth mentioning (up to 15\,\%),
which again assembles a~time-dependent block due to the diffusion coefficient.
Any other part takes up less than~5\,\% of the total computation time.
Although our code is not parallelized, \Matlab's built-in routines (in particular,~\code{mldivide}
and~\code{bsxfun}) make extensive use of multithreading.
Hence, these results are machine dependent and, especially on machines with a~different number of cores, the
runtime distribution might be different.
\par
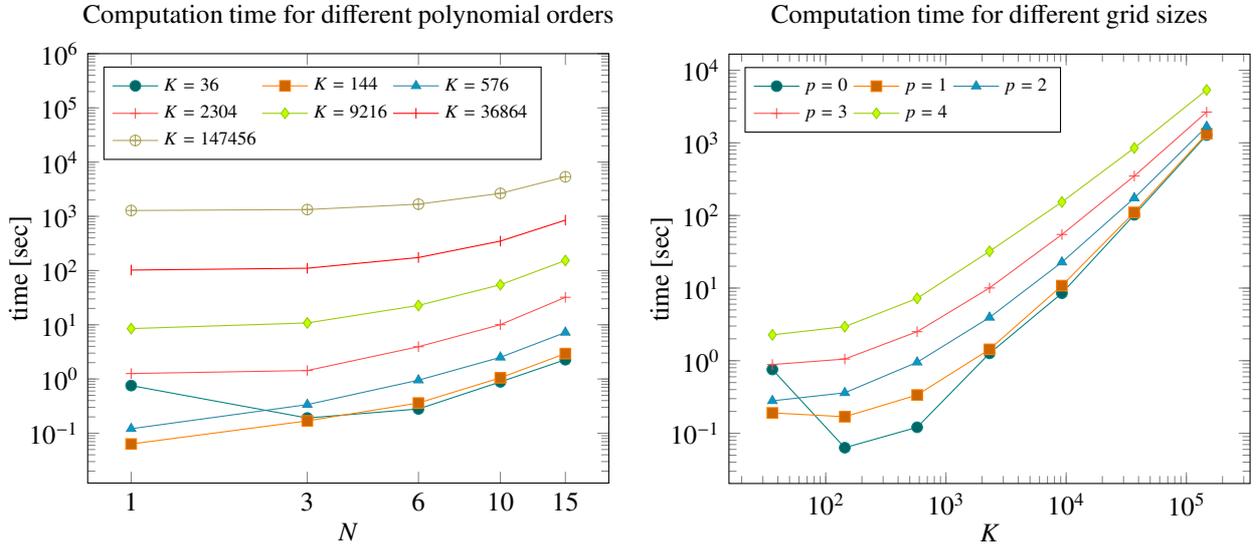
\begin{figure}[h!]
\begin{tabularx}{\linewidth}{@{}LL@{}}
\begin{tikzpicture}
  \begin{axis}[
    title=Computation time for different polynomial orders,
    xtick=data, xticklabels={1,3,6,10,15}, 
    ymode=log, ymax=1e6,
    xlabel={$N$}, ylabel={time [sec]}, 
    x label style={at={(.5,.02)}}, y label style={at={(.05,.5)}},
    legend entries={$K=36$,$K=144$,$K=576$,$K=2304$,$K=9216$,$K=36864$,$K=147456$},
    legend columns=3, legend cell align=left, legend style={font=\scriptsize},
    legend pos=north west,
    cycle list name=exotic,
    ]
    \addplot+[solid,mark=*]         table[x expr=ln(\thisrow{N})/ln(10), y=K36    ] 
      {tikzpicture/runtimes_N.dat};
    \addplot+[solid,mark=square*]   table[x expr=ln(\thisrow{N})/ln(10), y=K144   ] 
      {tikzpicture/runtimes_N.dat};
    \addplot+[solid,mark=triangle*] table[x expr=ln(\thisrow{N})/ln(10), y=K576   ] 
      {tikzpicture/runtimes_N.dat};
    \addplot+[solid,mark=+]         table[x expr=ln(\thisrow{N})/ln(10), y=K2304  ] 
      {tikzpicture/runtimes_N.dat};
    \addplot+[solid,mark=diamond*]  table[x expr=ln(\thisrow{N})/ln(10), y=K9216  ] 
      {tikzpicture/runtimes_N.dat};
    \addplot+[solid,mark=|]         table[x expr=ln(\thisrow{N})/ln(10), y=K36864 ] 
      {tikzpicture/runtimes_N.dat};
    \addplot+[solid,mark=oplus]     table[x expr=ln(\thisrow{N})/ln(10), y=K147456] 
      {tikzpicture/runtimes_N.dat};
  \end{axis}
\end{tikzpicture} & \begin{tikzpicture}
  \begin{axis}[
    title=Computation time for different grid sizes,
    xmode=log, ymode=log,
    xlabel={$K$}, ylabel={time [sec]},
    x label style={at={(.5,.02)}}, y label style={at={(.05,.5)}},
    legend entries={$p=0$, $p=1$, $p=2$, $p=3$, $p=4$},
    legend columns=3, legend cell align=left, legend style={font=\scriptsize},
    legend pos=north west,
    cycle list name={exotic},
    ]
    \addplot+[mark=*]         table[x=K, y=p0] {tikzpicture/runtimes_K.dat};
    \addplot+[mark=square*]   table[x=K, y=p1] {tikzpicture/runtimes_K.dat};
    \addplot+[mark=triangle*] table[x=K, y=p2] {tikzpicture/runtimes_K.dat};
    \addplot+[mark=+]         table[x=K, y=p3] {tikzpicture/runtimes_K.dat};
    \addplot+[mark=diamond*]  table[x=K, y=p4] {tikzpicture/runtimes_K.dat};
  \end{axis}
\end{tikzpicture}
\end{tabularx}
\caption{Computation times for approximation order \emph{(left)} and varying grid size \emph{(right)}.
Polynomial orders~$p=0,1,2,3,4$ correspond to~$N=1,3,6,10,15$ local degrees of freedom, respectively.
Measured on~$4\times$~Intel~Xeon~E7-4870~CPUs~(each 10~cores, 20~threads) with 500~GBytes~RAM
and \Matlab~R2014a~(8.3.0.532).}
\label{fig:runtimes}
\end{figure}

For a~sufficiently large number of time steps one can disregard the execution time of the initial
computations (generation of grid data, computation of basis function lookup tables and reference
element blocks, etc.); then the total runtime scales linearly with the number of time steps.
Therefore, we investigate the runtime behavior of the stationary test case as described in 
Sec.~\ref{sec:codeverification} and plot the computation times against the number of local degrees of freedom~$N$ and
the grid size in Fig.~\ref{fig:runtimes}.
This shows that the computation time increases with~$N$ but primarily 
depends on the number of elements~$K$.
We also observe that doubling the number of elements increases the computation time by more than
a~factor of two which is related to the fact that some steps of the algorithm (e.\,g., the linear
solver) have a~complexity that does not depend linearly on the number of degrees of freedom.

\subsection{Code verification}\label{sec:codeverification}
The code is verified by showing that the numerically estimated orders of convergences match the
analytically predicted ones for prescribed smooth solutions.  
Since the spatial discretization is more complex than the time discretization by far, 
we restrict ourselves to the stationary version of~\eqref{eq:diffusion} in this section. 
\par
We choose the exact solution~$c(\vec{x}) \coloneqq \cos(7x^1)\,\cos(7x^2)$ and the diffusion 
coefficient~$d(\vec{x}) \coloneqq \exp(x^1 + x^2)$ on the domain~$\Omega\coloneqq (0,1)^2$ with Neumann
boundaries at $x^2=0,1$ and Dirichlet boundaries elsewhere.  
The data~$c_\mathrm{D}$, $g_\mathrm{N}$, and~$f$ are derived algebraically by inserting~$c$ and~$d$ 
into~\eqref{eq:diffusion}.  
We then compute the solution~$c_{h_j}$ for a~sequence of 
increasingly finer meshes with element sizes~$h_j$, where the coarsest grid~$\setT_{h_0}$ covering~$\Omega$ is 
an~irregular grid, and each finer grid is obtained by regular refinement of its predecessor.
The discretization errors are computed according to~Sec.~\ref{sec:discerror}\,, and Tab.~\ref{tab:verification} 
contains the results demonstrating the (minimum) order of convergence~$\alpha$ in~$h$ estimated by
\begin{equation*}
\alpha\;\coloneqq\;\ln\Bigg(\frac{\|c_{h_{j-1}} - c\|_{L^2(\Omega)}}{\|c_{h_j} - c\|_{L^2(\Omega)}}\Bigg)\Bigg/\ln\Bigg(\frac{h_{j-1}}{h_j}\Bigg)\;.
\end{equation*}
\begin{table}[h]
\begin{tabularx}{\linewidth}{@{}cCcCcCcCcCc@{}}\toprule
$p$ & 0 & 0  & 1 & 1 & 2 & 2 & 3 & 3 & 4 & 4\\
$j$ & $\|c_h-c\|$ & $\alpha$  & $\|c_h-c\|$ & $\alpha$ & $\|c_h-c\|$ & $\alpha$ & $\|c_h-c\|$ & $\alpha$ & $\|c_h-c\|$ & $\alpha$\\\midrule
0  & 2.37E--1 & \phantom{--} --- &  7.70E--2 &  --- &  1.45E--2 &  --- &  1.97E--3  &  --- &  4.18E--4  &  --- \\
1  & 7.67E--2 & \phantom{--}1.63 &  2.46E--2 & 1.65 &  1.26E--3 & 3.52 &  1.53E--4  & 3.69 &  1.26E--5  & 5.05 \\
2  & 8.94E--2 &          -- 0.22 &  6.67E--3 & 1.88 &  1.11E--4 & 3.51 &  1.02E--5  & 3.90 &  3.78E--7  & 5.06 \\
3  & 9.39E--2 &          -- 0.07 &  1.71E--3 & 1.96 &  1.10E--5 & 3.33 &  6.42E--7  & 3.99 &  1.16E--8  & 5.03 \\
4  & 9.48E--2 &          -- 0.01 &  4.32E--4 & 1.99 &  1.22E--6 & 3.18 &  3.94E--8  & 4.03 &  3.59E--10 & 5.01 \\
5  & 9.49E--2 & \phantom{--}0.00 &  1.08E--4 & 2.00 &  1.44E--7 & 3.09 &  2.43E--9  & 4.02 &  1.11E--11 & 5.01 \\
6  & 9.49E--2 & \phantom{--}0.00 &  2.71E--5 & 2.00 &  1.75E--8 & 3.04 &  1.51E--10 & 4.01 &  3.64E--13 & 4.94 \\
\bottomrule
\end{tabularx}
\caption{Discretization errors measured in~$L^2(\Omega)$ and estimated orders of convergences using the penalty~$\eta=1$. We have~$h_j = \frac{1}{3\cdot2^j}$ and $K=36\cdot4^j$ triangles in the \mbox{$j$th} refinement level.}
\label{tab:verification}
\end{table}

\subsection{Visualization}\label{sec:visualization}
In order to get a~deeper insight into the data associated with a~grid~$\setT_h$ or with a~discrete 
variable from~\mbox{$\IP_p(\setT_h)$}, we provide the routines~\code{visualizeGrid}
and \code{visualizeDataLagr}, respectively (see~Sec.~\ref{sec:routines} for documentation).  
Since our code accepts arbitrary basis functions, in particular the modal basis functions 
of~Sec.~\ref{sec:basisfunctions}, we have to sample those at the Lagrangian points on each 
triangle (i.\,e.~the barycenter of~$T$ for~$\IP_0(T)$, the vertices of~$T$ for~$\IP_1(T)$ and the 
vertices and edge barycenters for~$\IP_2(T)$).  
This mapping from the DG~to the Lagrangian basis is realized in~\code{projectDataDisc2DataLagr}. 
The representation of a~discrete quantity in the latter basis which is as the 
DG~representation a~\mbox{$K\times N$}~matrix, is then used to generate a~VTK-file~\cite{VTKUsersGuide}.
These can be visualized and post-processed, e.\,g., by Paraview~\cite{Paraview}.
Unfortunately, the current \mbox{version~4.2.0} does not visualize quadratic functions as such but instead uses piecewise linear approximations consisting of four pieces per element.
\par
\begin{figure}[h]
\centering
\includegraphics[width=0.3\linewidth]{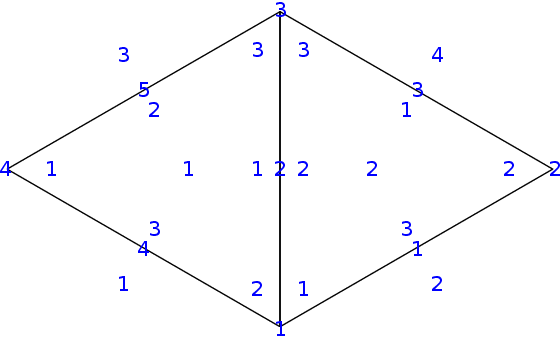}
\caption{Visualization of a~grid with vertices~$\vec{a}_{11}=\transpose{[-\sqrt{3},0]}$, 
$\vec{a}_{12}=\vec{a}_{21}=\transpose{[0,-1]}$, $\vec{a}_{13}=\vec{a}_{23}=\transpose{[0,1]}$, 
$\vec{a}_{22}=\transpose{[\sqrt{3},0]}$ (cf.~Sec.~\ref{sec:visualization}) with triangle 
numbers, global and local vertex numbers, and global and local edge
numbers. Global numbers are printed on, local numbers next to the respective mesh entity.
The boundary IDs which are used to associate parts of the boundary with specific boundary conditions are printed on the exterior of each boundary edge.}
\label{fig:gridvisualization}
\end{figure}
\begin{figure}[h]
\begin{tabularx}{\linewidth}{@{}CCC@{}}
\includegraphics[width=\linewidth]{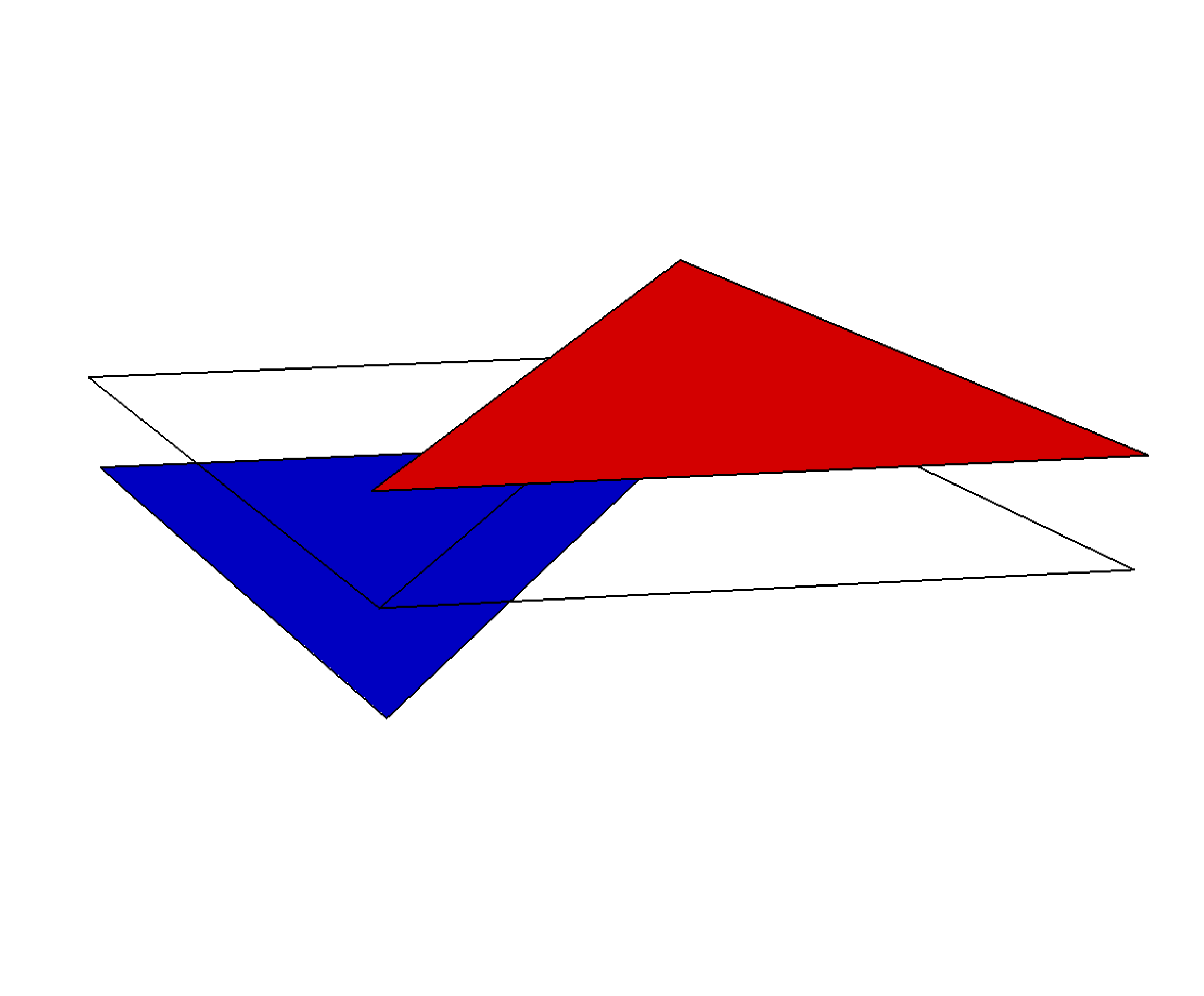} & \includegraphics[width=\linewidth]{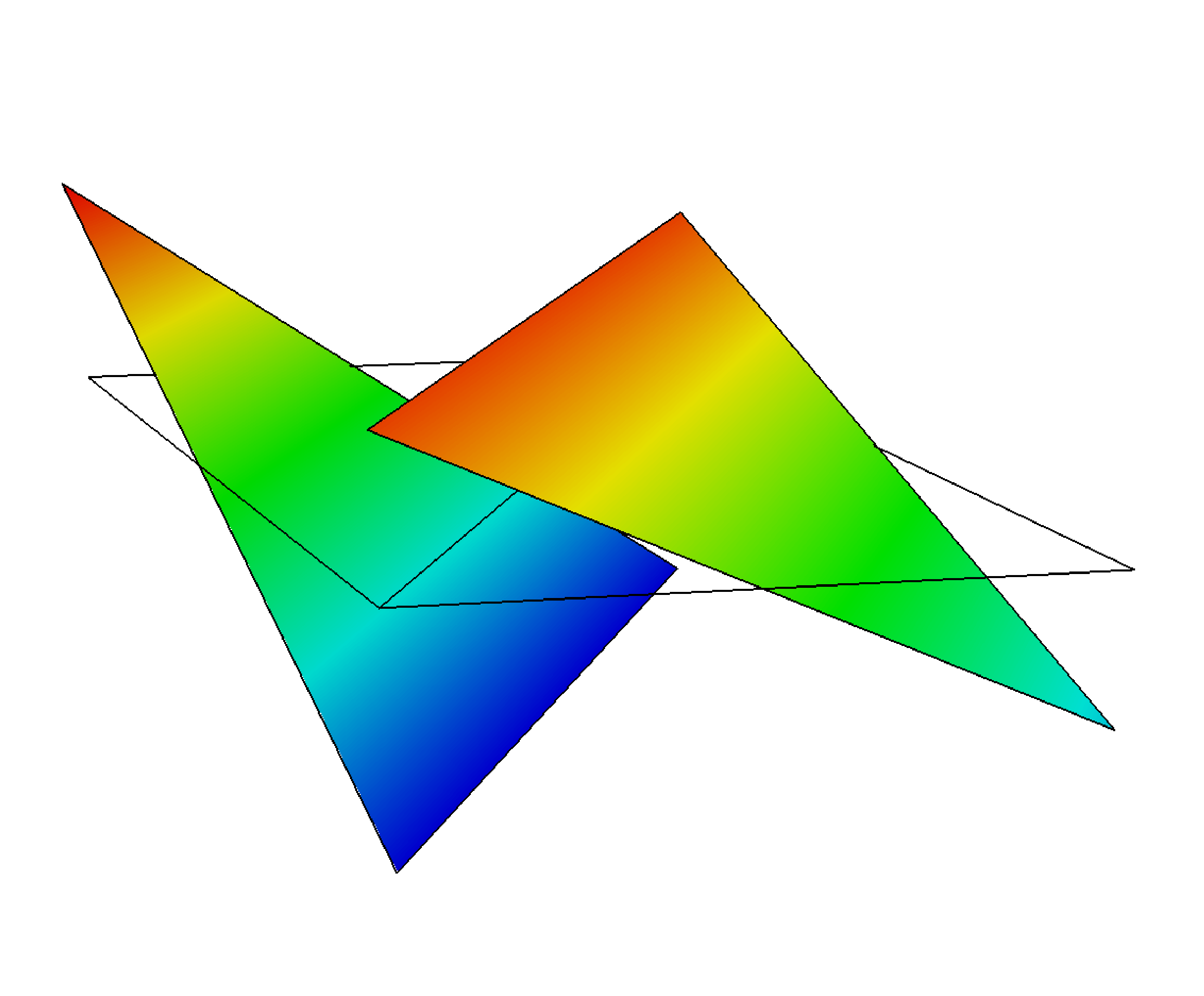} & \includegraphics[width=\linewidth]{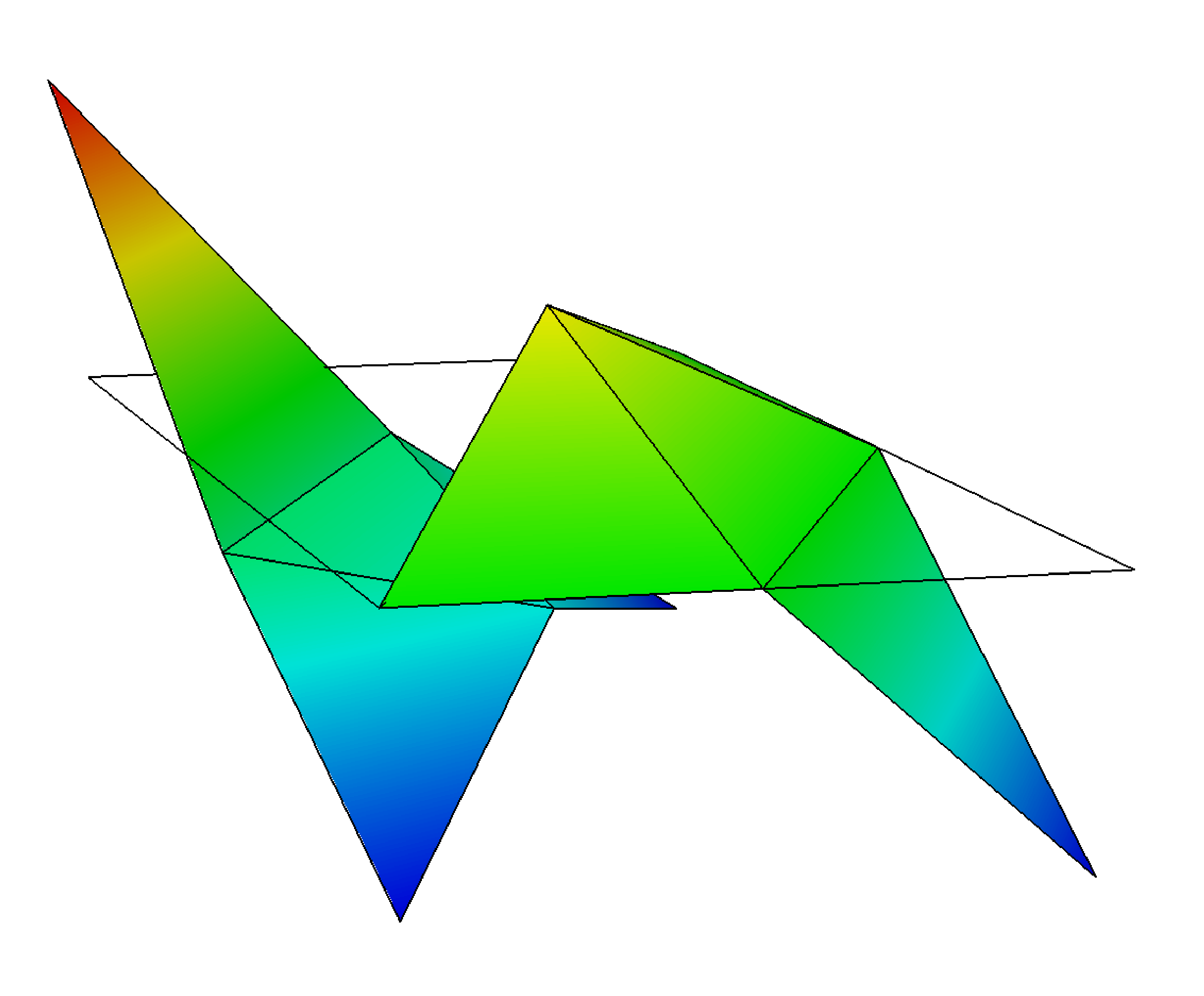}\\
\scriptsize range: $[-2/3,1/3]$ & \scriptsize range: $[-8/5,6/5]$ & \scriptsize range: $[-2,2]$ \\
\scriptsize \texttt{fDOF1.1.vtu} & \scriptsize \texttt{fDOF3.1.vtu} & \scriptsize \texttt{fDOF6.1.vtu}
\end{tabularx}
\caption{Visualization of a~discontinuous function, represented with constant, linear, and 
quadratic basis functions (left to right) using Paraview~(cf.~Sec.~\ref{sec:visualization}).
The underlying mesh is drawn in black.}
\label{fig:visualization}
\end{figure}

As an example, the following code generates a~grid with two triangles and visualizes it using
\code{visualizeGrid}; then a~quadratic, discontinuous function is projected into the DG~space
for~$p\in\{0,1,2\}$ and written to VTK-files using~\code{visualizeDataLagr}.
The resulting outputs are shown in~Fig.~\ref{fig:gridvisualization} and~Fig.~\ref{fig:visualization}.
\begin{lstlisting}
g = generateGridData([0,-1; sqrt(3),0; 0,1; -sqrt(3),0], [4,1,3; 1,2,3]);
g.idE = (abs(g.nuE(:,2))>0).*((g.nuE(:,1)>0)+(g.nuE(:,2)>0)*2+1);
visualizeGrid(g)
fAlg  = @(X1, X2) (X1<0).*(X1.^2-X2.^2-1) + (X1>=0).*(-X1.^2-X2.^2+1);
for N = [1, 3, 6]
  p = (sqrt(8*N+1)-3)/2;  quadOrd = max(2*p, 1);  computeBasesOnQuad(N);
  fDisc   = projectFuncCont2DataDisc(g, fAlg, quadOrd, integrateRefElemPhiPhi(N));
  fLagr = projectDataDisc2DataLagr(fDisc);
  visualizeDataLagr(g, fLagr, 'funname', ['fDOF', int2str(N)], 1)
end % for
\end{lstlisting}

\section{Register of Routines}\label{sec:routines}
We list here all routines of our implementation in alphabetic order.  
For the reason of compactness, we waive the check for correct function arguments, e.\,g., by means of 
routines as~\code{assert}.  However, it is strongly recommended to catch exceptions if the code is to be extended. 
The argument~\code{g}~is always a~\code{struct} representing the triangulation~$\setT_h$ 
(cf.~Sec.~\ref{sec:grid}), the argument~\code{N}~is always the number of local basis functions~$N$.  
A~script that demonstrates the application of the presented routines is given in~\code{main.m}.

\begin{lstlisting}[title={%
\code{ret = assembleMatEdgePhiIntPhiIntFuncDiscIntNu(g, markE0Tbdr, refEdgePhiIntPhiIntPhiInt, dataDisc)} assembles the matrices $\vecc{R}_\mathrm{D}^m$, $m\in\{1,2\}$ according to Sec.~\ref{sec:assembly:globRD}. It is essentially the same routine as the diagonal part of~\code{assembleMatEdgePhiPhiFuncDiscNu} but carried out only for the Dirichlet boundary edges marked by~\code{markE0Tbdr}.
}%
]
function ret = assembleMatEdgePhiIntPhiIntFuncDiscIntNu(g, markE0Tbdr, refEdgePhiIntPhiIntPhiInt, dataDisc)
[K, N] = size(dataDisc);  
ret = cell(2, 1); ret{1} = sparse(K*N, K*N); ret{2} = sparse(K*N, K*N);
for n = 1 : 3
  RDkn = markE0Tbdr(:,n) .* g.areaE0T(:,n);
  for l = 1 : N
    ret{1} = ret{1} + kron(spdiags(RDkn.*g.nuE0T(:,n,1).*dataDisc(:,l),0,K,K), refEdgePhiIntPhiIntPhiInt(:,:,l,n));
    ret{2} = ret{2} + kron(spdiags(RDkn.*g.nuE0T(:,n,2).*dataDisc(:,l),0,K,K), refEdgePhiIntPhiIntPhiInt(:,:,l,n));
  end
end % for
end % function
\end{lstlisting}

\begin{lstlisting}[title={%
\code{ret = assembleMatEdgePhiIntPhiInt(g, markE0Tbdr, refEdgePhiIntPhiInt)} assembles the matrix~$\vecc{S}_\mathrm{D}$ according to Sec.~\ref{sec:assembly:globSD}. It is similar to the diagonal part of~\code{assembleMatEdgePhiPhi} but carried out for Dirichlet boundary edges, which are marked in~\code{markE0Tbdr}.
}%
]
function ret = assembleMatEdgePhiIntPhiInt(g, markE0Tbdr, refEdgePhiIntPhiInt)
K = g.numT;  N = size(refEdgePhiIntPhiInt, 1);
ret = sparse(K*N, K*N);
for n = 1 : 3
  ret = ret + kron(spdiags(markE0Tbdr(:,n),0,K,K), refEdgePhiIntPhiInt(:,:,n));
end % for
end % function
\end{lstlisting}

\begin{lstlisting}[title={%
\code{ret = assembleMatEdgePhiIntPhiIntNu(g, markE0Tbdr, refEdgePhiIntPhiInt)} assembles the matrices~$\vecc{Q}_\mathrm{N}^m$,
$m\in\{1,2\}$ according to Sec.~\ref{sec:assembly:globQN}. It is essentially the same routine as the diagonal part of~\code{assembleMatEdgePhiPhiNu}, only with \code{markE0Tbdr} marking the Neumann boundary edges instead of interior edges.
}%
]
function ret = assembleMatEdgePhiIntPhiIntNu(g, markE0Tbdr, refEdgePhiIntPhiInt)
K = g.numT;  N = size(refEdgePhiIntPhiInt, 1);
ret = cell(2, 1); ret{1} = sparse(K*N, K*N); ret{2} = sparse(K*N, K*N);
for n = 1 : 3
  QNkn = markE0Tbdr(:,n) .* g.areaE0T(:,n);
  ret{1} = ret{1} + kron(spdiags(QNkn .* g.nuE0T(:,n,1), 0,K,K), refEdgePhiIntPhiInt(:,:,n));
  ret{2} = ret{2} + kron(spdiags(QNkn .* g.nuE0T(:,n,2), 0,K,K), refEdgePhiIntPhiInt(:,:,n));
end % for
end % function
\end{lstlisting}

\begin{lstlisting}[title={%
\code{ret = assembleMatEdgePhiPhi(g, markE0Tint, refEdgePhiIntPhiInt, refEdgePhiIntPhiExt)} assembles a~matrix containing integrals over interior edges of products of two basis functions. This corresponds to the matrix~$\vecc{S}$ according to Sec.~\ref{sec:assembly:globS}.  The arguments are the same as for~\code{assembleMatEdgePhiPhiNu}.
}%
]
function ret = assembleMatEdgePhiPhi(g, markE0Tint, refEdgePhiIntPhiInt, refEdgePhiIntPhiExt)
K = g.numT;  N = size(refEdgePhiIntPhiInt, 1);
ret = sparse(K*N, K*N);
for n = 1 : 3
  ret = ret + kron(spdiags(markE0Tint(:,n),0,K,K), refEdgePhiIntPhiInt(:,:,n));
end % for
for nn = 1 : 3
  for np = 1 : 3
    ret = ret - kron(g.markE0TE0T{nn, np}, refEdgePhiIntPhiExt(:,:,nn,np));
  end % for
end % for
end % function
\end{lstlisting}

\begin{lstlisting}[title={%
\code{ret = assembleMatEdgePhiPhiFuncDiscNu(g, markE0Tint, refEdgePhiIntPhiIntPhiInt, refEdgePhiIntPhiExtPhiExt, dataDisc)} assembles two matrices containing integrals over interior edges of products of two basis functions with a~discontinuous coefficient function and with a component of the edge normal. They are returned in a~$2\times1$~\code{cell}~variable. This corresponds to the matrices $\vecc{R}^m$, $m\in\{1,2\}$ according to Sec.~\ref{sec:assembly:globR}. The input arguments~\code{refEdgePhiIntPhiIntPhiInt} and~\code{refEdgePhiIntPhiExtPhiExt} store the local matrices~$\hat{\vecc{R}}^{m, \mathrm{diag}}$ and~$\hat{\vecc{R}}^{m, \mathrm{offdiag}}$ as defined in~\eqref{eq:hatRdiag} and~\eqref{eq:hatRoffdiag}, respectively. They can be computed by~\code{integrateRefEdgePhiIntPhiIntPhiInt} and~\code{integrateRefEdgePhiIntPhiExtPhiExt}. Each triangle's interior edges are marked in~\code{markE0Tint} as in~\code{assembleMatEdgePhiPhiNu}. A~representation of the diffusion coefficient in the polynomial space is stored in~\code{dataDisc} and can be computed by~\code{projectFuncCont2DataDisc}. The Hadamard product is carried out by~\code{bsxfun(@times,...)}. Note the transposed application of~$D_{kl}(t)$ in the second part of the routine as a~result of the diffusion coefficient being taken from the neighboring triangle~$T_{k^+}$.
}%
]
function ret = assembleMatEdgePhiPhiFuncDiscNu(g, markE0Tint, refEdgePhiIntPhiIntPhiInt, refEdgePhiIntPhiExtPhiExt, dataDisc)
[K, N] = size(dataDisc);  
ret = cell(2, 1); ret{1} = sparse(K*N, K*N); ret{2} = sparse(K*N, K*N);
for nn = 1 : 3
  Rkn = 0.5 * g.areaE0T(:,nn);
  for np = 1 : 3
    markE0TE0TtimesRkn1 = bsxfun(@times, g.markE0TE0T{nn, np}, Rkn.*g.nuE0T(:,nn,1));
    markE0TE0TtimesRkn2 = bsxfun(@times, g.markE0TE0T{nn, np}, Rkn.*g.nuE0T(:,nn,2));
    for l = 1 : N
      ret{1} = ret{1} + kron(bsxfun(@times, markE0TE0TtimesRkn1, dataDisc(:,l).'), refEdgePhiIntPhiExtPhiExt(:,:,l,nn,np));
      ret{2} = ret{2} + kron(bsxfun(@times, markE0TE0TtimesRkn2, dataDisc(:,l).'), refEdgePhiIntPhiExtPhiExt(:,:,l,nn,np));
    end % for
  end % for
  Rkn = Rkn .* markE0Tint(:, nn);
  for l = 1 : N
    ret{1} = ret{1} + kron(spdiags(Rkn.*g.nuE0T(:,nn,1).*dataDisc(:,l),0,K,K), refEdgePhiIntPhiIntPhiInt(:,:,l,nn));
    ret{2} = ret{2} + kron(spdiags(Rkn.*g.nuE0T(:,nn,1).*dataDisc(:,l),0,K,K), refEdgePhiIntPhiIntPhiInt(:,:,l,nn));
  end % for
end % for
end % function
\end{lstlisting}

\begin{lstlisting}[title={%
\code{ret = assembleMatEdgePhiPhiNu(g, markE0Tint, refEdgePhiIntPhiInt, refElemPhiIntPhiExt)} assembles two matrices containing integrals over interior edges of products of two basis functions with a component of the edge normal. They are returned in a~$2\times1$~\code{cell}~variable. This corresponds to the matrices $\vecc{Q}^m$, $m\in\{1,2\}$ according to Sec.~\ref{sec:assembly:globQ}. The input arguments~\code{refEdgePhiIntPhiInt} and~\code{refElemPhiIntPhiExt} store the local matrices~$\hat{\vecc{Q}}^{m, \mathrm{diag}}=\hat{\vecc{S}}^\mathrm{diag}$ and~$\hat{\vecc{Q}}^{m, \mathrm{offdiag}}=\hat{\vecc{S}}^\mathrm{offdiag}$ as given in~\eqref{eq:hatSdiag} and~\eqref{eq:hatSoffdiag}, respectively. They can be computed by~\code{integrateRefEdgePhiIntPhiInt} and~\code{integrateRefEdgePhiIntPhiExt}. Similarly to~\code{assembleVecEdgePhiIntFuncContNu}, the argument~\code{markE0Tint} is a~$K\times3$~\code{logical} array that marks each triangle's interior edges. Note the use of \code{bsxfun(@times,...)} to carry out the Hadamard product without building the
full coefficient matrix.
}%
]
function ret = assembleMatEdgePhiPhiNu(g, markE0Tint, refEdgePhiIntPhiInt, refEdgePhiIntPhiExt)
K = g.numT;  N = size(refEdgePhiIntPhiInt, 1);
ret = cell(2, 1); ret{1} = sparse(K*N, K*N);  ret{2} = sparse(K*N, K*N);
for nn = 1 : 3
  Qkn = 0.5 * g.areaE0T(:,nn);
  for np = 1 : 3
    ret{1} = ret{1} + ...
      kron(bsxfun(@times, g.markE0TE0T{nn,np}, Qkn .* g.nuE0T(:,nn,1)), refEdgePhiIntPhiExt(:,:,nn,np));
    ret{2} = ret{2} + ...
      kron(bsxfun(@times, g.markE0TE0T{nn,np}, Qkn .* g.nuE0T(:,nn,2)), refEdgePhiIntPhiExt(:,:,nn,np));
  end % for
  Qkn = markE0Tint(:,nn) .* Qkn;
  ret{1} = ret{1} + kron(spdiags(Qkn .* g.nuE0T(:,nn,1), 0,K,K), refEdgePhiIntPhiInt(:,:,nn));
  ret{2} = ret{2} + kron(spdiags(Qkn .* g.nuE0T(:,nn,2), 0,K,K), refEdgePhiIntPhiInt(:,:,nn));
end % for
end % function
\end{lstlisting}

\begin{lstlisting}[title={%
\code{ret = assembleMatElemDphiPhi(g, refElemDphiPhi)} assembles two matrices, each containing integrals of products of a basis function with a (spatial) derivative of a basis function. The matrices are returned in a~$2\times1$ \code{cell}~variable. This corresponds to the matrices~$\vecc{H}^m$, $m\in\{1,2\}$ according to Sec.~\ref{sec:assembly:globH}. The input argument~\code{refElemDphiPhi} stores the local matrices~$\hat{\vecc{H}}$ as defined in~\eqref{eq:hatH} and can be computed by~\code{integrateRefElemDphiPhi}.
}%
]
function ret = assembleMatElemDphiPhi(g, refElemDphiPhi)
K = g.numT; N = size(refElemDphiPhi, 1);
ret = cell(2, 1); ret{1} = sparse(K*N, K*N); ret{2} = sparse(K*N, K*N);
ret{1} = + kron(spdiags(g.B(:,2,2), 0,K,K), refElemDphiPhi(:,:,1)) ...
         - kron(spdiags(g.B(:,2,1), 0,K,K), refElemDphiPhi(:,:,2));
ret{2} = - kron(spdiags(g.B(:,1,2), 0,K,K), refElemDphiPhi(:,:,1)) ...
         + kron(spdiags(g.B(:,1,1), 0,K,K), refElemDphiPhi(:,:,2));
end % function
\end{lstlisting}

\begin{lstlisting}[title={%
\code{ret = assembleMatElemDphiPhiFuncDisc(g, refElemDphiPhiPhi, dataDisc)} assembles two matrices, each containing integrals of products of a basis function with a (spatial) derivative of a basis function and with a~discontinuous coefficient function whose coefficients are specified in~\code{dataDisc}. The matrices are returned in a~$2\times1$ \code{cell}~variable. This corresponds to the matrices~$\vecc{G}^m$, $m\in\{1,2\}$ according to Sec.~\ref{sec:assembly:globG}. The input argument~\code{refElemDphiPhiPhi} stores the local matrices~$\hat{\vecc{G}}$ (multidimensional array) as defined in~\eqref{eq:hatG} and can be computed by~\code{integrateRefElemDphiPhiPhi}. The coefficients of the projection of the algebraic diffusion coefficient $d$ into the broken polynomial space are stored in the input argument~\code{dataDisc} as explained in Sec.~\ref{sec:L2projection} and can be computed by~\code{projectFuncCont2dataDisc}.
}%
]
function ret = assembleMatElemDphiPhiFuncDisc(g, refElemDphiPhiPhi, dataDisc)
[K, N] = size(dataDisc);
ret = cell(2, 1);  ret{1} = sparse(K*N, K*N);  ret{2} = sparse(K*N, K*N);
for l = 1 : N
  ret{1} = ret{1} + kron(spdiags(dataDisc(:,l) .* g.B(:,2,2), 0,K,K), refElemDphiPhiPhi(:,:,l,1)) ...
                  - kron(spdiags(dataDisc(:,l) .* g.B(:,2,1), 0,K,K), refElemDphiPhiPhi(:,:,l,2));
  ret{2} = ret{2} - kron(spdiags(dataDisc(:,l) .* g.B(:,1,2), 0,K,K), refElemDphiPhiPhi(:,:,l,1)) ...
                  + kron(spdiags(dataDisc(:,l) .* g.B(:,1,1), 0,K,K), refElemDphiPhiPhi(:,:,l,2));
end % for
end % function
\end{lstlisting}

\begin{lstlisting}[title={%
\code{ret = assembleMatElemPhiPhi(g, hatM)} assembles a~matrix containing integrals of products of two basis functions, which corresponds to the global mass matrix $\vecc{M}$ according to Sec.~\ref{sec:assembly:globM}. The input argument~\code{hatM} stores the local matrices~$\hat{\vecc{M}}$ as defined in~\eqref{eq:hatM} and can be computed by~\code{integrateRefElemPhiPhi}.
}%
]
function ret = assembleMatElemPhiPhi(g, refElemPhiPhi)
K = g.numT;
ret = 2*kron(spdiags(g.areaT, 0, K, K), refElemPhiPhi);
end % function
\end{lstlisting}

\begin{lstlisting}[title={%
\code{ret = assembleVecEdgePhiIntFuncCont(g, markE0Tbdr, funcCont, N)} assembles a~vector containing integrals of products of a~basis function with a~continuous function. This corresponds to the contributions of Dirichlet boundaries~$\vec{K}_\mathrm{D}$ to the right-hand side of~\eqref{eq:spacediscretesystem:b} according to Sec.~\ref{sec:assembly:globKD}. Input arguments~\code{markE0Tbdr} and~\code{funcCont} are as described in~\code{assembleVecEdgePhiIntFuncContNu}.
}%
]
function ret = assembleVecEdgePhiIntFuncCont(g, markE0Tbdr, funcCont, N)
global gPhi1D
p = (sqrt(8*N+1)-3)/2;
qOrd = 2*p+1; [Q, W] = quadRule1D(qOrd);
Q2X1   = @(X1,X2) g.B(:,1,1)*X1 + g.B(:,1,2)*X2 + g.coordV0T(:,1,1)*ones(size(X1));
Q2X2   = @(X1,X2) g.B(:,2,1)*X1 + g.B(:,2,2)*X2 + g.coordV0T(:,1,2)*ones(size(X1));
ret = zeros(g.numT, N);
for n = 1 : 3
  [Q1, Q2] = gammaMap(n, Q);
  cDn = funcCont(Q2X1(Q1, Q2), Q2X2(Q1, Q2));
  for i = 1 : N
    ret(:,i) = ret(:,i) + markE0Tbdr(:,n) .* ( cDn * (W' .* gPhi1D{qOrd}(:,i,n)) );
  end % for
end % for
ret = reshape(ret', g.numT*N, 1);
end % function
\end{lstlisting}

\begin{lstlisting}[title={%
\code{ret = assembleVecEdgePhiIntFuncContNu(g, markE0Tbdr, funcCont, N)} assembles two vectors, each containing integrals of products of a basis function with a continuous function and with one of the components of edge normal. This corresponds to the contributions of Dirichlet boundaries $\vec{J}_\mathrm{D}^m$, $m\in\{1,2\}$ to the right-hand side of~\eqref{eq:spacediscretesystem:a} according to Sec.~\ref{sec:assembly:globJD}. The two vectors are returned in a $2\times1$ \code{cell}~variable. The input argument \code{markE0Tbdr} is a~$K\times3$~\code{logical} array marking for each triangle the local edges on the Dirichlet boundary. Assuming the Dirichlet boundary has the ID~$1$ then it can be computed as \code{markE0Tbdr = g.idE0T==1}. \code{funcCont} is a~\code{function_handle} for the algebraic representation of~$c_\mathrm{D}$~(cf.~\code{main.m}). Here, \lstinline$Q2X1$ and \lstinline$Q2X2$ carry out mappings $\vec{F}_k$ from the reference triangle $\hat{T}$ to physical triangles $T_k$ for all $k \in \{1,\ldots,K\}$ at once. This allows to evaluate $c_\mathrm{D}(t)$ at all quadrature points in the physical domain at the same time and, thus, vectorize the quadrature with respect to the triangle index.
}%
]
function ret = assembleVecEdgePhiIntFuncContNu(g, markE0Tbdr, funcCont, N)
global gPhi1D
K = g.numT; p = (sqrt(8*N+1)-3)/2;
qOrd = 2*p+1;  [Q, W] = quadRule1D(qOrd); 
Q2X1 = @(X1,X2) g.B(:,1,1)*X1 + g.B(:,1,2)*X2 + g.coordV0T(:,1,1)*ones(size(X1));
Q2X2 = @(X1,X2) g.B(:,2,1)*X1 + g.B(:,2,2)*X2 + g.coordV0T(:,1,2)*ones(size(X1));
ret = cell(2, 1);  ret{1} = zeros(K, N);  ret{2} = zeros(K, N);
for n = 1 : 3
  [Q1, Q2] = gammaMap(n, Q);
  cDn = funcCont(Q2X1(Q1, Q2), Q2X2(Q1, Q2));
  Jkn = markE0Tbdr(:,n) .* g.areaE0T(:,n);
  for i = 1 : N
    integral = cDn * ( W .* gPhi1D{qOrd}(:, i, n)' )';
    ret{1}(:,i) = ret{1}(:,i) + Jkn .* g.nuE0T(:,n,1) .* integral;
    ret{2}(:,i) = ret{2}(:,i) + Jkn .* g.nuE0T(:,n,2) .* integral;
  end % for
end % for
ret{1} = reshape(ret{1}',K*N,1);  ret{2} = reshape(ret{2}',K*N,1);
end % function
\end{lstlisting}
 
\begin{lstlisting}[title={%
\code{ret = assembleVecEdgePhiIntFuncDiscIntFuncCont(g, markE0Tbdr, dataDisc, funcCont)} assembles a~vector containing integrals of products of a~basis function with a~discontinuous coefficient fundtion and with a~continuous function. This corresponds to the contributions of Neumann boundaries~$\vec{K}_\mathrm{N}$ to the right-hand side of \eqref{eq:spacediscretesystem:b} according to Sec.~\ref{sec:assembly:globKN}. The input argument~\code{markE0Tbdr} is described in~\code{assembleVecEdgePhiIntFuncContNu}, \code{dataDisc} is the representation of the diffusion coefficient in the polynomial space (cf.~\code{projectFuncCont2DataDisc}), and \code{funcCont} is a~\code{function_handle} for the algebraic representation of~$g_\mathrm{N}$ (cf.~\code{main.m}).
}%
]
function ret = assembleVecEdgePhiIntFuncDiscIntFuncCont(g, markE0Tbdr, dataDisc, funcCont)
global gPhi1D
[K, N] = size(dataDisc);  p = (sqrt(8*N+1)-3)/2;
qOrd = 2*p+1;  [Q, W] = quadRule1D(qOrd);
Q2X1 = @(X1,X2) g.B(:,1,1)*X1 + g.B(:,1,2)*X2 + g.coordV0T(:,1,1)*ones(size(X1));
Q2X2 = @(X1,X2) g.B(:,2,1)*X1 + g.B(:,2,2)*X2 + g.coordV0T(:,1,2)*ones(size(X1));
ret = zeros(K, N);
for n = 1 : 3
  [Q1, Q2] = gammaMap(n, Q);
  funcAtQ = funcCont(Q2X1(Q1, Q2), Q2X2(Q1, Q2));
  Kkn = markE0Tbdr(:,n) .* g.areaE0T(:,n);
  for i = 1 : N
    for l = 1 : N
      integral = funcAtQ * ( W .* gPhi1D{qOrd}(:,i,n)' .* gPhi1D{qOrd}(:,l,n)' )';
      ret(:,i) = ret(:,i) + Kkn .* dataDisc(:,l) .* integral;
    end % for
  end % for
end % for
ret = reshape(ret',K*N,1);
end % function
\end{lstlisting}

\begin{lstlisting}[title={%
\code{computeBasesOnQuad(N)} evaluates the basis functions and their gradients in the quadrature points
on the reference triangle~$\hat{T}$ according to Sec.~\ref{sec:performance} for all required orders as described 
in Sec.~\ref{sec:quadrature}.
It computes the values of $\hat{\vphi}_i(\hat{\vec{x}})$, $\grad\hat{\vphi}_i(\hat{\vec{x}})$,
$\hat{\vphi}_i \circ \hat{\vec{\gamma}}(s)$, and 
$\hat{\vphi}_i \circ \mapEE_{n^-n^+} \circ \hat{\vec{\gamma}}(s)$ and stores them in~\code{global}
arrays \code{gPhi2D}, \code{gGradPhi2D}, \code{gPhi1D}, and \code{gThetaPhi1D}, respectively,
which are $(2p+1)\times1$~\code{cell} variables.
}%
]
function computeBasesOnQuad(N)
global gPhi2D gGradPhi2D gPhi1D gThetaPhi1D
p = (sqrt(8*N+1)-3)/2;
if p > 0, requiredOrders = [2*p, 2*p+1]; else requiredOrders = 1; end
gPhi2D = cell(max(requiredOrders),1);  gGradPhi2D  = cell(max(requiredOrders),1);
gPhi1D = cell(max(requiredOrders),1);  gThetaPhi1D = cell(max(requiredOrders),1);
for it = 1 : length(requiredOrders)
  ord = requiredOrders(it);
  [Q1, Q2, ~] = quadRule2D(ord);
  gPhi2D{ord}      = zeros(length(Q1), N);
  for i = 1 : N
    gPhi2D{ord}(:, i) = phi(i, Q1, Q2);
  end % for
  gGradPhi2D{ord}  = zeros(length(Q1), N, 2);
  for m = 1 : 2
    for i = 1 : N
      gGradPhi2D{ord}(:, i, m) = gradPhi(i, m, Q1, Q2);
    end % for
  end % for
  [Q, ~] = quadRule1D(ord);
  gPhi1D{ord} = zeros(length(Q), N, 3);
  for nn = 1 : 3
    [Q1, Q2] = gammaMap(nn, Q);
    for i = 1 : N
      gPhi1D{ord}(:, i, nn) = phi(i, Q1, Q2);
    end
    for np = 1 : 3
      [QP1,QP2] = theta(nn, np, Q1, Q2);
      for i = 1 : N
        gThetaPhi1D{ord}(:, i, nn, np) = phi(i, QP1, QP2);
      end % for
    end % for
  end % for
end % for
end % function
\end{lstlisting}

\begin{lstlisting}[title={%
\code{err = computeL2Error(g, dataDisc, funcCont, qOrd)} {computes the $L^2$-error according to Sec.~\ref{sec:discerror} using a~quadrature of order~\code{qOrd}, where \code{g}~is a~\code{struct} generated by~\code{generateGridData}, \code{dataDisc} is the $K\times N$~representation matrix of the approximate solution, \code{funcCont} is the algebraic formulation of the exact solution (see example below).
\\\emph{Example:} \code{N = 6;  h = 1/2;  qOrd = 5;}\\
\code{g = domainSquare(h);}  \code{computeBasesOnQuad(N);}  \code{hatM = integrateRefElemPhiPhi(N);}\\
\code{funcCont = @(X1, X2) sin(X1).*sin(X2);} \code{dataDisc  = projectFuncCont2DataDisc(g,funcCont,N,qOrd,hatM);}\\
\code{computeL2Error(g, dataDisc, funcCont, qOrd)}
}%
}%
]
function err = computeL2Error(g, dataDisc, funcCont, qOrd)
global gPhi2D
N = size(dataDisc, 2); qOrd = max(qOrd,1);
[Q1, Q2, W] = quadRule2D(qOrd);
R = length(W);
X1 = kron(g.B(:,1,1),Q1)+kron(g.B(:,1,2),Q2)+kron(g.coordV0T(:,1,1),ones(1,R));
X2 = kron(g.B(:,2,1),Q1)+kron(g.B(:,2,2),Q2)+kron(g.coordV0T(:,1,2),ones(1,R));
cExOnQuadPts = funcCont(X1, X2); % [K x R]
cApprxOnQuadPts = dataDisc*gPhi2D{qOrd}'; % [K x R] = [K x N] * [N x R]
err = sqrt(2*dot((cApprxOnQuadPts - cExOnQuadPts).^2 * W.', g.areaT)); 
end % function
\end{lstlisting}

\begin{lstlisting}[title={%
\code{g = domainCircle(h)} This method uses the grid generator Gmsh~\cite{Gmsh}. A~system call to Gmsh generates the ASCII file \code{domainCircle.mesh} based on geometry information of the domain~$\Omega$ stored in~\code{domainCircle.geo}.  The basic grid data is extracted from~\code{domainCircle.mesh} to call the routine~\code{generateGridData} (cf.~Sec.~\ref{sec:grid}) and to set the boundary IDs (from $1$ to~$4$).
}%
]
function g = domainCircle(h)
%% Generation of domainCircle.mesh using domainCircle.geo.
cmd = sprintf('gmsh -2 -format mesh -clscale %f -o "domainCircle.mesh" "domainCircle.geo"' , h);
system(cmd);
%% Extract data from domainCircle.mesh.
fid = fopen('domainCircle.mesh', 'r');
tline = fgets(fid);
while ischar(tline)
  if strfind(tline, 'Vertices')
    numV = fscanf(fid, '%d', [1, 1]);
    coordV = reshape(fscanf(fid, '%f'), 4, numV)';  coordV(:, 3:4) = [];
  end % if
  if strfind(tline, 'Edges')
    numEbdry = fscanf(fid, '%d', [1, 1]);
    tmp = reshape(fscanf(fid, '%f'), 3, numEbdry)';
    V0Ebdry = tmp(:, 1:2);  idEbdry = tmp(:, 3);
  end % if
  if strfind(tline, 'Triangles')
    numT = fscanf(fid, '%d', [1, 1]);
    V0T = reshape(fscanf(fid, '%f'), 4, numT)';  V0T(:, 4) = [];
  end % if
  tline = fgets(fid);
end % while
fclose(fid);
%% Generate lists and set boundary IDs.
g = generateGridData(coordV, V0T);
g.idE = zeros(g.numE, 1);
g.idE(g.V2E(sub2ind([g.numV,g.numV], V0Ebdry(:,1), V0Ebdry(:,2)))) = idEbdry;
g.idE0T = g.idE(g.E0T); % local edge IDs
end % function
\end{lstlisting}

\begin{lstlisting}[title={%
\code{domainCircle.geo} Geometry description of a~circle with center~$(0,0)$ and radius~$0.5$ serving as input for the grid generator Gmsh~\cite{Gmsh} which is called by the routine~\code{domainCircle}.  
}%
]
cx = 0.0;  cy = 0.0;  r = 0.5;
Point(0) = {cx, cy,   0.0, 1.0};  Point(1) = {cx-r, cy, 0.0, 1.0};  Point(2) = {cx, cy+r, 0.0, 1.0};
Point(3) = {cx+r, cy, 0.0, 1.0};  Point(4) = {cx, cy-r, 0.0, 1.0};
Circle(5) = {2, 0, 1};  Circle(6) = {3, 0, 2};  Circle(7) = {4, 0, 3};  Circle(8) = {1, 0, 4};
Line Loop(9) = {5, 6, 7, 8}; // boundary edge IDs will be 5 to 8
Plane Surface(10) = {9};
\end{lstlisting}

\begin{lstlisting}[title={%
\code{g = domainPolygon(X1, X2, h)} triangulates a polygonally bounded domain~$\Omega$, where~\code{X1} and~\code{X2} are lists of $x^1$ and $x^2$~coordinates of the corners. The output variable~\code{g} representing the triangulation~$\setT_h$ is of type~\code{struct} according to Sec.~\ref{sec:grid}. The boundary~IDs are set to the values produced by \code{initmesh}.  This method uses the \Matlab~grid generator~\code{initmesh}.
}%
]
function g = domainPolygon(X1, X2, h)
gd = [2; length(X1(:)); X1(:); X2(:)]; % geometry description
sf = 'polygon';                        % set formula
ns = double('polygon')';               % name space
[p, e, t] = initmesh(decsg(gd,sf,ns), 'Hmax', h);
g  = generateGridData(p', t(1:3, :)');
g.idE = zeros(g.numE, 1);
g.idE(g.V2E(sub2ind(size(g.V2E),e(1,:),e(2,:)))) = e(5,:);
g.idE0T = g.idE(g.E0T); % local edge IDs
end % function
\end{lstlisting}

\begin{lstlisting}[title={%
\code{g = domainSquare(h)} Friedrichs--Keller triangulation on the unit square. The input argument~\code{h} specifies the upper bound for the heights of the triangles (\emph{not} for the diameters).  The output variable~\code{g} representing the triangulation~$\setT_h$ is of type~\code{struct} according to Sec.~\ref{sec:grid}. The boundary~IDs are set from~$1$ to~$4$.
}%
]
function g = domainSquare(h)
dim = ceil(1/h); % number of edges per side of the unit square
h = 1/dim;
%% Build coordV.
[X, Y] = meshgrid(0:h:1);
Xlist = reshape(X, length(X)^2, 1);  Ylist = reshape(Y, length(X)^2, 1);
coordV = [Xlist, Ylist];
%% Build V0T.
pat1 = [1,dim+2,2]; % pattern of "lower-left" triangles
V0T1 = repmat(pat1, dim*(dim+1), 1) + repmat((0:dim*(dim+1)-1)', 1, 3);
V0T1(dim+1 : dim+1 : dim*(dim+1), :) = [];
pat2 = [dim+2,dim+3,2];
V0T2 = repmat(pat2, dim*(dim+1), 1) + repmat((0:dim*(dim+1)-1)', 1, 3);
V0T2(dim+1 : dim+1 : dim*(dim+1), :) = [];
%% Generate grid data and boundary IDs
g = generateGridData(coordV, [V0T1; V0T2]);
g.idE = zeros(g.numE, 1);
g.idE(g.baryE(:, 2) == 0) = 1; % south
g.idE(g.baryE(:, 1) == 1) = 2; % east
g.idE(g.baryE(:, 2) == 1) = 3; % north
g.idE(g.baryE(:, 1) == 0) = 4; % west
g.idE0T = g.idE(g.E0T); % local edge IDs
end % function
\end{lstlisting}

\begin{lstlisting}[title={%
\code{[X1, X2] = gammaMap(n, S)} {evaluates the parametrization of the \mbox{$n$th} edge~$\hat{E}_n$ of the reference triangle~$\hat{T}$ at parameter values specified by a~list of parameters~\code{S}, cf.~\eqref{eq:gammaMap}.}
}%
]
function [X1, X2] = gammaMap(n, S)
S = S(:)';
switch n
  case 1,  X1 = 1 - S;           X2 = S;
  case 2,  X1 = zeros(size(S));  X2 = 1 - S;
  case 3,  X1 = S;               X2 = zeros(size(S));
end % switch
end % function
\end{lstlisting}

\begin{lstlisting}[title={%
\code{g = generateGridData(coordV, V0T)}  assembles lists describing the geometric and topological properties of a~triangulation~$\setT_h$ according to Sec.~\ref{sec:grid} and stores them in the output variable~\code{g} of type~\code{struct}.  The input arguments are the array~\code{coordV} of dimension~$\card{\setV}\times 2$ that contains the $x^1$~and $x^2$~coordinates of the grid vertices (using a~\emph{global} index) and the array~\code{V0T} of dimension~$\card{\setT_h}\times 3$ storing the global vertex indices for each triangle with a counter-clockwise ordering.  A~usage example is found in Sec.~\ref{sec:visualization}.  Note that the lists~\code{g.idE} and~\code{g.idE0T} (cf.~Tab.~\ref{tab:lists}) storing the global and local edge indices are \emph{not} generated and have to be defined manually after calling~\code{generateGridData}.
}%
]
function g = generateGridData(coordV, V0T)
g.coordV = coordV;
g.V0T = V0T;
g.numT = size(g.V0T, 1);
g.numV = size(g.coordV, 1);
% The following implicitely defines the signs of the edges.
g.V2T  = sparse(g.V0T(:, [1 2 3 1 2 3 1 2 3]), g.V0T(:, [2 3 1 2 3 1 2 3 1]), ...
  [(1:g.numT)',zeros(g.numT,3),(1:g.numT)',zeros(g.numT,3),(1:g.numT)'],g.numV,g.numV);
% The following implicitely defines the edge numbers.
[r, c] = find(triu(g.V2T + g.V2T'));
g.V2E = sparse(r, c, 1 : size(r, 1), g.numV, g.numV);
g.V2E = g.V2E + g.V2E';
idxE = full(g.V2E(sub2ind([g.numV,g.numV],g.V0T(end:-1:1,[1,2,3]),g.V0T(end:-1:1,[2,3,1]))))';
g.V0E(idxE(:), 1) = reshape(g.V0T(end:-1:1, [1,2,3])', 3*g.numT, 1);
g.V0E(idxE(:), 2) = reshape(g.V0T(end:-1:1, [2,3,1])', 3*g.numT, 1);
g.T0E(idxE(:), 1) = reshape(full(g.V2T(sub2ind([g.numV,g.numV], ...
  g.V0T(end:-1:1,[1,2,3]), g.V0T(end:-1:1,[2,3,1]))))', 3*g.numT, 1);
g.T0E(idxE(:), 2) = reshape(full(g.V2T(sub2ind([g.numV,g.numV], ...
  g.V0T(end:-1:1,[2,3,1]), g.V0T(end:-1:1,[1,2,3]))))', 3*g.numT, 1);
g.numE = size(g.V0E, 1);
vecE = g.coordV(g.V0E(:, 2), :) - g.coordV(g.V0E(:, 1), :);
g.areaE = (vecE(:, 1).^2 + vecE(:, 2).^2).^(1/2);
g.nuE = vecE * [0,-1; 1,0] ./ g.areaE(:, [1, 1]);
g.areaT = ...
  ( g.coordV(g.V0T(:,1),1).*g.coordV(g.V0T(:,2),2) + g.coordV(g.V0T(:,2),1).*g.coordV(g.V0T(:,3),2) ...
  + g.coordV(g.V0T(:,3),1).*g.coordV(g.V0T(:,1),2) - g.coordV(g.V0T(:,1),1).*g.coordV(g.V0T(:,3),2) ...
  - g.coordV(g.V0T(:,2),1).*g.coordV(g.V0T(:,1),2) - g.coordV(g.V0T(:,3),1).*g.coordV(g.V0T(:,2),2) )/2;
g.baryT = (g.coordV(g.V0T(:,1),:)+g.coordV(g.V0T(:,2),:)+g.coordV(g.V0T(:,3),:))/3;
g.E0T = full(g.V2E(sub2ind([g.numV,g.numV],g.V0T(:,[2,3,1]),g.V0T(:,[3,1,2]))));
g.areaE0T = g.areaE(g.E0T);
sigE0T = 1-2*(bsxfun(@eq, reshape(g.T0E(g.E0T,2),g.numT,3),(1:g.numT)'));
g.baryE = 0.5 * (g.coordV(g.V0E(:, 1), :) + g.coordV(g.V0E(:, 2), :));
for n = 1 : 3
  for m = 1 : 2
    g.coordV0T(:, n, m) = g.coordV(g.V0T(:, n), m)';
    g.baryE0T(:, n, m) = g.baryE(g.E0T(:, n), m)';
    g.nuE0T(:, n, m) = g.nuE(g.E0T(:, n), m)'.* sigE0T(:, n)';
  end 
  Tn = sigE0T(:, n) == 1;  Tp = ~Tn;
  g.E0E(g.E0T(Tn, n), 1) = n;  g.E0E(g.E0T(Tp, n), 2) = n;
end % for
for m = 1 : 2
  g.B(:, m, 1) = g.coordV0T(:, 2, m) - g.coordV0T(:, 1, m);
  g.B(:, m, 2) = g.coordV0T(:, 3, m) - g.coordV0T(:, 1, m);
end % for
markEint = g.E0E(:, 2) ~= 0; % mark interior edges
g.markE0TE0T = cell(3, 3);
for nn = 1 : 3
  for np = 1 : 3
    g.markE0TE0T{nn,np} = sparse(g.numT, g.numT);
    markEn = g.E0E(:, 1) == nn;  markEp = g.E0E(:, 2) == np;
    idx = markEn & markEp & markEint;
    g.markE0TE0T{nn, np}(sub2ind([g.numT, g.numT], g.T0E(idx, 1), g.T0E(idx, 2))) = 1;
    markEn = g.E0E(:, 2) == nn;  markEp = g.E0E(:, 1) == np;
    idx = markEn & markEp & markEint;
    g.markE0TE0T{nn, np}(sub2ind([g.numT, g.numT], g.T0E(idx, 2), g.T0E(idx, 1))) = 1;
  end % for
end % for
end % function
\end{lstlisting}

\begin{lstlisting}[title={%
\code{ret = gradPhi(i, m, X1, X2)} {evaluates the $m$th component of the gradient of the $i$th basis function~$\hat{\vphi}_i$ on the reference triangle~$\hat{T}$ (cf.~Sec.~\ref{sec:basisfunctions}) at points specified by a~list of $\hat{x}^1$~coordinates~\code{X1} and $\hat{x}^2$~coordinates~\code{X2}.}
}%
]
function ret = gradPhi(i, m, X1, X2)
switch m
  case 1
    switch i
      case 1,  ret = zeros(size(X1));
      case 2,  ret = -6*ones(size(X1));
      case 3,  ret = -2*sqrt(3)*ones(size(X1));
      case 4,  ret = sqrt(6)*(20*X1 - 8);
      case 5,  ret = sqrt(3)*(10*X1 - 4);
      case 6,  ret = 6*sqrt(5)*(3*X1 + 4*X2 - 2);
      case 7,  ret = 2*sqrt(2)*(15+(-90+105*X1).*X1);
      case 8,  ret = 2*sqrt(6)*(13+(-66+63*X1).*X1+(-24+84*X1).*X2);
      case 9,  ret = 2*sqrt(10)*(9+(-30+21*X1).*X1+(-48+84*X1+42*X2).*X2);
      case 10, ret = 2*sqrt(14)*(3+(-6+3*X1).*X1+(-24+24*X1+30*X2).*X2);
      case 11, ret = sqrt(10)*(-24+(252+(-672+504*X1).*X1).*X1);
      case 12, ret = sqrt(30)*(-22+(210+(-504+336*X1).*X1).*X1+(42+(-336+504*X1).*X1).*X2);
      case 13, ret = 5*sqrt(2)*(-18+(138+(-264+144*X1).*X1).*X1 ...
                       +(102+(-624+648*X1).*X1+(-96+432*X1).*X2).*X2);
      case 14, ret = sqrt(70)*(-12+(60+(-84+36*X1).*X1).*X1 ...
                       +(132+(-456+324*X1).*X1+(-300+540*X1+180*X2).*X2).*X2);
      case 15, ret = 3*sqrt(10)*(-4+(12+(-12+4*X1).*X1).*X1 ...
                       +(60+(-120+60*X1).*X1+(-180+180*X1+140*X2).*X2).*X2);
    end
  case 2
    switch i
      case 1,  ret = zeros(size(X1));
      case 2,  ret = zeros(size(X1));
      case 3,  ret = -4*sqrt(3)*ones(size(X1));
      case 4,  ret = zeros(size(X1));
      case 5,  ret = 2*sqrt(3)*( -15*X2 + 6);
      case 6,  ret = 6*sqrt(5)*(4*X1 + 3*X2 - 2);
      case 7,  ret = zeros(size(X1));
      case 8,  ret = 2*sqrt(6)*(2+(-24+42*X1).*X1);
      case 9,  ret = 2*sqrt(10)*(6+(-48+42*X1).*X1+(-12+84*X1).*X2);
      case 10, ret = 2*sqrt(14)*(12+(-24+12*X1).*X1+(-60+60*X1+60*X2).*X2);
      case 11, ret = zeros(size(X1));
      case 12, ret = sqrt(30)*(-2+(42+(-168+168*X1).*X1).*X1);
      case 13, ret = 5*sqrt(2)*(-6+(102+(-312+216*X1).*X1).*X1 ...
                       +(12+(-192+432*X1).*X1).*X2);
      case 14, ret = sqrt(70)*(-12+(132+(-228+108*X1).*X1).*X1 ...
                       +(60+(-600+540*X1).*X1+(-60+540*X1).*X2).*X2);
      case 15, ret = 3*sqrt(10)*(-20+(60+(-60+20*X1).*X1).*X1 ...
                       +(180+(-360+180*X1).*X1+(-420+420*X1+280*X2).*X2).*X2);
    end % switch
end % switch
end % function
\end{lstlisting}

\begin{lstlisting}[title={%
\code{ret = integrateRefEdgePhiIntPhiExt(N)} computes a~multidimensional array of integrals over the edges of the reference triangle~$\hat{T}$, whose integrands consist of all permutations of two basis functions of which one belongs to a~neighboring element that is transformed using~$\mapEE$ (see~\eqref{eq:mapEE}). This corresponds to the local matrix~$\hat{\vecc{S}}^\mathrm{offdiag}$ as given in~\eqref{eq:hatSoffdiag}.
}%
]
function ret = integrateRefEdgePhiIntPhiExt(N)
global gPhi1D gThetaPhi1D
p = (sqrt(8*N+1)-3)/2;  qOrd = 2*p+1;  [~, W] = quadRule1D(qOrd); 
ret = zeros(N, N, 3, 3); % [N x N x 3 x 3]
for nn = 1 : 3
  for np = 1 : 3
    for i = 1 : N
      for j = 1 : N
        ret(i, j, nn, np) = sum( W' .* gPhi1D{qOrd}(:,i,nn) .* gThetaPhi1D{qOrd}(:,j,nn,np) );
      end % for
    end % for
  end % for
end % for
end % function
\end{lstlisting}

\begin{lstlisting}[title={%
\code{ret = integrateRefEdgePhiIntPhiExtPhiExt(N)} computes a~multidimensional array of integrals over the edges of the reference triangle~$\hat{T}$, whose integrands consist of all permutations of three basis functions of which two belong to a~neighboring element that are transformed using~$\mapEE$ (cf.~\eqref{eq:mapEE}). This corresponds to the local matrix~$\hat{\vecc{R}}^\mathrm{offdiag}$ as given in~\eqref{eq:hatRoffdiag}. 
}%
]
function ret = integrateRefEdgePhiIntPhiExtPhiExt(N)
global gPhi1D gThetaPhi1D
p = (sqrt(8*N+1)-3)/2;  qOrd = 2*p+1;  [~, W] = quadRule1D(qOrd);
ret = zeros(N,N,N,3,3); % [N x N x N x 3 x 3]
for nn = 1 : 3 % 3 edges
  for np = 1 : 3
    for l = 1 : N
      for i = 1 : N
        for j = 1 : N
          ret(i, j, l, nn,np) = sum( W'.*gPhi1D{qOrd}(:,i,nn) .* ...
            gThetaPhi1D{qOrd}(:,l,nn,np) .* gThetaPhi1D{qOrd}(:,j,nn,np) );
        end % for
      end % for
    end % for
  end % for
end % for
end
\end{lstlisting}

\begin{lstlisting}[title={%
\code{ret = integrateRefEdgePhiIntPhiInt(N)} computes a~multidimensional array of integrals over the edges of the reference triangle~$\hat{T}$, whose integrands consist of all permutations of two basis functions. This corresponds to the local matrix~$\hat{\vecc{S}}^\mathrm{diag}$ as given in~\eqref{eq:hatSdiag}. 
}%
]
function ret = integrateRefEdgePhiIntPhiInt(N)
global gPhi1D
p = (sqrt(8*N+1)-3)/2;  qOrd = 2*p+1;  [~, W] = quadRule1D(qOrd);
ret = zeros(N, N, 3); % [N x N x 3]
for n = 1 : 3 % 3 edges
  for i = 1 : N
    for j = 1 : N
      ret(i, j, n) = sum( W' .* gPhi1D{qOrd}(:,i,n) .* gPhi1D{qOrd}(:,j,n) );
    end % for
  end % for
end % for
end % function
\end{lstlisting}

\begin{lstlisting}[title={%
\code{ret = integrateRefEdgePhiIntPhiIntPhiInt(N)} computes a~multidimensional array of integrals over the edges of the reference triangle~$\hat{T}$, whose integrands consist of all permutations of three basis functions. This corresponds to the local matrix~$\hat{\vecc{R}}^\mathrm{diag}$ as given in~\eqref{eq:hatRdiag}.
}%
]
function ret = integrateRefEdgePhiIntPhiIntPhiInt(N)
global gPhi1D
p = (sqrt(8*N+1)-3)/2;  qOrd = max(2*p+1,1);  [~, W] = quadRule1D(qOrd);
ret = zeros(N, N, N, 3); % [N x N x N x 3]
for n = 1 : 3 % 3 edges
  for l = 1 : N % N basisfcts for D(t)
    for i = 1 : N
      for j = 1 : N
        ret(i,j,l,n) = sum(W' .* gPhi1D{qOrd}(:,i,n) .* gPhi1D{qOrd}(:,l,n) .* gPhi1D{qOrd}(:,j,n));
      end % for
    end % for
  end % for
end % for
end % function
\end{lstlisting}

\begin{lstlisting}[title={%
\code{ret = integrateRefElemDphiPhi(N)} computes a~multidimensional array of integrals on the reference triangle~$\hat{T}$, whose integrands consist of all permutations of a basis function with one of the (spatial) derivatives of a basis function. This corresponds to the local matrix~$\hat{\vecc{H}}$ as given in~\eqref{eq:hatH}.
}%
]
function ret = integrateRefElemDphiPhi(N)
global gPhi2D gGradPhi2D
ret = zeros(N, N, 2); % [ N x N x 2]
if N > 1 % p > 0
  p = (sqrt(8*N+1)-3)/2;  qOrd = max(2*p, 1);  [~,~,W] = quadRule2D(qOrd);
  for i = 1 : N
    for j = 1 : N
      for m = 1 : 2
        ret(i, j, m) = sum( W' .* gGradPhi2D{qOrd}(:,i,m) .* gPhi2D{qOrd}(:,j) );
      end % for
    end % for
  end % for
end % function
\end{lstlisting}

\begin{lstlisting}[title={%
\code{ret = integrateRefElemDphiPhiPhi(N)} computes a~multidimensional array of integrals on the reference triangle~$\hat{T}$, whose integrands consist of all permutations of two basis functions with one of the (spatial) derivatives of a basis function. This corresponds to the local matrix~$\hat{\vecc{G}}$ as given in~\eqref{eq:hatG}.
}%
]
function ret = integrateRefElemDphiPhiPhi(N)
global gPhi2D gGradPhi2D
ret = zeros(N, N, N, 2); % [N x N x N x 2]
if N > 1 % p > 0
  p = (sqrt(8*N+1)-3)/2;  qOrd = max(2*p, 1);  [~,~,W] = quadRule2D(qOrd);
  for i = 1 : N
    for j = 1 : N
      for l = 1 : N
        for m = 1 : 2
          ret(i,j,l,m) = sum( W' .* gGradPhi2D{qOrd}(:,i,m) .* gPhi2D{qOrd}(:,j) .* gPhi2D{qOrd}(:,l) );
        end % for
      end % for
    end % for
  end % for
end % function
\end{lstlisting}

\begin{lstlisting}[title={%
\code{ret = integrateRefElemPhiPhi(N)} computes a~multidimensional array of integrals on the reference triangle~$\hat{T}$, whose integrands consist of all permutations of two basis functions. This corresponds to the local matrix~$\hat{\vecc{M}}$ as given in~\eqref{eq:hatM}.
}%
]
function ret = integrateRefElemPhiPhi(N)
global gPhi2D
p = (sqrt(8*N+1)-3)/2;  qOrd = max(2*p, 1);  [~,~,W] = quadRule2D(qOrd);
ret = zeros(N); % [N x N]
for i = 1 : N
  for j = 1 : N
    ret(i, j) = sum( W' .* gPhi2D{qOrd}(:, i) .* gPhi2D{qOrd}(:, j) );
  end % for
end % for
end % function
\end{lstlisting}

\begin{lstlisting}[title={%
\code{main.m} {This is the main script to solve~\eqref{eq:diffusion} which can be used as a template for further modifications. Modifiable parameters are found in Lines~4--10, the problem data (initial condition, diffusion coefficient, right hand side and boundary data) is specified in Lines~16--20.}
}%
]
function main()
more off % disable paging of output
%% Parameters.
hmax        = 1/8;    % maximum edge length of triangle
p           = 2;      % local polynomial degree
tEnd        = pi;     % end time
numSteps    = 20;     % number of time steps
isVisGrid   = true;  % visualization of grid
isVisSol    = true;   % visualization of solution
eta         = 1;      % penalty parameter (eta>0)
%% Parameter check.
assert(p >= 0 && p <= 4, 'Polynomial order must be zero to four.')
assert(hmax > 0        , 'Maximum edge length must be positive.' )
assert(numSteps > 0    , 'Number of time steps must be positive.')
%% Coefficients and boundary data.
c0     = @(x1,x2) sin(x1).*cos(x2);
dCont  = @(t,x1,x2) (x1<3/4&x1>1/4&x2<3/4&x2>1/4) + 0.01;
fCont  = @(t,x1,x2) 0.1*t*(x1==x1);
cDCont = @(t,x1,x2) sin(2*pi*x2 + t);
gNCont = @(t,x1,x2) x2;
%% Triangulation.
g = domainSquare(hmax);
if isVisGrid,  visualizeGrid(g);  end
%% Globally constant parameters.
K           = g.numT;                      % number of triangles
N           = nchoosek(p + 2, p);          % number of local DOFs
tau         = tEnd/numSteps;               % time step size
markE0Tint  = g.idE0T == 0;                % [K x 3] mark local edges that are interior
markE0TbdrN = g.idE0T == 1 | g.idE0T == 3; % [K x 3] mark local edges on the Neumann boundary
markE0TbdrD = ~(markE0Tint | markE0TbdrN); % [K x 3] mark local edges on the Dirichlet boundary
%% Configuration output.
fprintf('Computing with polynomial order %d (%d local DOFs) on %d triangles.\n', p, N, K)
%% Lookup table for basis function.
computeBasesOnQuad(N);
%% Computation of matrices on the reference triangle.
hatM        = integrateRefElemPhiPhi(N);
hatG        = integrateRefElemDphiPhiPhi(N);
hatH        = integrateRefElemDphiPhi(N);
hatRdiag    = integrateRefEdgePhiIntPhiIntPhiInt(N);
hatRoffdiag = integrateRefEdgePhiIntPhiExtPhiExt(N);
hatSdiag    = integrateRefEdgePhiIntPhiInt(N);
hatSoffdiag = integrateRefEdgePhiIntPhiExt(N);
%% Assembly of time-independent global matrices.
globM  = assembleMatElemPhiPhi(g, hatM);
globH  = assembleMatElemDphiPhi(g, hatH);
globQ  = assembleMatEdgePhiPhiNu(g, markE0Tint, hatSdiag, hatSoffdiag);
globQN = assembleMatEdgePhiIntPhiIntNu(g, markE0TbdrN, hatSdiag);
globS  = eta * assembleMatEdgePhiPhi(g, markE0Tint, hatSdiag, hatSoffdiag);
globSD = eta * assembleMatEdgePhiIntPhiInt(g, markE0TbdrD, hatSdiag);
sysW = [ sparse(2*K*N,3*K*N) ; sparse(K*N,2*K*N), globM ];
%% Initial data.
cDisc = projectFuncCont2DataDisc(g, c0, 2*p, hatM);
sysY = [ zeros(2*K*N,1) ; reshape(cDisc', K*N, 1) ];
%% Time stepping.
fprintf('Starting time integration from 0 to %g using time step size %g (%d steps).\n', tEnd, tau, numSteps)
for nStep = 1 : numSteps
  t = nStep * tau;
  %% L2-projections of algebraic coefficients.
  dDisc = projectFuncCont2DataDisc(g, @(x1,x2) dCont(t,x1,x2), 2*p, hatM);
  fDisc = projectFuncCont2DataDisc(g, @(x1,x2) fCont(t,x1,x2), 2*p, hatM);
  %% Assembly of time-dependent global matrices.
  globG = assembleMatElemDphiPhiFuncDisc(g, hatG, dDisc);
  globR = assembleMatEdgePhiPhiFuncDiscNu(g, markE0Tint, hatRdiag, hatRoffdiag, dDisc);
  %% Assembly of Dirichlet boundary contributions.
  globRD = assembleMatEdgePhiIntPhiIntFuncDiscIntNu(g, markE0TbdrD, hatRdiag, dDisc);
  globJD = assembleVecEdgePhiIntFuncContNu(g, markE0TbdrD, @(x1,x2) cDCont(t,x1,x2), N);
  globKD = eta * assembleVecEdgePhiIntFuncCont(g, markE0TbdrD, @(x1,x2) cDCont(t,x1,x2), N);
  %% Assembly of Neumann boundary contributions.
  globKN = assembleVecEdgePhiIntFuncDiscIntFuncCont(g, markE0TbdrN, dDisc, @(x1,x2) gNCont(t,x1,x2));
  %% Assembly of the source contribution.
  globL = globM*reshape(fDisc', K*N, 1);
  %% Building and solving the system.
  sysA = [                        globM,              sparse(K*N,K*N), -globH{1}+globQ{1}+globQN{1};
                        sparse(K*N,K*N),                        globM, -globH{2}+globQ{2}+globQN{2};
           -globG{1}+globR{1}+globRD{1}, -globG{2}+globR{2}+globRD{2},                globS+globSD];
  sysV = [                   -globJD{1};                   -globJD{2};        globKD-globKN+globL];
  sysY = (sysW + tau*sysA) \ (sysW*sysY + tau*sysV);
  %% Visualization
  if isVisSol
    cDisc = reshape(sysY(2*K*N+1 : 3*K*N), N, K)';
    cLagr = projectDataDisc2DataLagr(cDisc);
    visualizeDataLagr(g, cLagr, 'c_h', 'solution', nStep)
  end % if
end % for
fprintf('Done.\n')
end % function
\end{lstlisting}

\begin{lstlisting}[title={%
\code{ret = phi(i, X1, X2)} {evaluates the $i$th basis function~$\hat{\vphi}_i$ on the reference triangle~$\hat{T}$ (cf.~Sec.~\ref{sec:basisfunctions}) at points specified by a~list of $\hat{x}^1$~coordinates~\code{X1} and $\hat{x}^2$~coordinates~\code{X2}.}
}%
]
function ret = phi(i, X1, X2)
switch i
  case 1,  ret = sqrt(2)*ones(size(X1));
  case 2,  ret = 2 - 6*X1;
  case 3,  ret = 2*sqrt(3)*(1 - X1 - 2*X2);
  case 4,  ret = sqrt(6)*((10*X1 - 8).*X1 + 1);
  case 5,  ret = sqrt(3)*((5*X1 - 4).*X1 + (-15*X2 + 12).*X2 - 1);
  case 6,  ret = 3*sqrt(5)*((3*X1 + 8*X2 - 4).*X1 + (3*X2 - 4).*X2 + 1);
  case 7,  ret = 2*sqrt(2)*(-1+(15+(-45+35*X1).*X1).*X1); 
  case 8,  ret = 2*sqrt(6)*(-1+(13+(-33+21*X1).*X1).*X1+(2+(-24+42*X1).*X1).*X2);
  case 9,  ret = 2*sqrt(10)*(-1+(9+(-15+7*X1).*X1).*X1+(6+(-48+42*X1).*X1+(-6+42*X1).*X2).*X2);
  case 10, ret = 2*sqrt(14)*(-1+(3+(-3+X1).*X1).*X1+(12+(-24+12*X1).*X1+(-30+30*X1+20*X2).*X2).*X2);
  case 11, ret = sqrt(10)*(1+(-24+(126+(-224+126*X1).*X1).*X1).*X1);
  case 12, ret = sqrt(30)*(1+(-22+(105+(-168+84*X1).*X1).*X1).*X1+(-2+(42+(-168+168*X1).*X1).*X1).*X2);
  case 13, ret = 5*sqrt(2)*(1+(-18+(69+(-88+36*X1).*X1).*X1).*X1+(-6+(102+(-312+216*X1).*X1).*X1 ...
                   +(6+(-96+216*X1).*X1).*X2).*X2);
  case 14, ret = sqrt(70)*(1+(-12+(30+(-28+9*X1).*X1).*X1).*X1+(-12+(132+(-228+108*X1).*X1).*X1 ...
                   +(30+(-300+270*X1).*X1+(-20+180*X1).*X2).*X2).*X2);
  case 15, ret = 3*sqrt(10)*(1+(-4+(6+(-4+X1).*X1).*X1).*X1+(-20+(60+(-60+20*X1).*X1).*X1 ...
                   +(90+(-180+90*X1).*X1+(-140+140*X1+70*X2).*X2).*X2).*X2);
end % switch
end % function
\end{lstlisting}

\begin{lstlisting}[title={%
\code{dataLagr = projectDataDisc2DataLagr(dataDisc)}
converts the representation matrix in the DG\,/\,modal basis to the respective representation matrix in a~Lagrange\,/\,nodal basis, both of size~\mbox{$K\times N$} for $p\in\{0,1,2\}$ (cf.~Sec.~\ref{sec:visualization}). For~$p>2$ the output argument has the size~\mbox{$K\times 6$} as \code{visualizeDataLagr} can visualize up to elementwise quadratics only. In this routine, the local basis functions~$\hat{\vphi}_i$ are sampled at the Lagrange nodes on the reference triangle~$\hat{T}$, whose $\hat{x}^1$ and $\hat{x}^2$~coordinates are stored in the variables~\code{L1} and~\code{L2}.
}%
]
function dataLagr = projectDataDisc2DataLagr(dataDisc)
[K, N] = size(dataDisc);
switch N
  case 1,    L1 = 1/3;                     L2 = 1/3;                    % locally constant
  case 3,    L1 = [0, 1, 0];               L2 = [0, 0, 1];              % locally linear
  otherwise, L1 = [0, 1, 0, 1/2, 0, 1/2];  L2 = [0, 0, 1, 1/2, 1/2, 0]; % locally quadratic
end % switch
dataLagr = zeros(K, length(L1));
for i = 1 : N
  dataLagr = dataLagr + dataDisc(:, i) * phi(i, L1, L2);
end % for
end % function
\end{lstlisting}

\begin{lstlisting}[title={%
\code{dataDisc = projectFuncCont2DataDisc(g, funcCont, R, refElemPhiPhi)}
computes the representation matrix~\code{dataDisc} of an~algebraic function in the DG\,/\,modal basis by performing the $L^2$-projection described in~Sec.~\ref{sec:L2projection}.  The non-obvious input arguments are as follows:
\code{funcCont} is of type~\code{function_handle} taking rows of coordinates and \code{refElemPhiPhi} the matrix~\mbox{$\hat{\vecc{M}}$} computed by~\code{integrateRefElemPhiPhi}.  
}%
]
function dataDisc = projectFuncCont2DataDisc(g, funcCont, ord, refElemPhiPhi)
global gPhi2D
ord = max(ord,1);  [Q1, Q2, W] = quadRule2D(ord);
K = g.numT; N = size(refElemPhiPhi, 1);
F1 = @(X1, X2) g.B(:,1,1)*X1 + g.B(:,1,2)*X2 + g.coordV0T(:,1,1)*ones(size(X1));
F2 = @(X1, X2) g.B(:,2,1)*X1 + g.B(:,2,2)*X2 + g.coordV0T(:,1,2)*ones(size(X1));
rhs = funcCont(F1(Q1, Q2), F2(Q1, Q2)) * (repmat(W.', 1, N) .* gPhi2D{ord});
dataDisc = rhs / refElemPhiPhi;
end % function
\end{lstlisting}

\begin{lstlisting}[title={%
\code{[Q, W] = quadRule1D(qOrd)}
returns a~list of Gauss--Legendre quadrature points~\code{Q} within the interval~\mbox{$(0,1)$} and associated (positive) weights~\code{W}.  The quadrature rule is exact for polynomials of order~\code{qOrd}.  The length of the  interval~\mbox{$[0,1]$} is incorporated in the weights.  A~rule with \mbox{$n$}~points is exact for polynomials up to order~\mbox{$2n-1$}.
\\\emph{Example:} \code{f = @(s) s.^2; [Q, W] = quadRule1D(2); dot(W, f(Q))}%
}%
]
function [Q, W] = quadRule1D(qOrd)
switch qOrd
  case {0, 1} % R = 1, number of quadrature points
    Q = 0;
    W = 2;
  case {2, 3} % R = 2
    Q = sqrt(1/3)*[-1, 1];
    W = [1, 1];
  case {4, 5} % R = 3
    Q = sqrt(3/5)*[-1, 0, 1];
    W = 1/9*[5, 8, 5];
  case {6, 7} % R = 4
    Q = [-1,-1,1,1].*sqrt(3/7+[1,-1,-1,1]*2/7*sqrt(6/5));
    W = 1/36*(18 + sqrt(30)*[-1,1,1,-1]);
  case {8, 9} % R = 5
    Q = [-1,-1,0,1,1].*sqrt(5+[2,-2,0,-2,2]*sqrt(10/7))/3;
    W = 1/900*(322+13*sqrt(70)*[-1,1,0,1,-1]+[0,0,190,0,0]);
  case {10, 11} % R = 6
    Q = [ 0.6612093864662645, -0.6612093864662645, -0.2386191860831969, ...
          0.2386191860831969, -0.9324695142031521,  0.9324695142031521];
    W = [ 0.3607615730481386,  0.3607615730481386,  0.4679139345726910, ...
          0.4679139345726910,  0.1713244923791704,  0.171324492379170];
  case {12, 13} % R = 7
    Q = [ 0.0000000000000000,  0.4058451513773972, -0.4058451513773972, ...
         -0.7415311855993945,  0.7415311855993945, -0.9491079123427585, ...
          0.9491079123427585];
    W = [ 0.4179591836734694,  0.3818300505051189,  0.3818300505051189, ...
          0.2797053914892766,  0.2797053914892766,  0.1294849661688697, ...
          0.1294849661688697];  
  case {14, 15} % R = 8
    Q = [-0.1834346424956498,  0.1834346424956498, -0.5255324099163290, ...
          0.5255324099163290, -0.7966664774136267,  0.7966664774136267, ...
         -0.9602898564975363,  0.9602898564975363]; 
    W = [ 0.3626837833783620,  0.3626837833783620,  0.3137066458778873, ...
          0.3137066458778873,  0.2223810344533745,  0.2223810344533745, ...
          0.1012285362903763,  0.1012285362903763];
  case {16, 17} % R = 9
    Q = [ 0.0000000000000000, -0.8360311073266358,  0.8360311073266358, ...
         -0.9681602395076261,  0.9681602395076261, -0.3242534234038089, ...
          0.3242534234038089, -0.6133714327005904,  0.6133714327005904];
    W = [ 0.3302393550012598,  0.1806481606948574,  0.1806481606948574, ...
          0.0812743883615744,  0.0812743883615744,  0.3123470770400029, ...
          0.3123470770400029,  0.2606106964029354,  0.2606106964029354];
end % switch
Q = (Q + 1)/2;  W = W/2; % transformation [-1; 1] -> [0, 1]
end % function
\end{lstlisting}

\begin{lstlisting}[title={%
\code{[Q1, Q2, W] = quadRule2D(qOrd)}
returns quadrature points~$\hat{\vec{q}}_r = \transpose{[\hat{q}_r^1, \hat{q}_r^2]}$ within the reference triangle~\mbox{$\hat{T}$} in lists of~\mbox{$\hat{x}^1$} and \mbox{$\hat{x}^2$}~coordinates~\code{Q1} and~\code{Q2}, respectively,
and the associated weights~\code{Q} (cf.~Sec.~\ref{sec:quadrature}).  The quadrature rule is exact for polynomials of order~\code{qOrd} (cf.~\cite{Stroud1971} for orders~$1,2,5$, \cite{Hillion1977} for order~$3$, and \cite{Cowper1973} for order~$4,6$). If \code{qOrd} is greater than~$6$, we call the third party function~\code{triquad}~\cite{triquad} that uses Gaussian quadrature points on a~square which is collapsed to a~triangle. The area~\mbox{$1/2$} of the reference triangle~$\hat{T}$ is incorporated in the weights such that the integral over one is~\mbox{$1/2$}.%
\\\emph{Example:} \code{[Q1, Q2, W] = quadRule2D(2); N = 3; M = eye(N);}
\\ \code{for i=1:N,for j=1:N,M(i,j)=sum(W.*phi(i,Q1,Q2).*phi(j,Q1,Q2));end,end}%
}%
]
function [Q1, Q2, W] = quadRule2D(qOrd)
switch qOrd
  case {0, 1} % R = 1
    Q1 = 1/3;              Q2 = 1/3;              W  = 1/2;
  case 2 % R = 3
    Q1 = [1/6, 2/3, 1/6];  Q2 = [1/6, 1/6, 2/3];  W  = [1/6, 1/6, 1/6];
  case 3 % R = 4
    Q1 = [0.666390246, 0.178558728, 0.280019915, 0.075031109];
    Q2 = [0.178558728, 0.666390246, 0.075031109, 0.280019915];
    W  = [0.159020691, 0.159020691, 0.090979309, 0.090979309];
  case 4 % R = 6
    Q1 = [0.445948490915965, 0.108103018168070, 0.445948490915965, ...
          0.091576213509771, 0.816847572980458, 0.091576213509771];
    Q2 = [0.108103018168070, 0.445948490915965, 0.445948490915965, ...
          0.816847572980458, 0.091576213509771, 0.091576213509771];
    W  = [0.111690794839005, 0.111690794839005, 0.111690794839005, ...
          0.054975871827661, 0.054975871827661, 0.054975871827661];
  case 5 % R = 7
    a1 = (6-sqrt(15))/21;     a2 = (6+sqrt(15))/21;
    w1 = (155-sqrt(15))/2400; w2 = (155+sqrt(15))/2400;
    Q1 = [1/3,     a1, 1-2*a1,     a1,     a2, 1-2*a2,     a2];
    Q2 = [1/3, 1-2*a1,     a1,     a1, 1-2*a2,     a2,     a2];
    W  = [9/80,    w1,     w1,     w1,     w2,     w2,     w2];
  case 6 % R = 12
    Q1 = [0.063089014491502, 0.873821971016996, 0.063089014491502, ...
          0.249286745170910, 0.501426509658179, 0.249286745170910, ...
          0.310352451033785, 0.053145049844816, 0.636502499121399, ...
          0.053145049844816, 0.636502499121399, 0.310352451033785];
    Q2 = [0.063089014491502, 0.063089014491502, 0.873821971016996, ...
          0.249286745170910, 0.249286745170910, 0.501426509658179, ...
          0.053145049844816, 0.310352451033785, 0.053145049844816, ...
          0.636502499121399, 0.310352451033785, 0.636502499121399];
    W  = [0.025422453185103, 0.025422453185103, 0.025422453185103, ...
          0.058393137863189, 0.058393137863189, 0.058393137863189, ...
          0.041425537809187, 0.041425537809187, 0.041425537809187, ...
          0.041425537809187, 0.041425537809187, 0.041425537809187];
  otherwise % use Gauss quadrature points on square
    [X,Y,Wx,Wy] = triquad(qOrd, [0 0; 0 1; 1 0]); % third party function, see references
    Q1 = X(:).';
    Q2 = Y(:).';
    W = Wx * Wy.';  W = W(:).';
end % switch
end % function
\end{lstlisting}

\begin{lstlisting}[title={%
\code{[XP1, XP2] = theta(nn, np, X1, X2)}
returns the mapped points from $n^-$th edge to the $n^+$th edge of the reference triangle~$\hat{T}$, cf.~\eqref{eq:mapEE}.
}%
]
function [XP1, XP2] = theta(nn, np, X1, X2)
switch nn
  case 1
    switch np
      case 1,  XP1 = 1-X1;  XP2 = 1-X2;
      case 2,  XP1 = zeros(size(X1));  XP2 = X2;
      case 3,  XP1 = X1;  XP2 = zeros(size(X1));
    end % switch
  case 2
    switch np
      case 1,  XP1 = 1-X2;  XP2 = X2;
      case 2,  XP1 = zeros(size(X1));  XP2 = 1-X2;
      case 3,  XP1 = X2;  XP2 = zeros(size(X1));
    end % switch
  case 3
    switch np
      case 1,  XP1 = X1;  XP2 = 1-X1;
      case 2,  XP1 = zeros(size(X1));  XP2 = X1;
      case 3,  XP1 = 1-X1;  XP2 = zeros(size(X1));
    end % switch
end % switch
end % function
\end{lstlisting}

\begin{lstlisting}[title={%
\code{visualizeDataLagr(g, dataLagr, varName, fileName, tLvl)} 
writes a~\code{.vtu}~file for the visualization of a~discrete quantity in~\mbox{$\IP_p(\setT_h)$, $p\in\{0,1,2\}$} defined on the triangulation~\code{g} according to~\code{generateGridData}.  The name of the generated file is \code{fileName.tLvl.vtu}, where \code{tLvl} stands for time level.  The name of the quantity within the file is specified by~\code{varName}.  The argument~\code{dataLagr} should be a~list of dimension~\mbox{$K\times N$} containing the Lagragian representation of the quantity.  The \mbox{$k$th}~row of~\code{dataLagr} has to hold the value on~$T_k$ for~$\IP_0(\setT_h)$, the values on the vertices of~$T_k$ for~$\IP_1(\setT_h)$, and the values on the vertices and on the edge barycenters of~$T_k$ for~$\IP_2(\setT_h)$. Note that we treat functions of~$\IP_0(\setT_h)$ as if they were in~$\IP_1(\setT_h)$ as we assign the constant value on~$T_k$ to the vertices~$\vec{a}_{k1},\vec{a}_{k2},\vec{a}_{k3}$. An~alternative was using \code{CellData} instead of~\code{PointData}. A~usage example is found in Sec.~\ref{sec:visualization}.
}%
]
function visualizeDataLagr(g, dataLagr, varName, fileName, tLvl)
[K, N] = size(dataLagr);
%% Open file.
fileName = [fileName, '.', num2str(tLvl), '.vtu'];
file     = fopen(fileName, 'wt'); % if this file exists, then overwrite
%% Header.
fprintf(file, '<?xml version="1.0"?>\n');
fprintf(file, '<VTKFile type="UnstructuredGrid" version="0.1" byte_order="LittleEndian" compressor="vtkZLibDataCompressor">\n');
fprintf(file, '  <UnstructuredGrid>\n');
%% Points and cells.
switch N
  case {1, 3}
    P1          = reshape(g.coordV0T(:, :, 1)', 3*K, 1);
    P2          = reshape(g.coordV0T(:, :, 2)', 3*K, 1);
    numP        = 3; % number of local points
    id          = 5; % vtk ID for linear polynomials
  case 6
    P1          = reshape([g.coordV0T(:,:,1), g.baryE0T(:,[3,1,2],1)]',6*K,1);
    P2          = reshape([g.coordV0T(:,:,2), g.baryE0T(:,[3,1,2],2)]',6*K,1);
    numP        = 6; % number of local points
    id          = 22; % vtk ID for quadratic polynomials
end % switch
fprintf(file, '    <Piece NumberOfPoints="%d" NumberOfCells="%d">\n',K*numP,K);
fprintf(file, '      <Points>\n');
fprintf(file, '        <DataArray type="Float32" NumberOfComponents="3" format="ascii">\n');
fprintf(file, '          %.3e %.3e %.3e\n',  [P1, P2, zeros(numP*K, 1)]');
fprintf(file, '        </DataArray>\n');
fprintf(file, '      </Points>\n');
fprintf(file, '      <Cells>\n');
fprintf(file, '        <DataArray type="Int32" Name="connectivity" format="ascii">\n');
fprintf(file, '           '); fprintf(file,'%d ', 0:K*numP-1);
fprintf(file, '\n        </DataArray>\n');
fprintf(file, '        <DataArray type="Int32" Name="offsets" format="ascii">\n');
fprintf(file, '           %d\n', numP:numP:numP*K);
fprintf(file, '        </DataArray>\n');
fprintf(file, '        <DataArray type="UInt8" Name="types" format="ascii">\n');
fprintf(file, '           %d\n', id*ones(K, 1));
fprintf(file, '        </DataArray>\n');
fprintf(file, '      </Cells>\n');
%% Data.
switch N
  case 1 % locally constant
    dataLagr = kron(dataLagr, [1;1;1])';
  case 3 % locally quadratic
    dataLagr = reshape(dataLagr', 1, K*N);
  case 6 % locally quadratic (permutation of local edge indices due to vtk format)
    dataLagr = reshape(dataLagr(:, [1,2,3,6,4,5])', 1, K*N);
end % switch
fprintf(file, '      <PointData Scalars="%s">\n', varName);
fprintf(file, '        <DataArray type="Float32" Name="%s" NumberOfComponents="1" format="ascii">\n', varName);
fprintf(file, '          %.3e\n', dataLagr);
fprintf(file, '        </DataArray>\n');
fprintf(file, '      </PointData>\n');
%% Footer.
fprintf(file, '    </Piece>\n');
fprintf(file, '  </UnstructuredGrid>\n');
fprintf(file, '</VTKFile>\n');
%% Close file.
fclose(file);
disp(['Data written to ' fileName])
end % function
\end{lstlisting}

\begin{lstlisting}[title={%
\code{visualizeGrid(g)} visualizes the triangulation~$\setT_h$ along with global and local indices and edge normals, cf.~Fig.~\ref{fig:visualization}.  The input argument~\code{g} is the output of the routine~\code{generateGridData}.
}%
]
function visualizeGrid(g)
figure('Color', [1, 1, 1]); % white background
hold('on'),  axis('off')
daspect([1, 1, 1]) % adjust aspect ration, requires Octave >= 3.8
textarray = @(x1,x2,s) arrayfun(@(a,b,c) text(a,b,int2str(c),'Color','blue'), x1, x2, s);
%% Triangle boundaries.
trisurf(g.V0T,g.coordV(:,1),g.coordV(:,2),zeros(g.numV,1), 'facecolor', 'none');
%% Local edge numbers.
w = [1/12, 11/24, 11/24; 11/24, 1/12, 11/24; 11/24, 11/24, 1/12];
for kE = 1 : 3
  textarray(reshape(g.coordV(g.V0T,1),g.numT,3)*w(:,kE), ...
    reshape(g.coordV(g.V0T,2),g.numT,3)*w(:,kE), kE*ones(g.numT, 1))
end % for
%% Global vertex numbers.
textarray(g.coordV(:,1), g.coordV(:,2), (1:g.numV)');
%% Local vertex numbers.
w = [5/6, 1/12, 1/12; 1/12, 5/6, 1/12; 1/12, 1/12, 5/6];
for kV = 1 : 3
  textarray(reshape(g.coordV(g.V0T,1),g.numT,3)*w(:,kV), ...
    reshape(g.coordV(g.V0T,2),g.numT,3)*w(:,kV), kV*ones(g.numT, 1))
end % for
%% Global edge numbers.
textarray(g.baryE(:,1), g.baryE(:,2), (1:g.numE)');
%% Triangle numbers.
textarray(g.baryT(:,1), g.baryT(:,2), (1:g.numT)');
%% Edge IDs.
markEext = g.idE ~= 0; % mark boundary edges
textarray(g.baryE(markEext,1) + g.nuE(markEext,1).*g.areaE(markEext)/8, ...
  g.baryE(markEext,2) + g.nuE(markEext,2).*g.areaE(markEext)/8, g.idE(markEext))
end % function
\end{lstlisting}

\section{Conclusion and Outlook}

The \MatOct~toolbox described in this work represents the first step of 
a~multi-purpose package that will include performance optimized discretizations 
for a~range of standard problems, first, from the CFD, and then, conceivably, from 
other application areas. 
Several important features will be added in the upcoming parts of this paper and 
in the new releases of the toolbox: 
slope limiters based on the nodal slope limiting procedure proposed 
in~\cite{Aizinger2011, Kuzmin2010}, convection terms, nonlinear advection operators, and
higher order time solvers.
Furthermore, the object orientation capabilites provided by \Matlab's \code{classdef},
which are expected to be supported in future \Octave~versions, will be exploited
to provide a more powerful and comfortable user interface.

\subsection*{\bf Acknowledgments}
The work of B.~Reuter was supported by the German Research Foundation (DFG) under grant AI 117/1-1.

\section*{Index of notation}
\noindent
\begin{footnotesize}
\begin{tabularx}{\linewidth}{@{}lX@{}}\toprule
\textbf{Symbol}    & \textbf{Definition}\\\midrule
$\avg{\,\cdot\,}$,\;  $\jump{\cdot}$   & average and jump on an~edge\\
$\diag(\vecc{A},\vecc{B})$ & $\coloneqq \begin{bmatrix}\vecc{A}&\quad\\\quad&\vecc{B}\end{bmatrix}$, block-diagonal matrix with blocks~$\vecc{A}$, $\vecc{B}$\\
$\card{\mathcal{M}}$ & cardinality of a~set~$\mathcal{M}$\\
$\vec{a}\cdot\vec{b}$& $\coloneqq \sum_{m=1}^2a_mb_m$, Euclidean scalar product in~$\IR^2$\\
$\circ$             & composition of functions or Hadamard product\\
$\otimes$           & Kronecker product\\
$\vec{a}_{ki}$, $\hat{\vec{a}}_i$ & $i$th vertex of the physical triangle~$T_k$, $i$th vertex of the reference triangle~$\hat{T}$\\
$c$                 & concentration (scalar-valued unknown)\\
$d$                 & diffusion coefficient\\
$\eta$              & penalty parameter\\
$\delta_{m\in\mathcal{M}}$ & $\coloneqq \{1~\text{if}~m\in\mathcal{M},~0~\text{if}~m\notin\mathcal{M}\}$, Kronecker delta\\
$E_{kn}$, $\hat{E}_n$          & $n$th edge of the physical triangle~$T_k$, 
$n$th edge of the reference triangle~$\hat{T}$\\
$\setV$,\; $\setE$,\; $\setT$   & sets of vertices, edges, and triangles\\
$\setE_{\mathrm{D}}$,\; $\setE_\mathrm{N}$  & set of boundary edges, $\setE_{\partial\Omega}= \setE_{\mathrm{D}}\cup\setE_\mathrm{N}$\\
$\setE_{\Omega}$,\; $\setE_{\partial\Omega}$ & set of interior edges, set of boundary edges\\
$\vec{F}_k$         & affine mapping from $\hat{T}$ to $T_k$\\
$h$                 & mesh fineness\\
$h_T$               & $\coloneqq \diam(h)$, diameter of triangle~$T\in\setT_h$\\
$J$                 & $\coloneqq (0,t_\mathrm{end})$, open time interval\\
$K$                 & $\coloneqq \card{\setT_h}$, number of triangles\\
$\vec{\nu}_{T}$     & unit normal on~$\partial T$ pointing outward of~$T$\\
$\vec{\nu}_k$       & $\coloneqq\vec{\nu}_{T_k}$\\
$N$                 & $\coloneqq (p+1)(p+2)/2$, number of local degrees of freedom\\
$\omega_r$          & quadrature weight associated with~$\hat{\vec{q}}_r$\\
$\Omega$,\; $\partial\Omega$          & spatial domain in two dimensions, boundary of $\Omega$\\
$\partial\Omega_\mathrm{D}$,\; $\partial\Omega_\mathrm{N}$& Dirichlet and Neumann boundaries,  $\partial\Omega = \partial\Omega_\mathrm{D}\cup\partial\Omega_\mathrm{N}$\\
$p$                 & $= (\sqrt{8N+1}-3)/2$, polynomial degree\\
$\vphi_{ki}$,\;  $\hat{\vphi}_i$      & $i$th hierarchical basis function on~$T_k$, $i$th hierarchical basis function on~$\hat{T}$\\
$\IP_p(T)$          & space of polynomials of degree at most~$p$\\
$\IP_p(\setT_h)$    & $\coloneqq \{ w_h:\overline{\Omega}\rightarrow \IR\,;\forall T\in\setT_h,\, {w_h}|_T\in\IP_p(T)\}$\\
$\hat{\vec{q}}_r$   & $r$th quadrature point in~$\hat{T}$\\
$R$                 & number of quadrature points\\
$\IR^+$,\; $\IR_0^+$  & set of (strictly) positive real numbers, set of nonnegative real numbers\\
$t$                 & time variable\\
$t^n$               & $n$th time level\\
$t_\mathrm{end}$    & end time\\
$\mapEE_{n^-n^+}$   & mapping from $\hat{E}_{n^-}$ to $\hat{E}_{n^+}$\\
$\Delta t^n$        & $\coloneqq t^{n+1} - t^n$, time step size\\
%
$T_k$,\; $\partial T_k$               & $k$th physical triangle, boundary of~$T_k$\\
$\hat{T}$           & bi-unit reference triangle\\
$\vec{x}$           & $=\transpose{[x^1,x^2]}$, space variable in the physical domain~$\Omega$\\
$\hat{\vec{x}}$     & $=\transpose{[\hat{x}^1, \hat{x}^2]}$, space variable in the reference triangle~$\hat{T}$\\
$\vec{z}$           & diffusion mass flux (vector-valued unknown)\\
\bottomrule
\end{tabularx}
\end{footnotesize}

\bibliography{FESTUNG}
\bibliographystyle{elsarticle-num}

\end{document}